\documentclass[11pt]{book}
\usepackage{amsfonts,amssymb,amsmath,epsf,graphics,graphpap}

\usepackage[english]{babel}
\newcommand{\halmos}{\rule{1ex}{1.4ex}}

\newtheorem{ittheorem}{Theorem}
\newtheorem{itlemma}{Lemma}
\newtheorem{itproposition}{Proposition}
\newtheorem{itdefinition}{Definition}
\newtheorem{itremark}{Remark}

\newenvironment{theorem}{\addtocounter{equation}{1}
\begin{ittheorem}}{\end{ittheorem}}

\newenvironment{lemma}{\addtocounter{equation}{1}
\begin{itlemma}}{\end{itlemma}}

\newenvironment{proposition}{\addtocounter{equation}{1}
\begin{itproposition}}{\end{itproposition}}

\newenvironment{definition}{\addtocounter{equation}{1}
\begin{itdefinition}}{\end{itdefinition}}

\newenvironment{remark}{\addtocounter{equation}{1}
\begin{itremark}}{\end{itremark}}

\newenvironment{corollary}{\addtocounter{equation}{1}
\begin{itcorollary}}{\end{itcorollary}}

\newenvironment{proof}{\noindent {\em Proof}.\,\,\,}
{\hspace*{\fill}$\halmos$\medskip}

\newcommand{\beq}{\begin{eqnarray}}
\newcommand{\eeq}{\end{eqnarray}}

\newcommand{\beqt}{\begin{eqnarray*}}
\newcommand{\eeqt}{\end{eqnarray*}}

\newcommand{\be}{\begin{equation}}
\newcommand{\ee}{\end{equation}}

\newcommand{\bl}{\begin{lemma}}
\newcommand{\el}{\end{lemma}}

\newcommand{\br}{\begin{remark}}
\newcommand{\er}{\end{remark}}

\newcommand{\bt}{\begin{theorem}}
\newcommand{\et}{\end{theorem}}

\newcommand{\bd}{\begin{definition}}
\newcommand{\ed}{\end{definition}}

\newcommand{\bp}{\begin{proposition}}
\newcommand{\ep}{\end{proposition}}

\newcommand{\bc}{\begin{corollary}}
\newcommand{\ec}{\end{corollary}}

\newcommand{\bpr}{\begin{proof}}
\newcommand{\epr}{\end{proof}}

\newcommand{\bi}{\begin{itemize}}
\newcommand{\ei}{\end{itemize}}

\newcommand{\ben}{\begin{enumerate}}
\newcommand{\een}{\end{enumerate}}

\newcommand{\Z}{\mathbb Z}
\newcommand{\R}{\mathbb R}
\newcommand{\N}{\mathbb N}

\newcommand{\E}{\mathbb E}

\newcommand{\s}{\ensuremath{\mathcal{S}}}

\newcommand{\om}{\ensuremath{\omega}}

\newcommand{\la}{\ensuremath{\Lambda}}
\newcommand{\si}{\ensuremath{\sigma}}

\begin{document}

\title{{\bf  Introduction to (generalized) Gibbs measures\footnote{ Lectures given at the {\em Semana de Mec\^anica Estat\'\i stica},
(Departamento de Matem\'atica, Universidade Federal de Minas Gerais, Belo Horizonte, February
5-9, 2007) and at the Instituto de Matem\'atica de UFRGS
(Universidade Federal do Rio Grande do Sul, Porto Alegre, Abril
2-6,  2007).}}}

\author{Arnaud Le Ny\footnote{Laboratoire de math\'ematiques, \'equipe de statistiques
 et mod\'elisation stochastique, b\^atiment 425, universit\'e de Paris-Sud,
 91405 Orsay Cedex, France. E-mail: arnaud.leny@math.u-psud.fr.}}

\maketitle



\normalsize

\chapter{Introduction}

These notes have been written to complete a mini-course ''Introduction to (generalized) Gibbs measures''
given at the universities UFMG ({\em Universidade Federal de Minas Gerais}, Belo Horizonte) and UFRGS
({\em Universidade Federal do Rio Grande do Sul}, Porto Alegre) during the first semester 2007.
 The main goal of the lectures was
 to describe  Gibbs and generalized Gibbs measures on
lattices at a rigorous mathematical level, as equilibrium states
of systems of a huge number of particles in interaction. In
particular, our main message is that although the historical
approach based on potentials has been rather successful from a
physical point of view,
 one has to insist on (almost sure) continuity properties of conditional probabilities to get a proper mathematical framework.

Gibbs measures are ''probably'' the central object of {\em
Equilibrium statistical mechanics}, a branch of probability theory
that takes its origin from  Boltzmann (\cite{Bo}, 1876) and Gibbs
(\cite{Gi}, 1902), who introduced a statistical approach to
thermodynamics that allows to deduce collective macroscopic
behaviors from individual microscopic information. Starting from
the observation that true physical systems with a very disordered
microscopic structure, like gases, ferromagnets like irons etc.,
could present a more ordered, non-fluctuating, macroscopic
behavior, they start to consider the microscopic components as
random variables and macroscopic equilibrium states as probability
measures that concentrate on the ''most probable'' states among
the possible ''configurations'' of the microscopic system, in a
sense consistent with the second laws of thermodynamics. Of
course, they did not use these modern probabilistic terms at that
time, and it is one of the tasks of mathematical statistical
mechanics to translate their intuitions in a more modern and
rigorous formalism.

 These ideas have been first introduced
and justified by Boltzmann in his introduction of statistical
entropy \cite{Bo2} and have been thereafter used by Gibbs as a
postulate to introduce his {\em microcanonical, canonical} and
{\em grand canonical ensembles} \cite{Gi}, providing three
different ways of describing equilibrium states, which would
nowadays be called ''probability measures'',  at the macroscopic
level. The main goal of modern mathematical statistical mechanics
is thus to describe rigorously
 these concepts in the standard framework of probability and
measure theory that has been developed during the century
following Boltzmann's ideas, pursuing two main goals: To describe
these {\em ensembles} as proper probability measures allowing the
modelization of phase transitions phenomena, and to interpret
them as equilibrium states in a probabilistic sense that would
incorporate ideas taken from the second law of thermodynamics.

For this purpose of describing  phase transitions phenomena,
 roughly seen as the possibility to get
 different macroscopic structures for the same microscopic
 interaction (e.g. gas versus liquid, positive or negative
 magnetization of iron, etc.),
we shall see that an infinite-volume formalism, which can be
loosely justified by the large number of microscopic components in
any macroscopic part of interacting systems, is required. For the
sake of simplicity, and because it already incorporates many of
the most interesting features of the theory, we shall  focus on
systems where the whole space is modelled by a discrete infinite
lattice (mainly $\mathbb{Z}^d$), with, attached at each site, a
microscopic element modelled by a finite value (e.g. $+1$ for a
positive ''microscopic magnetization'' in iron).

To describe equilibrium states and to model phase transitions
phenomena in  such a framework, we are led to construct
probability measures on an infinite product probability space in
an alternative way to the standard Kolmogorov's construction. This
alternative ''DLR'' construction, rigorously introduced in the
late sixties by Dobrushin \cite{DOB} and Lanford/Ruelle \cite{LR},
makes use of systems of compatible conditional probabilities with
respect to the outside of finite subsets, when the outside is
fixed in a {\em boundary condition}, to reach thereafter
infinite-volume quantities. This  {\em DLR approach} can also be
seen as an extension of the Markov chains formalism and to
describe Gibbs measures we shall focus on an extension of the
Markov property, {\em quasilocality}, closely related to
topological properties of conditional probability measures.

As we shall see, this approach allows  to model phase transitions
and the related critical phenomena. In that case,  a qualitative
change of the macroscopic system at a  ''critical point'' is
physically observed, together with a very chaotic critical
behavior. This criticality is physically interpreted as a highly
correlated system without any ''proper scale'', namely where a
physical quantity called {\em the correlation length} should
diverge, and such a system  should be thus reasonably
scaling-invariant. These considerations have led to the use of the
so-called {\em Renormalization group (RG) transformations}, which
appeared to be a very powerful tool in the theoretical physics of
critical phenomena \cite{RGT2,RGT}. It also gave rise to
ill-understood phenomena, the {\em RG pathologies}, detected in
the early seventies by Griffiths/Pearce \cite{GP} and Israel
\cite{Is2}, and interpreted a few decades later  by van Enter {\em
et al.} \cite{VEFS} as the manifestation of the occurrence of
non-Gibbsianness. This last observation was the starting point of
the {\em Dobrushin program of restoration of Gibbsianness},
launched by Dobrushin in 1995 in a talk in Renkum \cite{DOB1} and
consisting in two main goals: Firstly, to provide alternative
(weaker) definitions of Gibbs measures that would be stable under
the natural scaling transformations of the RG, and secondly to
restore the thermodynamics properties for these new notions in
order to still be  able to interpret them as equilibrium states.
This gave rise
to {\em generalized Gibbs measures}.\\

These notes are organized as follows: We introduce in Chapter 2
the necessary mathematical background, focusing on topological and
measurable properties of functions and measures on an infinite
product probability space; we also recall there important
properties of conditional expectations and introduce regular
versions of conditional probabilities to describe the DLR
construction of measures on infinite probability product spaces, mainly following \cite{F,Ge}.
We describe then the general structure of the set of DLR measures in the realm of convexity theory
and mention a few general consequences and examples at the end of the same chapter.
 We introduce Gibbs measures in the context
of quasilocality and describe the main features of the set of
Gibbs measures for a given interacting system in Chapter 3. The
interpretation of Gibbs and quasilocal measures as equilibrium
states is rigorously established in a general set-up in Chapter 4,
and we describe renormalization group pathologies and generalized
Gibbs measures in Chapter 5.

\chapter{Topology and measures on
product spaces}

\section{Configuration space: set-up and notations}
\subsection{Lattices}
For the sake of simplicity,  which could also be loosely justified
by the very discrete nature of physics, the physical space will be
modelled by a lattice $S$, which in our examples will mostly be
the $d$-dimensional regular lattice $\Z^d$. It is endowed with a
canonical distance $d$ and its elements, called {\em sites}, will
be designed by Latin letters $i,j,x,y$, etc. A pair of sites
$\{i,j\}$ such that $d(i,j)=1$ will be called {\em nearest
neighbor} (n.n.) and denoted by $\langle i j \rangle$. Finite
subsets of the lattice $S$ will play an important role for us and
will be generically denoted by capital Greek letters $\Lambda,
\Lambda', \Delta$, etc. We denote  the set of these finite subsets
of $S$ by

$$ \mathcal{S}=\Big\{\Lambda \in \mathcal{S},\; | \Lambda |< \infty \Big\} = \Big\{ \Lambda \subset \subset S \Big\}
$$ where $| \Lambda |$
denotes the cardinality of $\Lambda$ and $\subset \subset$ means
inclusion of a finite set in a bigger set. This notation $\mid
\cdot \mid$ will be used for many different purposes without
giving  its exact meaning when it is  obvious. It will be moreover
mostly sufficient to work with increasing sequences
$(\Lambda_n)_{n \in \N}$ of {\em cubes},  defined e.g. when the
lattice is $\mathbb{Z}^d$ by $\Lambda_n= [-n,n]^d \cap S$ for all
$n \in \mathbb{N}$.
\subsection{Single-spin state spaces}
To each (microscopic) site $i$ of the lattice we attach the same finite\footnote{This theory  also holds, modulo a few adaptations, for more general measurable spaces, compact \cite{Ge} or even non-compact \cite{DOB0,Lebo}, but the simpler finite case  already gets the main features of the theory.}
measurable space $(E,\mathcal{E},\rho_0)$, of cardinality $e:=| E |$. The {\em a
priori measure} $\rho_0$ will then be chosen to be the normalized
uniform counting measure on the $\sigma$-algebra
$\mathcal{E}=\mathcal{P}(E)$, formally defined in terms of Dirac
measures by $\delta_0=\frac{1}{e} \; \sum_{q \in E} \; \delta_q$.

In our guiding example, the {\em Ising model of ferromagnetism} \cite{KO,Peie}, this set is $E=\{-1,+1\}$, but other models might be considered. At each site $i$ of
the lattice will be thus attached a random variable $\sigma_i \in
E$, called $spin$ to keep in mind this seminal Ising
model.

This finite measurable space E will be called the {\em single-spin
state space} and will be endowed with the discrete topology, for
which the singleton sets are open, so that the open sets are all
the subsets of $E$.
\subsection{Configuration space}
The microscopic states are then represented by the collections of
random variables $\sigma=(\sigma_i)_{i \in S}$, living  in the
infinite product space $(\Omega, \mathcal{F}, \rho)=(E^S,
\mathcal{E}^{\otimes S},\rho_0^{\otimes S})$ called {\em the
configuration space}. (Infinite-volume) {\em configurations} will be
denoted by Greek letters $\sigma,\omega$, etc.

 For any $\Lambda
\in \mathcal{S}$, the finite product space
$\Omega_\Lambda=E^\Lambda$  comes with a  finite collection of the
random variables $\sigma_i$ for the sites $i \in \Lambda$ and for
any $\sigma \in \Omega$, one denotes by
$\sigma_\Lambda=(\sigma_i)_{i \in \Lambda}$ this configuration at
finite volume $\Lambda$. We also define concatenated
configurations at infinite-volume by prescribing values on
partitions of $S$, writing e.g. $\sigma_\Lambda
\omega_{\Lambda^c}$ for the configuration which agrees with a
configuration $\sigma$ in $\Lambda$ and with another configuration
$\omega$ outside $\Lambda$.
\section{Measurable properties of the configuration space}
The {\em product $\sigma$-algebra}
$\mathcal{F}=\mathcal{E}^{\otimes S}$ is the smallest
$\sigma$-algebra generated by the set of cylinders
$C_{\sigma_\Delta} = \big\{ \omega \in \Omega :
\omega_\Delta=\sigma_\Delta \big\}$, when $\sigma_\Delta$ runs
over $\Omega_\Delta$ and $\Delta$ runs over $\mathcal{S}$. We also
write $\mathcal{C}=\big\{ (C_{\sigma_\Delta}), \sigma_\Delta \in
\Omega_\Delta, \Delta \in \mathcal{S} \big\}$ and
$\mathcal{C}_\Lambda=\big\{ (C_{\sigma_\Delta}), \sigma_\Delta \in
\Omega_\Delta, \Delta \subset \subset \Lambda \big\}$ for the
family of cylinders restricted to any sub-lattice $\Lambda \subset
S$, not necessarily finite.  Alternatively, one defines for all
sites $i$ of the lattice, the canonical projection $\Pi_i : \Omega
\longrightarrow E$ defined for all $\omega \in \Omega$ by
$\pi_i(\omega)=\omega_i$, and  denotes by $\Pi_\Lambda$ the
canonical projection from $\Omega$ to $\Omega_\Lambda$ for all
$\Lambda \in \mathcal{S}$, defined for $\omega \in \Omega$ by
$\Pi_\Lambda(\omega)=\omega_\Lambda:=(\omega_i)_{i \in \Lambda}$.
Then, using the following rewriting of the cylinders,
$$C_{\sigma_\Lambda}=\Pi_\Lambda^{-1} (\{\sigma_\Lambda\}), \; \forall \sigma \in \Omega$$
one gets that $\mathcal{F}$ is also  the
smallest $\sigma$-algebra that makes the projections measurable.

The macroscopic states will be represented by {\em random fields},
i.e. probability measures on $(\Omega,\mathcal{F})$, whose set
will be denoted by $\mathcal{M}_1^+(\Omega,\mathcal{F})$, or more
briefly $\mathcal{M}_1^+(\Omega)$. The simplest one is the a
priori product measure $\rho=\rho_0^{\otimes S}$ defined as the
product of $\rho_0$ on the cylinders and extended to the whole
lattice  by virtue of the Kolmogorov's extension theorem
\cite{Bill,fell}, recalled later in this chapter. This particular
random field models the equilibrium state of a non-interacting
particle system, for which the spins are independent random
variables.

In order to mathematically describe microscopic and macroscopic
behaviors, one would like to distinguish {\em local and non-local}
events. The local ones are the elements of a sub-$\sigma$-algebra
$\mathcal{F}_\Lambda$ for a finite $\Lambda \in \mathcal{S}$,
where  $\mathcal{F}_\Lambda$ is the $\sigma$-algebra generated by
the finite cylinders $\mathcal{C}_\Lambda$ defined above. A function $f:\omega \longrightarrow \R$ is said to be
$\mathcal{F}_\Lambda$-measurable if and only if (iff) ''it depends
only on the spins in $\Lambda$'':
$$f \in \mathcal{F}_\Lambda \; \Longleftrightarrow \; \big (\omega_\Lambda=\sigma_\Lambda \Longrightarrow
f(\omega)=f(\sigma) \big ).
$$

\begin{definition}[Local functions]
A function $f:\Omega \longrightarrow \mathbb{R}$ is said to be
{\em local} if it is $\mathcal{F}_\Lambda$-measurable for some
$\Lambda \in \mathcal{S}$. The set of local functions will be
denoted by $\mathcal{F}_{\rm{loc}}$.
\end{definition}

We shall use the same notation $f \in \mathcal{F}$
for the measurability w.r.t a $\sigma$-algebra or $f \in \mathcal{H}$
for the membership in a space $\mathcal{H}$ of functions.

Another important sub-$\sigma$-algebra concerns macroscopic
non-local events. It is the so-called $\sigma$-{\em algebra at
infinity}, of {\em tail} or {\em asymptotic events}
$\sigma$-algebra, formally defined by
$$
\mathcal{F}_{\infty} = \bigcap_{\Lambda \in \mathcal{S}}
\mathcal{F}_{\Lambda^c}.
$$
Equivalently, it is the $\sigma$-algebra (countably) generated by the tail cylinders $\mathcal{C}_\infty:= \cap_{\Lambda \in \mathcal{S}} \; \mathcal{C}_{\Lambda^c}$.
It consists of events that do not depend on what happens in
microscopic subsets of the systems; they are typically defined by
some limiting procedure. In our description of the Ising model, we
shall encounter for example the tail events $B_m$, defined, for $m \in
[-1,+1]$, by
\begin{equation}\label{tailm}
B_m =\Big\{\omega: \;\lim_{n\to\infty}\; \frac{1}{|\Lambda_n
|}\sum_{i \in \Lambda_n}\omega_i=\; m \Big\}
\end{equation}
that will help to distinguish  the  physical phases of the system.
Similarly, a function $g$ is $\mathcal{F}_\infty$-measurable ($g \in \mathcal{F}_\infty$) if it
does not depend on the spins in any finite region, i.e. iff
$$
\exists \Lambda \in \mathcal{S} \; \rm{s.t.} \;
\sigma_{\Lambda^c}=\omega_{\Lambda^c} \; \Longrightarrow
g(\omega)=g(\sigma).
$$

These functions will be important later on to characterize
macroscopic quantities and to detect non-Gibbsianness. They are
also generally defined by some limiting procedure, the following
function being e.g. tail-measurable:

\begin{displaymath}
\forall \omega \in \Omega, \; g(\omega)=\left\{
\begin{array}{lll}\; \lim_{n\to\infty}\; \frac{1}{\mid \Lambda_n \mid} \sum_{i \in \Lambda_n} \omega_i \; & &\rm{if \; the \; limit \; exists.}\\
\\
\; {\rm anything} \;  \; & &\textrm{otherwise.}\\
\end{array} \right.
\end{displaymath}

\medskip

Similar tail $\sigma$-algebras are also used in ergodic theory or
in classical probability theory, in some $0$-$1$-laws for example
\cite{BR,Will}. To connect with these fields, we introduce here
the basic notion of translation-invariance, which  will also be
important for physical interpretations later on. For simplicity,
we introduce this notion on the lattice $S=\mathbb{Z}^d$ but it
could be easily extended to other lattices.

First one defines {\em translations} on the lattice as a
family of invertible transformations $(\tau_x)_{x \in
\mathbb{Z}^d}$ indexed by the sites of the lattice and defined for
all $x \in \mathbb{Z}^d$ by $$\tau_x :  y \longmapsto \tau_x y = y + x \in \mathbb{Z}^d$$
where additions and subtractions on the lattice are standard.
 They induce translations on $\Omega$: The translate by $x$
of $\omega \in \Omega$ is the configuration
$\tau_x \omega$ defined for all $i \in S$ by
$$
(\tau_x \omega)_i=\omega_{\tau_{-x} i}=\omega_{i-x}.
$$
It also extends naturally to measurable sets (our ''events''),
measurable functions and  measures. In particular, the set of {\em
translation-invariant probability measures} on
$(\Omega,\mathcal{F})$ is denoted by
$\mathcal{M}_{1,\;\rm{inv}}^+(\Omega)$ and the $\sigma$-algebra
generated by the translation-invariant functions is the {\em
translation-invariant} $\sigma$-{\em algebra} denoted by
$\mathcal{F}_{\rm{inv}}$.\\

 Let us briefly leave the field of lattices to consider another
 framework that links our approach to exchangeability in
the context of the so-called {\em mean-field models}. When
$S=\mathbb{N}$, one can define first a group $I_n$ of {\em
permutations at finite volume $n$}, that are bijections that leave
invariant the sites $i >n$, and define its union $I=\cup_{n \in
\mathbb{N}} \; I_n$ to be the {\em group of all permutations of
finitely many coordinates}. The $I$-invariant probability measures
on $(\Omega,\mathcal{F})$ form the set $\mathcal{M}_I$ of {\em
exchangeable probability measures}. For  $n \in \mathbb{N}$, the
$\sigma$-algebra of the events invariant under permutations of
order $n$ is defined to be
$$
\mathcal{I}_n= \big\{ A \in I_n : \pi^{-1} A = A ,\; \forall \pi
\in I_n \big\}
$$
and its intersection is the $\sigma$-algebra of {\em symmetric or
permutation-invariant} events
\be \label{exch}
\mathcal{I}=\bigcap_{n \in
\mathbb{N}} \mathcal{I}_n.
\ee

\section{Topological properties of the configuration space}
\subsection{Product topology}
As we shall see, the notion of Gibbs measures is based on the
interplay between topology and measure theory, and to relate these notions
we need to introduce a topology $\mathcal{T}$ whose Borel $\sigma$-algebra coincides with $\mathcal{F}$.
 The latter and  $\mathcal{T}$ are then said to be  {\em
compatible} in the sense
that  both open sets and continuous functions are
then  measurable. Thus, the topology $\mathcal{T}$ is endowed with
the same generators as those of $\mathcal{F}$ and $\mathcal{T}$ is the
smallest topology on $\Omega$ containing the cylinders  or making
the projections  continuous. To do so, we consider on the whole
configuration space $(\Omega,\mathcal{F})$ the product topology
$\mathcal{T}= \mathcal{T}_0^{\otimes_S}$ of the discrete
topology $\mathcal{T}_0$ on  $E$. Endowed with these
topological and measurable structures, our configuration space has
the  following nice properties:

\begin{theorem}\cite{F,Ge}
 The topological space  $(\Omega,\mathcal{T})$ is compact, its Borel
$\sigma$-algebra coincides with the product $\sigma$-algebra $\mathcal{F}$, and the measurable space $(\Omega,\mathcal{F})$ is
a Polish space, i.e. metrizable, separable and complete.
\end{theorem}

Compactness follows from Tychonov's theorem and will be helpful in
proving existence results and to simplify the topological
characterizations of Gibbs measures.

As a metric, one can choose $\delta : \Omega \times \Omega
\longrightarrow \mathbb{R}^+$, defined for all $\omega,\sigma \in
\Omega$ by
\begin{displaymath}
\delta (\omega,\sigma)=\sum_{i \in S} 2^{-n(i)}
\mathbf{1}_{\{\omega_i \neq \sigma_i\}}
\end{displaymath}
where $ n \colon S \longrightarrow \mathbb N$ is any bijection
assumed to be fixed and known. With this topology, open sets are finite unions of
cylinders and in  particular, a typical neighborhood of $\omega \in
\Omega$ is given by a cylinder  for $\Lambda
\in \mathcal{S}$ denoted in this context by
\begin{displaymath}
\mathcal{N}_{\Lambda}(\omega)=\big \{\sigma \in \Omega :
 \sigma_{\Lambda}=\omega_{\Lambda} , \sigma_{\Lambda^c} \; \textrm{arbitrary} \; \big \}.
\end{displaymath}
Similarly, when $S=\mathbb{Z}^d$, a basis of neighborhoods of a
configuration $\omega \in \Omega$ is given by the family of
cylinders $(\mathcal{N}_{\Lambda_n}(\omega))_{n \in \mathbb{N}}$,
for a sequence of cubes $(\Lambda_n)_{n \in \mathbb{N}}$. Thus,
two configurations are closed in this topology if they coincide
over large finite regions, and the larger the region is, the
closer they are\footnote{This topological framework is standard
also when one consider Cantor sets and dyadic expansions of
reals.}. Moreover, the set of asymptotic events is dense for this
topology, because they are insensitive to changes in finite
regions. In particular, the set of configurations that are
asymptotically constant is  a countable and dense subset, leading
thus to separability of the product topology by compatibility of
the latter with the measurable structure.

\subsection{Quasilocality for functions}

This nice topological setting allows us to provide different
equivalent characterizations of microscopic quantities. Firstly,
we find it natural to say that a microscopic function  $f$ on
$\Omega$ is arbitrarily "close" to functions which depend on
finitely many coordinates, i.e.  local functions. This leads
to the important concept of a {\em quasilocal function}:

\begin{definition}\label{quasilocalfcts}
A function $f:\Omega \longrightarrow \mathbb{R}$ is said to be
{\em quasilocal} if it can be uniformly approximated by local
functions, i.e. if   for each $\epsilon>0$, there exists  $f_\epsilon \in \mathcal{F}_{{\rm loc}}$ s.t.
$$
\sup_{\omega \in \Omega}  \Big| f(\omega) - f_\epsilon(\omega)
\Big| \; < \; \epsilon.
$$
\end{definition}

We denote by $\mathcal{F}_{\rm{qloc}}$ the set of quasilocal
functions. It is the uniform closure of $\mathcal{F}_{\rm{loc}}$
in the sup-norm, and by compactness is automatically bounded.
Moreover, due to the Polish and compact structure of $\Omega$, one
can use sequences and makes coincide continuity and uniform
continuity. Quasilocal functions are continuous while asymptotic
tail-measurable functions are discontinuous. Then, using the
metric $\delta$ or the basis of neighborhoods described above, it
is a simple exercise to prove that quasilocal functions are in
fact the (uniformly) continuous functions on $\Omega$, and we use
it in the next lemma to give alternative definitions of
quasilocality. When we do not use sequences, we shall deal with
the following convergence:

\begin{definition}[Convergence along a net  directed
  by inclusion]\label{net}
$$\lim_{\Lambda \uparrow \mathcal{S}} F(\Lambda) = a$$ means convergence
of a set-function $F : \mathcal{S}  \longrightarrow \mathbb{R}$ along a set
$\mathcal{S}$ directed by inclusion:
\begin{displaymath}
\forall \epsilon > 0 , \exists\  K_{\epsilon} \in \mathcal{S} \;
\textrm{s.t.} \; \mathcal{S} \ni  \Lambda  \supset K_{\epsilon}\
\Longrightarrow  \Big| F(\Lambda) - a \Big| \leq \epsilon.
\end{displaymath}
\end{definition}
\begin{lemma}\cite{F,Ge} A function $f:\Omega \longrightarrow \mathbb{R}$ is quasilocal iff
one of the following holds:
\begin{itemize}
\item{\bf Continuity:} It is continuous at every $\omega \in
\Omega$, i.e. $\forall \omega \in \Omega,\forall \epsilon
>0$, $\exists n \in \mathbb{N}$ s.t.
$$
\sup_{\sigma \in \Omega} \Big| f(\omega_{\Lambda_n}
\sigma_{\Lambda_n^c}) - f(\omega) \Big| < \epsilon.
$$
\item{\bf Uniform limit of local functions:} There exists  $(f_n)_{n \in \mathbb{N}}$ s.t. $\forall n \in \mathbb{N}$, $f_n \in \mathcal{F}_{\Lambda_n}$ and
$$
\lim_{n \to \infty} \sup_{\omega \in \Omega} \Big| f_n(\omega) -
f(\omega) \Big| = 0.
$$
\item{\bf Sequential uniform continuity:} For each $\epsilon >0$,
there exists $n \in \mathbb{N}$ s.t.
$$
 \sup_{\sigma, \omega \in \Omega} \Big|
f(\omega_{\Lambda_n} \sigma_{\Lambda_n^c}) - f(\omega) \Big|<
\epsilon.
$$
\item{{\bf Uniform continuity:}} $$\lim_{\Lambda \uparrow
\mathcal{S}} \sup_{\omega , \sigma \in \Omega ,\omega_{\Lambda} =
\sigma_{\Lambda}} \Big| f(\omega)-f(\sigma) \Big| = 0.$$
\end{itemize}
\end{lemma}

An important consequence of this lemma is that a non-constant
tail-measurable function can never be quasilocal. For example, let
us consider the event
\begin{displaymath}
B_0=\Big\{\omega: \lim_{n \to \infty} \frac{1}{\mid \Lambda_n
\mid}\sum_{i \in \Lambda_n}\omega_i=0\Big\}.
\end{displaymath}
The indicator function $f$ of this event is tail-measurable, non-constant and non-quasilocal. Take for example
$\Omega=\{0,1\}^{\mathbb Z}$. The configuration
$\omega=\mathbf{0}$, null everywhere, belongs to $B_0$ and
$f(\mathbf{0})=1$. Let $\mathcal{N}$ be a neighborhood of this
null configuration, and choose it to be
$\mathcal{N}_{\Lambda_n}(\mathbf{0})$ for some $n>0$. There exists
then $\sigma \in \mathcal{N}_{\Lambda_n}(\mathbf{0})$ such that
$\sigma_{{\Lambda_n^c}}=\mathbf{1}_{\Lambda_{n}^{c}}$, where
$\mathbf{1} \in \Omega$ is the configuration which value is $1$
everywhere. For this configuration, $\lim_{n \to \infty}
\frac{1}{\mid \Lambda_n
  \mid}\sum_{i \in \Lambda_n}\sigma_i =1$ and thus $f(\sigma)=0$:
This proves that $f$ is discontinuous and thus non-quasilocal.
Non-quasilocal functions will be important to detect non-Gibbsian
measures in Chapter 5.

\subsection{Weak
convergence of probability measures}

We have already introduced the space $\mathcal{M}_1^+(\Omega,\mathcal{F})$ of probability measures
on $(\Omega,\mathcal{F})$ that
 represents the macroscopic possible
states of our systems. Before introducing different ways of
constructing such measures on our infinite product spaces, we need
 a proper notion for the convergence of probability
measures, i.e.  to introduce a topology on
$\mathcal{M}_1^+(\Omega)$. For any $\mu \in
\mathcal{M}_1^+(\Omega)$ and $f \in \mathcal{F}_{\rm{qloc}}$, we
write $\mu[f]=\int f d\mu$ for the expectation of $f$ under $\mu$.
A strong way to do so is to consider the topology inherited from
the so-called {\em total-variation norm} but this is indeed   too
strong a notion of convergence due to our willing of describing
''non-chaotic'' equilibrium states: Physically, this convergence
means that expected values  converges, {\em uniformly} for all
bounded or continuous observables, i.e. microscopic in our point
of view, and this occurs rarely in physical situations. We shall
thus require a topology whose convergence mainly concerns
non-uniform expectations of microscopic variables. This is the
famous {\em weak convergence of probability measures}, which is
indeed weaker than most  ways of convergence, see
\cite{Bill,fell,Will}.

\begin{definition}[Weak convergence]
A sequence $(\mu_n)_{n \in \mathbb{N}}$ in
$\mathcal{M}_1^+(\Omega)$ is said to {\em converge} weakly to $\mu
\in \mathcal{M}_1^+(\Omega)$ if  expectations of continuous
functions converge:
$$
\mu_n \; \stackrel{W}{\longrightarrow} \; \mu \;
\Longleftrightarrow \; \lim_{n \to \infty} \; \mu_n[f] \;= \;
\mu[f], \; \forall f \in \mathcal{F}_{\rm{qloc}}.
$$
\end{definition}

This convergence gives no information on the convergence of the
expectations of discontinuous (macroscopic, asymptotic)
quantities. This will be important for our purpose of modelling
phase transitions phenomena by working at finite but larger and
larger sets through some infinite-volume limit. By definition, the
set of local functions is dense in $\mathcal{F}_{\rm{qloc}}$, so
it is enough to test this convergence on $\mathcal{F}_{\rm{loc}}$
or on cylinders.

To describe, at the end of the chapter, the general convex structure of the set of Gibbs measures in case of phase transitions, we shall also need to deal with probability measures on spaces of probability measures, and we first endow such spaces with a canonical measurable structure. For any subset of probability measures $\mathcal{M} \subset  \mathcal{M}_1^+(\Omega,\mathcal{F})$, the natural way to do so is to evaluate any $\mu \in \mathcal{M}$ via the numbers $\big \{\mu(A), \; A \in \mathcal{F} \big \}$. One introduces then the {\em evaluation maps} on $\mathcal{M}$ defined for all $A \in \mathcal{F}$ by
\begin{equation}\label{eval}
e_A : \mathcal{M} \longrightarrow [0,1]; \mu \longmapsto e_A(\mu)=\mu(A).
\end{equation}

The {\em evaluation} $\sigma$-{\em algebra} $e(\mathcal{M})$
is then the smallest $\sigma$-algebra on $\mathcal{M}$ that makes measurable
these evaluation maps,
 or equivalently the $\sigma$-algebra generated by the sets $\{e_A \leq c\}$ for all $A \in \mathcal{F}, \; 0 \leq c \leq 1$.
 For any bounded measurable function $f \in \mathcal{F}$,
  the map $$e_f : \mathcal{M} \longrightarrow [0,1]; \mu \longmapsto e_f(\mu):=\mu[f]$$ is then $e(\mathcal{M})$-measurable.

\section{Probability theory on infinite product
spaces}
\subsection{Kolmogorov's consistency}

The standard way to construct probability measures on an infinite
product measurable space is to start from a {\em consistent system
of finite dimensional marginals}, following a terminology of
Kolmogorov \cite{Bill}:

\begin{definition}
A family $(\mu_\Lambda)_{\Lambda \in \mathcal{S}}$ of probability
measures on $(\Omega_\Lambda,\mathcal{F}_\Lambda)$ is said to be
{\em consistent} in the sense of Kolmogorov iff for all $\Lambda
\subset \Lambda'  \in \mathcal{S}$,
$$
\mu_{\Lambda}(A)=\mu_{\Lambda'} \Big(\big(
\Pi^{\Lambda'}_{\Lambda}\big)^{-1}(A) \Big), \; \forall A \in
\mathcal{F}_\Lambda
$$
where $\Pi_\Lambda^{\Lambda'}$ is the natural projection from
$\Omega_{\Lambda'}$ to $\Omega_\Lambda$.
\end{definition}

Given a consistent family of conditional probabilities, it is
possible to extend it, in our mild framework, to the whole
configuration space:
\begin{theorem}[Kolmogorov's extension theorem]
Let $(\mu_\Lambda)_{\Lambda \in \mathcal{S}}$ be a consistent
family of marginal distributions  on a {\em Polish} infinite-
product probability space $(\Omega,\mathcal{F})$. Then there
exists {\em a unique} probability measure $\mu \in
\mathcal{M}_1^+(\Omega)$ s.t. for all $\Lambda \in \mathcal{S}$,
$$
\forall A \in \mathcal{F}_\Lambda,\; \mu
\big(\Pi_\Lambda^{-1}(A)\big)=\mu(A)
$$
where $\Pi_\Lambda^{-1}(A)$ is the pre-image of $A$ by the
projection from $\Omega$ to $\Omega_\Lambda$, defined by
$$
\Pi_\Lambda^{-1}(A)=\Big\{ \sigma \in \Omega :
\Pi_\Lambda(\sigma)=\sigma_\Lambda \in A \Big\}.
$$
\end{theorem}

The main example of application of this theorem is the
construction of the a priori  product measure $\rho$ on
$(E,\mathcal{E})$. Consider the counting measure $\rho_0 \in
\mathcal{M}_1^+(\Omega)$ on the single-site state space and the
finite product measure $\rho_\Lambda=\rho_0^{\otimes \Lambda}$ on
any of the finite product probability spaces
$(\Omega_\Lambda,\mathcal{F}_\Lambda)$, defined for all $\Lambda
\in \mathcal{S}$ on the cylinders by  $$\rho_\Lambda
(\sigma_\Lambda) = \prod_{i \in \Lambda} \rho_0(i), \; \forall
\sigma_\Lambda \in \Omega_\Lambda$$ and extended on $
\mathcal{F}_\Lambda$ by requiring, for all $A \in
\mathcal{F}_\Lambda$,
$$
\rho_\Lambda(A)=\sum_{\sigma_\Lambda \in A} \rho_\Lambda
(\sigma_\Lambda).
$$

 The system $(\rho_\Lambda)_{\Lambda \in \mathcal{S}}$ is
trivially a consistent family of marginals, and our configuration
space being a Polish space, it  extends into a unique
probability measure
$\rho=\rho_0^{\otimes S}$.\\

Hence, we know how to build elements $\mu$ of the (convex) set
$\mathcal{M}_1^+(\Omega,\mathcal{F})$, which are interpreted  as a
macroscopical description of the physical phases of the systems in
our settings. The problem now is that we also want to model phase
transitions, i.e. to get different infinite-volume measures
corresponding to the same finite volume description, and to do so we
have to proceed differently and work with systems of conditional
probabilities consistent in a different sense than that of
Kolmogorov, based on successive conditionings w.r.t. decreasing
sub-$\sigma$-algebras. Let us first recall a few important
properties of conditional probabilities on infinite product
probability Polish spaces.
\subsection{Regular versions of
conditional probabilities}
\begin{definition}\label{condexpect}[Conditional expectation] Let $(\Omega,\mathcal{F})$ be a measurable space, $\mu \in
\mathcal{M}_1^+(\Omega,\mathcal{F})$, $\mathcal{G}$ a
sub-$\sigma$-algebra of $\mathcal{F}$  and $f \in \mathcal{F}$,
$\mu$-integrable. A conditional expectation of $f$ given
$\mathcal{G}$, w.r.t. $\mu$, is a function $\mathbb{E}_\mu[f \mid
\mathcal{G}] : \Omega \longrightarrow, \;
 \mathbb{R}; \omega \longmapsto \mathbb{E}_\mu[f \mid \mathcal{G}](\omega)$ such that
\begin{enumerate}
\item $\mathbb{E}_\mu[f \mid \mathcal{G}]$ is
$\mathcal{G}$-measurable. \item For any $g \in \mathcal{G}$ bounded,
$\ \ \ \ \ \ \int g \cdot \mathbb{E}_\mu[f \mid \mathcal{G}] \; d
\mu=  \int g
\cdot f d \mu$\\

and in particular $\ \ \ \ \ \ \; \; \ \ \ \  \int \mathbb{E}_\mu [f
\mid \mathcal{G}] \; d \mu = \int f d\mu.$
\end{enumerate}
\end{definition}

The existence of such functions is insured by the Radon-Nikod\'ym
theorem \cite{Dell, Will}. Nevertheless, a $\mu$-integral being involved in
point 2. of the definition above, such a conditional
expectation is not  unique, but  two
different {\em versions} of it can only differ at most on a set of
$\mu$-measure zero. Thus, Definition \ref{condexpect} does not
define a unique function, but measure-zero modifications are
however the only one possible: The conditional expectation
$\mathbb{E}_\mu[f \mid \mathcal{G}]$ is thus defined "$\mu$-a.s."\\

At this point, in the purpose of defining a probability via a
prescribed system of conditional probabilities w.r.t. the outside
of finite sets, one could get into troubles when trying without
care to give a sense to conditional probabilities w.r.t. a
sub-$\sigma$-algebra. In the same settings as the definition
above, the good candidate for such an almost-surely defined
conditional probability $\mu_{\mathcal{G}}(\cdot \mid \omega)$
would be defined to be, for all $A \in \mathcal{F}$,
$$
\mu_{\mathcal{G}}(A \mid \omega):=\E_\mu[\mathbf{1}_A \mid
\mathcal{G}](\omega), \; \; \mu\rm{-a.s.}(\omega).
$$
Indeed the following characterizing properties of a
probability measure are (a.s.) true:
\begin{itemize}
\item $\mu_\mathcal{G} (\Omega \mid \cdot)=1, \; \mu$-a.s. and
$\mu_\mathcal{G} (\emptyset \mid \cdot)=0, \;\mu$-a.s. \item For
all $A \in \mathcal{F}, \; 0 \leq \mu_\mathcal{G} (A \mid \cdot)
\leq 1,\; \mu$-a.s. \item For any countable collection $(A_i)_{i \in I}$ of
pairwise disjoints elements of $\mathcal{F}$,
$$
\mu_\mathcal{G} (\cup_i A_i \mid \cdot \; )=\sum_i \mu_\mathcal{G}
(A_i \mid \cdot \; ),\; \mu{\rm -a.s.}
$$
\end{itemize}
and one also has that for all $B \in \mathcal{G}$,
$\mu_\mathcal{G} [B
\mid \cdot]=\mathbf{1}_B(\cdot),\; \mu$-a.s. \\

The problem in this definition comes from the fact that the previous
properties are only valid almost surely, and that the sets of
measure zero that appear depend on the sets $A$ and $(A_i)_{i \in I}
\in \mathcal{F}$ considered; the later being uncountably many, we
cannot say that we have defined these conditional probabilities
$\mu$-almost everywhere\footnote{To see how to construct
counter-examples, consult e.g. \cite{Stoy}.}. What is needed to say
so is to get a unique
 set
of $\mu$-measure zero, independent of the sets $A$, outside which
the above properties are true. In such a case one says that there
exists {\em a regular version of the conditional probabilities} of
$\mu$ w.r.t. sub-$\sigma$-algebras of $\mathcal{F}$. More
precisely, this occurs when there exists a probability kernel (see
next subsection) $\mu_\mathcal{G}$ from $(\Omega,\mathcal{F})$ to
itself such that
$$
 \;
\mu{\rm -a.s.},\; \mu_\mathcal{G}[f \mid \cdot]=\mu[f \mid \mathcal{G}](\cdot), \forall f \in \mathcal{F} \; {\rm bounded}
$$
where the ''$\mu$-a.s.'' means that there exists a (mostly abstract)
measurable set of full $\mu$-measure $\Omega_\mu$ where the above
characterizing properties of a probability measure hold for all
$\omega \in \Omega_\mu$, independently of the measurable set $A \in
\mathcal{F}$. In our framework, this is hopefully granted:
\begin{theorem}\label{regular}\cite{Dell}
Any measure on a Polish probability space $(\Omega,\mathcal{F})$
admits a regular conditional probability w.r.t. any
sub-$\sigma$-algebra of $\mathcal{F}$.
\end{theorem}

 We also mention here two direct
consequences of Definition \ref{condexpect} which will be useful to
characterize measures in terms of systems of regular conditional
probabilities. Keeping the same settings, one has, $\mu$-almost
surely, for any  bounded  $\mathcal{G}$-measurable function $g$, any
bounded measurable function $f$ and any sub-$\sigma$-algebra
$\mathcal{G}' \subset \mathcal{G}$, $$ \mu_\mathcal{G}[g \cdot f
\mid \cdot ]= \E_\mu[g \cdot f \mid \mathcal{G}](\cdot)= g \cdot
\E_\mu[f \mid \mathcal{G}]=g \cdot \mu_\mathcal{G}[f \mid \cdot ]
$$ and \be \label{E2} \mu_{\mathcal{G}'}\big[\mu_\mathcal{G}[f \mid
\cdot] \cdot \big]= \E_\mu \big[ \E_\mu[f \mid \mathcal{G}] \mid
 \mathcal{G}'\big] (\cdot) = \E_\mu[f \mid \mathcal{G'}](\cdot)= \mu_{\mathcal{G}'}[f \mid \cdot].\ee

We recall now the useful concept of {\em probability kernel}  to
describe the alternative way of defining  probability measures on
infinite product probability spaces introduced in the late sixties
by Dobrushin, Lanford and Ruelle  to model phase transitions.

\begin{definition}\label{Probabker}
A {\em probability kernel} from a probability space
$(\Omega,\mathcal{F})$ to a probability space
$(\Omega',\mathcal{F}')$ is a
map $\gamma(\cdot \mid \cdot) : \mathcal{F}' \times \Omega \to [0,1]$ such that
\begin{itemize}
\item For all $\omega \in \Omega$, $\gamma( \cdot \mid \omega)$ is
a probability measure on $(\Omega',\mathcal{F}')$. \item For all
$A' \in \mathcal{F}'$, $\gamma(A'\mid \cdot)$ is
$\mathcal{F}$-measurable.
\end{itemize}
\end{definition}

The simplest example is  the map
$\gamma(A\mid\omega)=\mathbf{1}_A(\omega)$ defined  for any
probability space $(\Omega,\mathcal{F})$, any $A \in \mathcal{F}$,
any $\omega \in \Omega$. It is a probability kernel from
$(\Omega,\mathcal{F})$ into itself. More interesting examples
concern regular versions of conditional probabilities, Markov
transition kernels etc. We extend this notion in order  to
introduce the concept of specification and to prescribe
conditional probabilities of a measure to try to define it. To do
so, we  state a few definitions.

\begin{definition} Let $\gamma$ be a probability kernel from $(\Omega,\mathcal{F})$
to $(\Omega',\mathcal{F}')$. For any  function $f \in \mathcal{F}'$, we define $\gamma f \in \mathcal{F}$ to be the function defined for all
$\omega \in \Omega$ by
$$
\gamma f(\omega)= \int_{\Omega'} f(\sigma) \;  \gamma(d \sigma \mid
\omega).
$$
We also define for  any  $\mu \in \mathcal{M}_1^+(\Omega,\mathcal{F})$
the measure $\mu \gamma \in \mathcal{M}_1^+(\Omega',\mathcal{F}')$
by

$$
\forall A' \in \mathcal{F}',\;\; \mu \gamma(A')=\int_{\Omega} \gamma (A' \mid \omega) \; \mu (d
\omega).
$$
\end{definition}

A little bit more has to be required for a kernel to represent a
regular version of a conditional probability. In order to
illustrate the "double-conditioning" stability (\ref{E2}) of
conditional probabilities  we also introduce the notion of  {\em
product} (or composition) of kernels.

\begin{definition} Let $\gamma$ be a kernel from $(\Omega,\mathcal{F})$ to
$(\Omega',\mathcal{F}')$ and $\gamma'$ a kernel from
$(\Omega',\mathcal{F}')$ to $(\Omega'',\mathcal{F}'')$. Then the {\em product}
$\gamma \gamma'$ is the kernel from $(\Omega,\mathcal{F})$ to
$(\Omega'',\mathcal{F}'')$ s.t.
$$\forall A'' \in \mathcal{F}'', \forall \omega \in \Omega, \; \;
\gamma \gamma' (A'' \mid \omega) = \int_{\Omega'} \gamma'(A'' \mid
\sigma) \gamma(d \sigma \mid \omega).
$$
\end{definition}

We are now ready to  give the more formal

\begin{definition}[Regular version of conditional probability]
Let $(\Omega,\mathcal{F},\mu)$ be a probability space and
$\mathcal{G}$ a sub-$\sigma$-algebra of $\mathcal{F}$. A {\em
regular (version of) conditional probability} of $\mu$ given $\mathcal{G}$ is a
probability kernel $\mu_\mathcal{G}(\cdot \mid \cdot)$ from
$(\Omega,\mathcal{G})$ to $(\Omega,\mathcal{F})$ s.t.
$$
\mu {\rm-a.s.}, \; \; \; \mu_\mathcal{G} [f \mid \cdot] = \mathbb{E}_\mu[ f \mid
\mathcal{G}](\cdot), \; \forall f  \in \mathcal{F} \; {\rm and} \; \mu-{\rm integrable}
$$
\end{definition}

Using the action of a kernel to a measure and the definition of
the conditional expectation, it is also possible to characterize
it in a more closed form, which will lead soon to a consistency
condition different from the Kolmogorov one.

\begin{definition}[Regular conditional probability II]
In the same settings as above,  {\em regular conditional
probability} of $\mu$ given $\mathcal{G}$ is a probability kernel
$\mu_\mathcal{G} (\cdot \mid \cdot)$ from $(\Omega,\mathcal{F})$
to itself s.t. for all $f$ $\mathcal{F}$-measurable and
$\mu$-integrable,
\begin{enumerate}
\item $\mu_\mathcal{G}[f \mid \cdot]$ is $\mathcal{G}$-measurable.
\item $\mu$-a.s.,\;  $\mu_\mathcal{G} [g \cdot f \mid \cdot]=g
\cdot \mu_\mathcal{G} [f \mid \cdot]$, for each bounded $g \in
\mathcal{G}$ \item The kernel leaves invariant the probability
measure $\mu$: $\; \; \; \; \mu \mu_\mathcal{G} = \mu$.
\end{enumerate}
\end{definition}

This last definition, coupled to the fact that every measure on a
Polish space has regular conditional probabilities, enables also
to describe the double conditioning property in terms
of kernels and will give rise to the concept of specification.

\begin{definition}[System of regular conditional probabilities]
Let $(\Omega,\mathcal{F},\mu)$ be a probability space and
$(\mathcal{F}_i)_{i\in I}$ a family of sub-$\sigma$-algebras of
$\mathcal{F}$. A system of regular conditional probabilities of
$\mu$ given $(\mathcal{F}_i)_{i \in I}$ is a family of probability
kernels $\big(\mu_{\mathcal{F}_i}\big)_{i \in I}$ on
$(\Omega,\mathcal{F})$ s.t.
\begin{enumerate}
\item For each $i \in I$, $\mu_{\mathcal{F}_i}$ is a regular
conditional probability of $\mu$ given $\mathcal{F}_i$. \item If
$i,j \in I$ are such that $\mathcal{F}_i \subset \mathcal{F}_j$, then $\; \; \; \mu_{\mathcal{F}_i} \mu_{\mathcal{F}_j} =
\mathcal{\mu}_{\mathcal{F}_i}, \; \; \mu{\rm-a.s.}$
\end{enumerate}
\end{definition}

In statistical mechanics we go in the opposite direction: Starting
from a regular system of conditional probabilities, one wants to
reconstruct probability measures, and so one aims at removing the
"$\mu-$a.s." of the last definition, as one wants to obtain the
same conditional probabilities for different measures.

\subsection{DLR-consistency  and specifications}

Around 1970, Dobrushin \cite{DOB}, Lanford/Ruelle \cite{LR} have introduced a new way to construct probability measures
on infinite product probability spaces that does not immediately
yield uniqueness in the case of a Polish space, leaving the door
open to the modelling of phase transitions in mathematical
statistical mechanics. The key-point of their approach is to
replace a system of marginals consistent in the sense of Kolmogorov
 by a system of regular conditional probabilities
with respect to the outside of any finite set, giving rise to
 {\em finite volume versions of conditional probabilities
with prescribed boundary condition(s)}.

Let us consider now $\Lambda' \subset \Lambda \in \mathcal{S}$ and  the family of
sub-$\sigma$-algebra's $(\mathcal{F}_{\Lambda^c})_{\Lambda \in
\mathcal{S}}$, directed by inclusion in the sense that if
$\Lambda' \subset \Lambda \in \mathcal{S}$, one has

$$
\mathcal{F}_{\Lambda^c} \subset \mathcal{F}_{\Lambda'^c} \; {\rm
and} \; \bigcap_{\Lambda \in \mathcal{S}} \mathcal{F}_{\Lambda^c}
\; = \; \mathcal{F}_\infty.
$$

A system of regular conditional probabilities of $\mu \in
\mathcal{M}_1^+(\Omega)$ w.r.t. the mentioned filtration exists,
according to preceding section. To remove the ''$\mu$-a.s''
dependency and describe candidates to represent this system in the
case of an equilibrium state, Preston \cite{Pres2} has introduced
the concept of {\em specification}.

\begin{definition}[Specification]\label{specification}

 A specification is a family $\gamma
=(\gamma_{\Lambda})_{\Lambda \in
  \mathcal{S}}$ of probability kernels from $(\Omega,\mathcal{F})$
  into itself such that
\begin{enumerate}
\item For all $A \in \mathcal{F},\gamma_{\Lambda}(A|\cdot)$ is
$\mathcal{F}_{\Lambda^{c}}$-measurable. \item (Properness) For all $\omega \in \Omega,\; B \in
\mathcal{F}_{\Lambda^{c}}$ ,
$\gamma_{\Lambda}(B|\omega)=\mathbf{1}_{B}(\omega)$
\item (Consistency) \be\label{DLRcons} \Lambda' \subset \Lambda
\in \mathcal{S} \; \Longrightarrow \gamma_{\Lambda}
\gamma_{\Lambda'} = \gamma_{\Lambda}.\ee
\end{enumerate}
\end{definition}
We recall that $\gamma_{\Lambda} \gamma_{\Lambda'}$ is the map on
$\Omega \times \mathcal{F}$ defined by
\begin{displaymath}
\gamma_{\Lambda} \gamma_{\Lambda'}(A|\omega) =
\int_{\Omega}\gamma_{\Lambda'}(A|\omega')\; \gamma_{\Lambda}(d\omega'
| \omega).
\end{displaymath}

Specifications are thus the appropriate objects to describe
conditional probabilities; an important point is that they are
defined everywhere on $\Omega$, for the convenient reason that we
want to deal with objects defined everywhere (not $\mu$-a.s.),
and characterize $\mu$ afterwards. This will allow the description
of different measures for a single specification, that is  to model
{\em phase transitions}  in our settings.  We also emphasize that
for all $\sigma, \omega \in \Omega$, for all $\Lambda \in
\mathcal{S}, \gamma_{\Lambda}(\sigma|\omega)$ depends only on
$\sigma_{\Lambda}$ and $\omega_{\Lambda^{c}}$. For this reason, because only the components of $\omega$
outside $\Lambda$ (beyond the boundary) are involved,
$\omega$ is often called {\em boundary condition} in statistical
mechanics, and we shall use this term frequently in the next
sections. Moreover, Conditions 1. and 2. of the definition of a specification can be removed
by requiring  $\gamma_\Lambda$ to be kernels from $(\Omega_{\Lambda^c},\mathcal{F}_{\Lambda^c})$ to $(\Omega_\Lambda,\mathcal{F}_\Lambda)$.\\

In fact, our main goal in this course is precisely to describe
the set  the measures satisfying the consistency relation when
$\Lambda$ becomes the whole lattice.
\begin{definition}[DLR measures]
Let $\gamma$ be a specification on $(\Omega,\mathcal{F})$. The set
of DLR measures for $\gamma$ is the set
\begin{equation}\label{DLR1}
\mathcal{G}(\gamma)=\Big\{\mu \in
\mathcal{M}^{+}_{1}(\Omega,\mathcal{F}):
\forall \Lambda \in \mathcal{S}, \; \mu[A \mid
\mathcal{F}_{\Lambda^{c}}](\cdot)=\gamma_{\Lambda}(A \mid \cdot), \;
\mu \textrm{-a.s.},\; \forall A \in \mathcal{F}\Big\}
\end{equation}
of the probability measures \emph{consistent} with $\gamma$.
Equivalently,
\begin{equation}\label{DLR2}
\mu \in \mathcal{G}(\gamma) \; \Longleftrightarrow \; \mu \gamma_\Lambda = \mu, \; \forall \Lambda \in \mathcal{S}.
\end{equation}
A {\em DLR measure} is thus a measure specified by some
specification $\gamma$.
\end{definition}
This definition reminds one of the one of the Kolmogorov: Instead
of dealing with the family of marginals of the measure, we deal
with its system of conditional probabilities. It will be of
importance when one models equilibrium states. Indeed, for such a
DLR measure, the consistency relation implies that integrating out
with respect to boundary conditions typical for the
''equilibrium'' DLR state outside a finite volume does not change
the state in the finite volume. We shall be more precise about
equilibrium properties in Chapter 4. On a Polish space, the
Kolmogorov compatibility yields existence and uniqueness of the
consistent measure, whereas the set $\mathcal{G}(\gamma)$ could
have a very different structure, the latter being a very important
fact for our  purpose of modelling the phenomenon of phase
transitions. Indeed, in contrast to what occurs in Kolmogorov's
consistency theorem, here neither existence, nor uniqueness needs
to occur. Before describing more precisely various sets of
DLR-measures, we provide  a few examples describing these
different possible structures.

\subsection{Examples}

We begin  by two examples that illustrate the negative side of
this description, the possibility of non-existence of measures
specified by a specification. This will help us to extract the
topological properties required to build a satisfactory framework
describing Gibbs measures as equilibrium states of interacting
particle systems. Thereafter, we provide as an example of
uniqueness followed by an example of non-uniqueness interpreted as
the occurrence of a  {\em phase transition}, the standard Ising
model.
\begin{enumerate}
\item{{\bf One-dimensional random walk}}:\\

This analysis goes back to  Spitzer but our presentation is
inspired by \cite{Pr}. The single-spin state $E=\mathbb Z$ is not
compact and this is the reason for the non-existence of a DLR
measure. The lattice is the time, modelled by $S=\mathbb{Z}$. The
symmetric n.n. random walk on $\mathbb{Z}$, $Y=(Y_n)_{n \in
\mathbb{Z}}$, is then a random element of  the configuration space
$\Omega = \mathbb{Z}^{\mathbb{Z}}$ and we denote by $\mathbb{P}$
its
  law on $(\Omega,\mathcal{F})$,  canonically built using Kolmogorov's extension
  theorem.

  Let us  try to define a specification $\gamma$ with which $\mathbb{P}$ would be
  consistent. Using the Markov property for random walks,
  the candidate
  is given  by the kernel $\gamma_{\Lambda_n}$,
defined for all cube $\Lambda_n=[-n,n] \cap \Z$ and for all $\sigma , \omega \in \Omega$ by
\begin{displaymath}
\gamma_{\Lambda_n}(\sigma \mid \omega)=\mathbb{P}\big[Y_{\Lambda_n}=\sigma_{\Lambda_n}|Y_{-n-1}=\omega_{-n-1},Y_{n+1}=\omega_{n+1}\big]
\end{displaymath}
where the event in the conditioning is the cylinder $C_{\omega_{\{n-1,n+1\}}}$, of positive $\mathbb{P}$-measure.
It is straightforward to extend it to any finite $\Lambda \in
\mathcal{S}$ to get a family of proper kernels
$\gamma=(\gamma_\Lambda)_{\Lambda \in \mathcal{S}}$ that is indeed
a specification. Let us  assume that there exists $\mu \in \mathcal{G}(\gamma)$. Then we claim that $\mu$ cannot be a probability measure because
for all $k \in \mathbb{Z}$, and for all $\epsilon
>0$, $\mu[Y_0= k] < \epsilon$. The reason for this is that for
all $n \in \mathbb N,\; S_n=Y_n - Y_0$ follows a binomial law. It
is unbounded and thus for all $\epsilon
>0$ and $k \in \mathbb{Z}$ and $n$ big enough, $\mathbb{P}[S_n =
k] < \epsilon$. Using then the consistency relation $\mu \gamma_{\Lambda_n} = \mu$
to evaluate $\mu[Y_0= k]$ in terms of conditional probabilities of
$\mathbb{P}$, one gets this result of "escape of mass to
infinity": If $\mu \in \mathcal{G}(\gamma)$, then for all $k \in
\mathbb{Z}$, and for all $\epsilon >0$, $\mu[Y_0= k] \leq
\epsilon$, and thus $\mu$ cannot be a probability measure: $\mathcal{G}(\gamma)=\emptyset$.
\smallskip
\item{{\bf Totally random single-particle}}:\\

 This example has
been provided by Georgii \cite{Ge}. Consider the case of a lattice
gas, i.e. $E=\{0,1\}$, $S=\mathbb{Z}^d$, denote ${\bf 0}$ the
configuration null everywhere and for any  $a
\in S$ consider the configuration $\sigma^a$
characterizing a single particle localized at the site $a$ defined
for all $i \in S$ by $(\sigma^a)_i=1$ iff  $i=a$ (and zero otherwise). To model a single particle evolving totally at random in a lattice
gas in this DLR framework, introduce the following kernel,
proper by construction:

\begin{displaymath}
\forall \Lambda \in \mathcal{S},\; \forall \omega \in \Omega, \forall A \in \mathcal{F}\;,
\gamma_\Lambda(A \mid \omega)=\left\{
\begin{array}{lll}\; \frac{1}{\mid \Lambda \mid} \sum_{a \in \Lambda} \mathbf{1}_A(\sigma^a) \;
 \; \; \; \rm{if} \; \; \;  \omega_{\Lambda^c} = {\bf 0}_{\Lambda^c}\\
\\
\; \mathbf{1}_A({\bf 0}_\Lambda \omega_{\Lambda^c}) \;  \; \; \; \; \; \;\; \; \;  \; \rm{otherwise.}\\
\end{array} \right.
\end{displaymath}

Firstly, one can check that the corresponding $\gamma$ is a specification and  that for any sequence $(\mu_n)_n$ of probability measures on
$(\Omega_{\Lambda_n},\mathcal{F}_{\Lambda_n})$, the sequence of
probability measures $(\mu_n \gamma_{\Lambda_n})_n$ converges
weakly towards the Dirac measure $\delta_{\bf 0}$ on the null configuration, and thus that the latter is a good candidate to
be in $\mathcal{G}(\gamma)$. Nevertheless this is not the case and
assuming that a  measure $\mu$ belongs to $\mathcal{G}(\gamma)$,
one  proves \cite{Ge} that
$$ \mu \Big[\sum_{i \in \Z^d} \omega_i > 1 \Big] =  \mu \Big[\sum_{i \in \Z^d} \omega_i = 1 \Big] =
\mu \Big[\sum_{i \in \Z^d} \omega_i = 0 \Big] = 0 $$
 using different techniques and expressions of the kernel in
the three cases. Thus $\mu(\Omega)=0$,  $\mu$ cannot be a
probability measure and
$\mathcal{G}(\gamma)$ is empty.\\

In this case, the non-existence comes  from the dependence of  the
kernel  on what happens at infinity and this cannot be controlled
by the topology of weak convergence\footnote{We shall see later
that under extra topological properties, one can construct
measures in $\mathcal{G}(\gamma)$ by considering weak limits of
sequences of finite volume probability measures, with random
boundary conditions.}. A good framework to insure existence would
be specifications where this influence is shield out, and the main
one corresponds to specifications that transform local functions
into quasilocal ones, giving rise to the concept of {\em
quasilocal specification}, central in this theory of
infinite-volume Gibbs
measures  as we shall see next chapter.\\

\medskip

\item{{\bf An example of existence and uniqueness: reversible
Markov chain}}:\\

 Let us describe
reversible Markov chains by means of specifications, following
again a presentation of \cite{Pr}. Consider
$\Omega=\{-1,+1\}^{\mathbb{Z}}$ and  a stochastic matrix
\begin{displaymath}
M= \left( \begin{array}{cc}
p&1-p\\
1-q&q\\
\end{array}\right)
\end{displaymath}
with $p>0,\;q>0$ such that $M$ is irreducible and aperiodic. Thus\footnote{Our standard reference for Markov Chains is \cite{CH}.}
\begin{displaymath}
\exists \; \textrm{unique}\; \nu \in
\mathcal{M}^{+}_{1}(E,\mathcal{E}) \; \; \; \textrm{such} \;
\textrm{that} \; \; \; \nu M=\nu.
\end{displaymath}
This defines  an ergodic Markov chain $X=(X_n)_{n \in \mathbb{N}}$
and by Kolmogorov's existence theorem, one defines a unique $\mathbb{P}_\nu \in
\mathcal{M}^{+}_{1}(\Omega)$ s.t. for all $\omega \in
\Omega, k, i_1,\cdots ,i_k \;
\in \mathbb{N}$,
\begin{displaymath}
\mathbb{P}_\nu\big[C_{\omega_{\{i_1,
\cdots,i_k\}}}\big]=\nu(\omega_{i_1}) .
M^{i_2-i_1}(\omega_{i_1},\omega_{i_2}).\cdots.M^{i_k-i_{k-1}}(\omega_{i_{k-1}},\omega_{i_k}).
\end{displaymath}
The ergodicity of the chain is crucial to get uniqueness, using that \cite{CH}
\begin{displaymath}
\forall j,k \in E,\; \lim_{n \to \infty}M^{n}(j,k)=\nu(k)>0.
\end{displaymath}

The considered Markov chain  $X=(X_n)_{n \in \mathbb N}$ is then
the  sequence of random variables on $\big(\{-1,+1\}^{\mathbb
N},\mathcal{E}^{\otimes{\mathbb N}}\big)$ of law $\mathbb{P}_\nu$. Writing its elementary cylinders in the form $C_{\omega_k}=\{X_k(\omega)=i\}=\{\omega_k=i\}$, one has for all
$\omega \in \Omega,\;  k \in \mathbb N ,\;
i_1,\cdots ,i_k\; \in \mathbb N$,
$$\mathbb{P}_\nu\big[C_{\omega_{\{i_1, \cdots,i_k}\}}\big]=\nu(\omega_{i_1}). M^{i_2-i_1}(\omega_{i_1},\omega_{i_2}). \cdots . M^{i_k-i_{k-1}}(\omega_{i_{k-1}},\omega_{i_k})
$$
and gets  the Markov property: $\forall k \in \mathbb N,\;\forall
i,j,\epsilon_{k-1},\cdots,\epsilon_{0}\; \in E$
$$
\mathbb{P}_\nu\big[\omega_{i_{k+1}}=i|\omega_{i_k}=j,\cdots,\omega_0=\epsilon_0\big]\\
=\mathbb{P}_\nu\big[\omega_{i_{k+1}}=i|\omega_{i_k}=j\big]=M(j,i).$$
This Markov chain is  also reversible: $\forall k \in \mathbb
N,\;\forall l \in \mathbb N,\;\forall
i,j,\epsilon_{k+1},\cdots,\epsilon_{k+l}\; \in E$
\begin{eqnarray*}
& & \mathbb{P}_\nu[\omega_{i_{k}} = i|\omega_{i_{k+1}}=j,\cdots,\omega_{k+l}=\epsilon_{k+l}]\\
&=&\frac{\mathbb{P}_\nu[\omega_{i_{k}}=i,\omega_{i_{k+1}}=j,\cdots,\omega_{k+l}=\epsilon_{k+l}]}{\mathbb{P}_\nu[\omega_{i_{k+1}}=j,\cdots,\omega_{k+l}=\epsilon_{k+l}]}
=\frac{\nu(i)M(i,j)\cdots
M(\epsilon_{k+l-1},\epsilon_{k+l})}{\nu(j)M(j,\epsilon_{k+2})\cdots
  M(\epsilon_{k+l-1},\epsilon_{k+l})}\\
&=&\frac{\nu(i)M(i,j)}{\nu(j)}:=N(j,i)
\end{eqnarray*}
where $N$ is then the  stochastic matrix associated to the reverse
chain. Hence, we can extend this chain on
$\Omega=\{-1,+1\}^{\mathbb{Z}}$, and in particular it is still
ergodic. Introduce now a specification $\gamma$ such that
$\mathbb{P}_\nu \in \mathcal{G}(\gamma)$, and compute
\begin{eqnarray*}
\mathbb{P}_\nu\big[\sigma_{\Lambda_n}|\sigma_{\Lambda_{n}^{c}}=\omega_{\Lambda_{n}^{c}}\big]
&=&\frac{\mathbb{P}_\nu[\omega_{]-\infty,-n-1]}\sigma_{\Lambda_n}\omega_{[n+1,+\infty[}]}{\mathbb{P}_\nu[\omega_{]-\infty,-n-1]}\omega_{[n+1,+\infty[}]}
=\frac{\mathbb{P}_\nu[\omega_{-n-1}\sigma_{\Lambda_n}\omega_{n+1}]}{\mathbb{P}_\nu[\omega_{-n-1}\omega_{n+1}]}\\
&=&\frac{\nu(\omega_{-n-1})M(\omega_{-n-1},\sigma_{-n}) \cdots
   M(\sigma_{n},\omega_{n+1})}{\nu(\omega_{-n-1})M^{2n+2}(\omega_{-n-1},\omega_{n+1})}\\
   &=&\frac{M(\omega_{-n-1},\sigma_{-n}) \cdots
   M(\sigma_{n},\omega_{n+1})}{M^{2n+2}(\omega_{-n-1},\omega_{n+1})}.\\
\end{eqnarray*}
Denote then the (finite) normalization by $\mathbf{Z}_{\Lambda_n}(\omega)=M^{2n+2}(\omega_{-n-1},\omega_{n+1})$
and define  a proper kernel $\gamma_{\Lambda_n}$ on
$(\Omega,\mathcal{F})$, for all $\sigma \in \Omega$,  by
$$
\gamma_{\Lambda_n}(\sigma|\omega)=\frac{1}{\mathbf{Z}_{\Lambda_n}(\omega)} \cdot M(\omega_{-n-1},\sigma_{-n})
\cdots M(\sigma_{n},\omega_{n+1}).
$$
One can check that we define thus a specification $\gamma$ such
that $\mathbb{P}_\nu \in \mathcal{G}(\gamma)$: In contrast to the
modelization of the simple random walk described above, the
existence of a DLR measure is then insured. Let us consider now
any $\mu \in \mathcal{G}(\gamma)$ and prove that
$\mu=\mathbb{P}_\nu$. To do so, it is enough to prove for all
$\omega \in \Omega,\; k , i_1,\cdots ,i_k\; \in \mathbb N$,
\begin{displaymath}
\mathbb{P}_\nu\big[C_{\omega_{\{i_1, \cdots,i_k}\}}\big]=\mu\big[(\omega_{i_1},\cdots,\omega_{i_k})\big].
\end{displaymath}
Let us prove it for the one-dimensional cylinder using  Markov
property and consistency: We have for all $x \in \E$ and $n \in \mathbb N$
\begin{eqnarray*}
\mu[\sigma_0=x]&= &\sum_{i \in E,j \in
  E}\mu[\sigma_0=x|\sigma_{-n-1}=i,\sigma_{n+1}=j] \cdot \mu[\sigma_{-n-1}=i,\sigma_{n+1}=j]\\
&=& \sum_{i \in E,j \in
  E}\frac{M^{n+1}(i,x) \cdot M^{n+1}(x,j)}{M^{2n+2}(i,j)} \cdot \mu[\sigma_{-n-1}=i,\sigma_{n+1}=j].
\end{eqnarray*}
Taking now the limit when $n$ goes to infinity, one gets
\begin{eqnarray*}
\mu[\sigma_0=x]&=&\sum_{i \in E,j \in
  E}\frac{\nu(x)\cdot \nu(j)}{\nu(j)} \cdot \mu[\sigma_{-n-1}=i,\sigma_{n+1}=j]\\
&=&\nu(x)\sum_{i \in E,j \in
  E}\mu[\sigma_{-n-1}=i,\sigma_{n+1}=j]\\
&=&\nu(x) =\mathbb{P}_\nu[\sigma_0=x].
\end{eqnarray*}
We obtain the equality of these measures on the other cylinders in
the same way. Thus $\mathcal{G}(\gamma)$ is the
singleton\footnote{Remark that the terms involving the measure
$\nu$ cancel out in  the specification.
Nevertheless, its  invariant character is
encoded in the conditioning yielding a  DLR
measure that depends on $\nu$.}
$\{\mathbb{P}_\nu\}$.\\

\medskip

\item{{\bf An example of phase transition: Ferromagnetic 2d-Ising
model}}:

It is the archetype of original Gibbs specification and we present
it briefly at dimension $d=2$, temperature $T=\frac{1}{\beta} >0$
and no external field, as originally introduced by Lenz to model
ferromagnetism\footnote{Ising analyzed this model in one dimension
in its thesis supervised by Lenz in 1922 \cite{KO}. His only
higher dimensional contribution was the wrong interpretation that
just as in $d=1$, in higher dimension there is no phase
transition.}. To do so, one considers microscopic magnets
$\sigma_i \in E=\{-1,+1\}$ at each site $i \in \mathbb{Z}^2$,  and
to express the fact that two neighbors have a tendency to align,
the following {\em nearest neighbor potential}. It is a family
$\Phi=(\Phi_A)_{A \in \mathcal{S}}$ of $\mathcal{F}_A$-measurable
functions defined for all $\omega \in \Omega$ by

\begin{displaymath}
\Phi_A(\omega)=\left\{
\begin{array}{lll}\; -J \omega_i \omega_j  \;  \; \; \; {\rm if}  \; \; \;  \; A=\langle ij \rangle\\
\\
\; 0 \;  \;  \; \; \;  \; \; \; \rm{otherwise.}\\
\end{array} \right.
\end{displaymath}

The corresponding Hamiltonian at finite volume $\Lambda \in
\mathcal{S}$,  temperature $\beta^{-1} >0$, coupling $J >0$ and
boundary condition $\omega \in \Omega$ is the
well-defined\footnote{It is a finite sum  here but this is not the
case in general. One has usually to check summability conditions
on the potential to define Gibbs measures, as we shall see next
section.} function on $\Omega \times \Omega$ defined by

$$
H_\Lambda^{\beta \Phi}(\sigma \mid \omega)= \sum_{A \in \mathcal{S},A \cap \Lambda
\neq  \emptyset} \; \Phi_A(\sigma_\Lambda \omega_{\Lambda^c}).
$$

The Gibbs specification at inverse temperature $\beta >0$ is the
family of probability kernels $\gamma^{\beta \Phi}=(\gamma^{\beta
\Phi}_{\Lambda})_{\Lambda \in \mathcal{S}}$ given for all $\Lambda
\in \mathcal{S},\; \sigma,\omega \in \Omega$ by
$$
\gamma^{\beta \Phi}_{\Lambda} (\sigma \mid
\omega)=\frac{1}{Z_\Lambda^{\beta \Phi}(\omega)} \; e^{- \beta
H_\Lambda(\sigma \mid \omega)}
$$
where the {\em partition function} ${Z_\Lambda^{\beta
\Phi}(\omega)}$ is a standard normalization depending on the
boundary condition $\omega$. It is indeed a specification due to
the expression of the Hamiltonian in terms of a sum over local
potential terms (see \cite{Ge} or next section). Intensive studies
have established the following:

\begin{theorem}[Phase transition at low T]\label{thmIsing}
There exists $\beta_c >0$ s.t.
\begin{itemize}
\item There exists a unique measure consistent with $\gamma^{\beta
\Phi}$ at high temperatures $\beta < \beta_c$. \item At low
temperatures $\beta > \beta_c$, the set $\mathcal{G}(\gamma^{\beta
\Phi})$ is the Choquet simplex $[\mu_\beta^-,\mu_\beta^+]$ whose
extremal elements are mutually singular and can be selected by the weak limits
$$
\mu_\beta^\pm(\cdot) = \lim_{\Lambda \uparrow \mathcal{S}}
\gamma_\Lambda (\cdot \mid \pm).
$$
 with the magnetizations satisfying
$$
\mu_\beta^+[\sigma_0]=-\mu_\beta^-[\sigma_0]=M_0 >0
$$
\end{itemize}
\end{theorem}

The existence of the weak limits is usually proved here using {\em
correlation} or related (GKS, FKG,etc.) {\em inequalities} valid
for some ferromagnetic systems. The existence of a critical
temperature has been qualitatively established by Peierls in 1936
\cite{Peie,Grif}, using a geometrical computation based on the
energy of {\em contours}, that are circuits in the dual of the
lattice associated to a configuration and whose lengths are
related to its energy. His analysis  gave rise to the powerful
{\em Pirogov-Sinai theory of phase transitions} for more general
models \cite{PS,VEFS}. The exact value of $\beta_c$ is due to
Kramers and Wannier in 1941 \cite{KW}, while Yang got the
magnetization in 1951 \cite{Ya}, both using algebraic
tools\footnote{The magnetization has been conjectured but
unpublished by Onsager in 1949. He also rigorously derived the
free energy in  in 1944 \cite{Ons}.}. The full convex picture,
restricted to translation-invariant measures, has been
independently proved by Aizenmann \cite{Aiz}
 and Higuchi \cite{Hig}, both inspired by considerations on
 percolation raised by Russo \cite{Ru}, described in \cite{GeHi}. This picture has been recently extended to higher dimension by Bodineau \cite{Bod}.
  The fact that $\mathcal{G}(\gamma^{\beta
\Phi})$ is the Choquet simplex $[\mu_\beta^-,\mu_\beta^+]$ means
that any measure $\mu \in \mathcal{G}(\gamma^{\beta \Phi})$ is
uniquely determined by a convex combination of the extreme phases
$\mu^\pm_\beta$, i.e. that there exists a unique $\alpha \in
[0,1]$ s.t. $\mu= \alpha \mu_\beta^- + (1- \alpha) \mu_\beta^+$.
The situation is more complex in higher dimension or on other lattices, as we shall see.
\end{enumerate}
\section{Convexity theory of DLR-measures}

Before introducing  Gibbs measures properly speaking within the
nice topological framework of quasilocality, we first study the
general structure of the set $\mathcal{G}(\gamma)$ of probability
measures consistent with a specification $\gamma$ on a general
Polish probability space $(\Omega,\mathcal{F})$, product of
\footnote{In a more general set-up, $(E,\mathcal{F})$ has to be a
{\em standard Borel space}, see \cite{Ge}.} a finite single-site
state space $(E,\mathcal{E})$. In this case, the family
$\mathcal{C}$ of cylinders has the nice following property of
being a {\em countable core}, following a terminology of
\cite{Ge}:

\begin{definition}[Countable core]\label{Core}
A countable family $\mathcal{C} \subset \mathcal{F}$ is said to be a countable core if it has the following properties:
\begin{enumerate}
\item $\mathcal{C}$ generates $\mathcal{F}$ and is stable under finite
 intersections\footnote{This property corresponds to a $\pi$-system, see next section.}.
\item If $(\mu_n)_{n \in \mathbb{N}}$ is a sequence of $\mathcal{M}_1^+(\Omega)$ such that $\lim_{n \to \infty} \mu_n(C)$
exists for any cylinder $C \in \mathcal{C}$, then there exists a unique  $\mu \in \mathcal{M}_1^+(\Omega)$ that
coincides with this limit on $\mathcal{C}$:
\be \label{coreprop}
\forall C \in \mathcal{C},\; \mu(C)=\lim_{n \to \infty} \mu_n(C).
\ee
\end{enumerate}
\end{definition}
The proof that the family of cylinders is indeed a core relies on Carath\'eodory`s extension theorem \cite{Ge,Will}.
The main purpose of this section is to use this property to provide a general description of $\mathcal{G}(\gamma)$
when it is not an empty set\footnote{Conditions insuring existence are described in Chapter 3.}.

\subsection{Choquet simplex of DLR-measures}

\begin{theorem}\label{pipoint}
Assume that  $\mathcal{G}(\gamma)\neq \emptyset$. Then
$\mathcal{G}(\gamma)$ is a convex subset of
$\mathcal{M}_1^+(\Omega,\mathcal{F})$
 whose extreme boundary is denoted ${\rm ex}  \mathcal{G}(\gamma)$, and satisfies  the following properties:
\begin{enumerate}
\item The extreme elements of $\mathcal{G}(\gamma)$ are the probability measures $\mu \in \mathcal{G}(\gamma)$ that
are trivial on the tail $\sigma$-field $\mathcal{F}_\infty$:
\begin{equation}\label{tailext}
{\rm ex}\mathcal{G}(\gamma)= \Big \{ \mu \in \mathcal{G}(\gamma): \mu(B)=0 \; {\rm or} \; 1, \; \forall B \in    \mathcal{F}_\infty \Big \}.
\end{equation}
Moreover, distinct extreme elements $\mu,\nu \in {\rm ex}\mathcal{G}(\gamma)$ are mutually
singular: $\exists B \in  \mathcal{F}_\infty$, $\mu(B)=1$ and $\nu(B)=0$, and more  generally,
 each $\mu \in \mathcal{G}(\gamma)$ is uniquely determined within $\mathcal{G}(\gamma)$ by its restriction to $\mathcal{F}_\infty$.
\item $\mathcal{G}(\gamma)$ is a  {\em Choquet simplex}: Any $\mu \in \mathcal{G}(\gamma)$ can be written in a unique way as
\begin{equation}\label{Choquet}
\mu = \int_{{\rm ex}\mathcal{G}(\gamma)} \; \nu \cdot \alpha_\mu(d \nu)
\end{equation}
where $\alpha_\mu \in \mathcal{M}_1^+\Big({\rm ex}\mathcal{G}(\gamma),e({\rm ex}\mathcal{G}(\gamma))\Big)$ is defined
for all $M \in e({\rm ex}\mathcal{G}(\gamma))$ by
\be \label{weights}
\alpha_\mu (M) = \mu \Big[ \big \{ \omega \in \Omega: \exists \nu \in M,\; \lim_n \gamma_{\Lambda_n} (C | \omega) = \nu(C) \;{\rm for \; any \; cylinder} \; C \big \} \Big].
\ee
\end{enumerate}
\end{theorem}
In particular, when $M$ is a singleton $\{\nu\} \in e({\rm ex}\mathcal{G}(\gamma))$, (\ref{weights}) reads
\be \label{weightsnu}
\alpha_\mu (\{ \nu \})=\mu \Big[ \big \{ \omega \in \Omega:  \lim_n \gamma_{\Lambda_n} (\cdot | \omega) = \nu(\cdot) \big \} \Big].
\ee

The convexity of $\mathcal{G}(\gamma)$ is trivial and 1. will be a
direct consequence of the following lemma, also crucial for  2. We
recall that $\mu \in \mathcal{G}(\gamma)$ is {\em extreme} iff
$$
\mu=\alpha \nu + (1-\alpha) \bar{\nu} \;\; \;  {\rm with} \; \; \; \alpha \in ]0,1[ \; \; \; {\rm and} \; \; \; \nu,\bar{\nu} \in \mathcal{G}(\gamma) \; \Longrightarrow \; \nu=\bar{\nu}=\mu.
$$

\begin{lemma}\label{lemtaildensity}
Assume that $\mu \in \mathcal{G}(\gamma)$ is such that
$$
\mu=\alpha \nu + (1-\alpha) \bar{\nu},\; \;  {\rm with} \; \; \alpha \in ]0,1[,\; \nu, \bar{\nu} \in \mathcal{M}_1^+(\Omega).
$$
Then $\nu << \mu,\; \bar{\nu} << \mu$ and
\begin{equation}\label{ACtail}
\nu \in \mathcal{G}(\gamma)\; \Longleftrightarrow \; f:\frac{d \nu}{ d \mu} \; \in \mathcal{F}_\infty.
\end{equation}
\end{lemma}

{\bf Proof}: Let $\mu,\nu,\bar{\nu}$ and $\alpha$ as above. The
absolute continuity of $\nu$ w.r.t. $\mu$ comes trivially from the
positiveness of probability measures: Take $A \in \mathcal{F}$ with
$\mu(A)=0$, then $\alpha \nu(A) + (1-\alpha) \bar{\nu}(A)=0$ implies
$\nu(A)=\bar{\nu}(A)=0$ because $0 < \alpha < 1$, and thus $\nu <<
\mu$, $\bar{\nu} << \mu$. Now let us prove the important statement
(\ref{ACtail}). We follow\footnote{More generally, the proof comes
from  \cite{Dy}, but we have rewritten it to avoid the introduction
of too many concepts. These ideas are also related to the so-called
{\em desintegration of measures}, see \cite{Dell,BR,BR2}.} mostly
the  proof of \cite{Ge} and introduce first two $\sigma$-algebra압
related to the specification $\gamma$: The $\sigma$-algebra of
$\gamma$-invariant measurable sets
$$
\mathcal{F}_\gamma = \Big \{ A \in \mathcal{F}: \gamma_\Lambda(A | \cdot) = \mathbf{1}_A(\cdot), \forall \Lambda \in \mathcal{S} \Big \}
$$
and\footnote{They are indeed $\sigma$-algebra압 \cite{Ge}.}
the $\sigma$-algebra of $\mu$-almost surely $\gamma$-invariant measurable sets
$$
\mathcal{F}_\gamma(\mu) = \Big \{ A \in \mathcal{F}: \gamma_\Lambda(A | \cdot) = \mathbf{1}_A(\cdot) \;\;  \mu{\rm -a.s.}, \forall \Lambda \in \mathcal{S} \Big \}.
$$
We first prove that for a given $\mu \in \mathcal{G}(\gamma)$,
\be\label{fgam2}
\nu \in \mathcal{G}(\gamma) \; \Longleftrightarrow \; f:=\frac{d \nu}{ d \mu} \in \mathcal{F}_\gamma(\mu).
\ee
To prove the first part of (\ref{fgam2}),
 it is enough to prove that for all $c \in [0,1]$
 the event $\{ f \geq c \} \in \mathcal{F}_\gamma(\mu)$. For $\Lambda \in \mathcal{S}$,
 we want to prove that, whenever $\mu,\nu \in \mathcal{G}(\gamma)$ and $f=\frac{d \nu}{d \mu}$,
$$
\gamma_\Lambda(f \geq c | \cdot) = \mathbf{1}_{f \geq c}(\cdot), \; \; \;  \mu-a.s.
$$
or equivalently that when $g=\mathbf{1}_{f \geq c}$, one has
\begin{equation}\label{star}
\gamma_\Lambda g=g, \; \; \; \;  \mu-{\rm a.s.}
\end{equation}
Now, $\mu \in \mathcal{G}(\gamma)$ implies $\mu\big[\gamma_\Lambda g - g\big]=0$, so to prove (\ref{star}) it is enough to prove
\begin{equation}\label{sstar}
 \gamma_\Lambda g \leq g, \; \; \;  \; \mu-{\rm a.s.}
\end{equation}
Writing  $\gamma_\Lambda g = (\gamma_\Lambda g) \cdot
\mathbf{1}_{f \geq c} + (\gamma_\Lambda g) \cdot \mathbf{1}_{f <
c} = g \cdot  \gamma_\Lambda g +   \mathbf{1}_{f <  c} \cdot
\gamma_\Lambda g$, one gets that proving \be \label{enough}
(\gamma_\Lambda g) \cdot \mathbf{1}_{f <  c} = 0, \; \; \;\;
\mu-{\rm a.s.} \ee will be enough to get (\ref{sstar}) and thus
(\ref{star}). To do so, let us prove that \be\label{enough2}
\int_{\{f <c\}} (f-c) \cdot \gamma_\Lambda g \; d \mu \geq 0 \ee
i.e.  that \be\label{enough3} \int_{\{f <c\}} f \cdot
\gamma_\Lambda g \; d \mu \geq c \cdot \int_{\{f <c\}}
\gamma_\Lambda g \; d \mu. \ee Using $\nu= f \mu,\; \nu
\gamma_\Lambda = \nu,\; \mu =  \mu \gamma_\Lambda$ and the
expression for $g$, one writes
\begin{eqnarray*}
\int_{f<c} f \cdot \gamma_\Lambda g \;d \mu &=& \int_\Omega f \cdot  \gamma_\Lambda g \;  d \mu  - \int_{f\geq c} f \cdot \gamma_\Lambda g \; d \mu
= \int_\Omega f \cdot \gamma_\Lambda g \; d \mu -  \int_\Omega f \cdot g \cdot \gamma_\Lambda g \; d \mu\\
&=& \int_\Omega  \gamma_\Lambda g \; d \nu -  \int_\Omega f \cdot  g \cdot \gamma_\Lambda g \; d \mu
= \int_\Omega   g \cdot d \nu -  \int_\Omega f \cdot g \cdot \gamma_\Lambda  g \; d \mu\\
&=& \int_\Omega   f \cdot g \cdot  d \mu -  \int_\Omega f \cdot g \; \gamma_\Lambda g \; d \mu
=\int_\Omega   f \cdot g \cdot (1- \gamma_\Lambda g) \;d \mu.
\end{eqnarray*}
But  $f \cdot g = f \cdot \mathbf{1}_{f \geq c} \geq c \cdot \mathbf{1}_{f \geq c} = c \cdot g$, and
 because $0 \leq \gamma_\Lambda g \leq 1$,  one gets
$$ \int_{f<c} f \cdot \gamma_\Lambda g \;  d \mu \geq c \cdot \int_\Omega g \cdot (1- \gamma_\Lambda g) \;d \mu = c \cdot \mu[g] - c \cdot \mu[g \gamma_\Lambda g]=  c \cdot \int_{f <c} \gamma_\Lambda g \; d \mu$$
where the last equality has been obtained using the consistency relation $\mu[g]=\mu[\gamma_\Lambda g]$.
 So (\ref{enough3}) holds, which in turns implies (\ref{enough2}) and then (\ref{enough}) because trivially $(f-c)$
 is strictly negative on the event $\{f < c\}$,  implying thus $\gamma_\Lambda g = 0$  on the same event, that is exactly
 (\ref{sstar}). Thus one has $\gamma_\Lambda g = g \;\; \; \mu$-a.s., and eventually that the density $f \in \mathcal{F}_\gamma(\mu)$.\\

Let us now prove the converse statement, i.e. that when $\mu \in
\mathcal{G}(\gamma)$ and $f=\frac{d \nu}{d \mu}$,
$$
f \in \mathcal{F}_\gamma(\mu) \Longrightarrow \nu= f \cdot \mu \in \mathcal{G}(\gamma).
$$

It is enough to prove it for a step function $f=\mathbf{1}_A$, with $A \in \mathcal{F}_\gamma(\mu)$,
 so let us prove that for all $\Lambda \in \mathcal{S}$ and $A \in \mathcal{F}_\gamma(\mu)$,
 the measure $\nu(\cdot):=\mathbf{1}_A(\cdot) \mu(\cdot)=\mu(\cdot \cap A)$ satisfies $\nu \gamma_\Lambda=\nu$.
By the defining properties of a specification, one can write, for all  $D \in \mathcal{F}$ and $\Lambda \in \mathcal{S}$,
\begin{eqnarray*}
\nu  \gamma_\Lambda (D) =  \int_\Omega \gamma_\Lambda (D | \cdot ) \; d \nu &=& \nu \big[\gamma_\Lambda(A \cap D | \cdot)\big] + \nu \big[\gamma_\Lambda(D \setminus A | \cdot) \big]\\
& \leq& \mu \big[ \gamma_\Lambda(A \cap D | \cdot)\big] + \mu \big[\mathbf{1}_A(\cdot) \gamma_\Lambda(\Omega \setminus A | \cdot ) \big]\\
&=& \mu(A \cap  D) + \mu \big[\mathbf{1}_A(\cdot) \mathbf{1}_{\Omega \setminus A}(\cdot) \big] = \nu(D).
\end{eqnarray*}
Working similarly on  $D^c$ one also gets the domination of $\nu
\gamma_\Lambda (D^c)$ by  $\nu(D^c)$, and together with
$$
\nu \gamma_\Lambda (D) + \nu \gamma_\Lambda (D^c) = 1 = \nu(D) + \nu(D^c)
$$
this implies $\nu \gamma_\Lambda (D)=\nu(D), \; \forall D \in \mathcal{F}$
and $\Lambda \in \mathcal{S}$, so $\nu \in \mathcal{G}(\gamma)$ and    (\ref{fgam2}) holds.\\

To conclude, realize first that $\mathcal{F}_\infty$ is exactly
 the  $\gamma$-invariant $\sigma$-algebra  $\mathcal{F}_\gamma$: Any $A \in \mathcal{F}_\infty$ is $\gamma$-invariant by properness,
  and reciprocally, any $\gamma$-invariant set $A \in \mathcal{F}_\infty$,
   because it can be written  $A=\big\{\gamma_\Lambda(A | \cdot) = 1 \big\}$
   for any $\Lambda \in \mathcal{S}$. Eventually, tail-triviality is obtained
   because  $\mu$ is trivial on $\mathcal{F}_\gamma(\mu)$ if and only if it is
    trivial on  $\mathcal{F}_\gamma=\mathcal{F}_\infty$, the $\mu$-completion of
     the latter being exactly\footnote{It is not  the case for general kernels, properness  of the specification is crucial \cite{Ge}.} $\mathcal{F}_\gamma (\mu)$, see \cite{Ge}. \\


{\bf Proof of Theorem \ref{pipoint}}: \\

{\bf 1.}  It is straightforward to check that
$\mathcal{G}(\gamma)$ is a convex subset of
$\mathcal{M}_1^+(\Omega)$. Then ${\rm ex} \mathcal{G}(\gamma)$ is
non-empty\footnote{Consider for example a regular version of
$\mu[\cdot | \mathcal{F}_\infty]$ for $\mu \in
\mathcal{G}(\gamma)$.} and suppose  $\mu$ is one of its extreme
elements, and that  there exists $B \in \mathcal{F}_\infty$ with
$0 < \mu(B) <1$. Then the conditional probabilities w.r.t. $B$ and
its complement $B^c$ are well-defined as probability measures on
$(\Omega,\mathcal{F})$ in such a way that
$$
\mu(\cdot) = \mu(\cdot | B) \mu(B) + \mu(\cdot | B^c) \mu(B^c)
$$
or equivalently
\begin{equation}\label{contra}
\mu(\cdot)=\alpha \; \mu(\cdot | B) + (1- \alpha) \; \mu(\cdot |
B^c)
\end{equation}
with $\alpha=\mu(B) \in \ ]0,1[$. Denote $\nu(\cdot)=\mu(\cdot | B)$ and rewrite
$$
\nu(\cdot)=\frac{\mu(\cdot \cap B)}{\mu(B)}=\frac{\mathbf{1}_B(\cdot)}{\mu(B)} \;\cdot  \mu(\cdot)
$$
in such a way that $\nu << \mu$, with a density $\frac{d \nu}{ d
\mu}= \frac{\mathbf{1}_B(\cdot)}{\mu(B)}$ that belongs to
$\mathcal{F}_{\infty}$ because $B$ is a tail event. Lemma
\ref{lemtaildensity} proves thus that $\nu(\cdot)=\mu(\cdot | B)$
and $\bar{\nu}(\cdot)=\mu(\cdot | B^c)$ are distinct elements of
$\mathcal{G}(\gamma)$, and together with (\ref{contra}) and
$0<\alpha<1$,  this contradicts the extremality of $\mu$. Thus,
such a tail event $B \in \mathcal{F}_\infty$ cannot exist and one
gets the first part of item 1. of this theorem:

$$\; \; \; \mu \in {\rm ex} \mathcal{G}(\gamma) \; \Longrightarrow \; \mu(B)=0 \; {\rm or} \; 1,\; \forall B \in \mathcal{F}_\infty.$$

\smallskip

To prove the converse statement,
consider $\mu \in \mathcal{G}(\gamma)$, trivial on $\mathcal{F}_\infty$ and such that there exists $\alpha \in ]0,1[$
and $\nu, \bar{\nu} \in {\rm ex} \mathcal{G}(\gamma)$ with $\mu=\alpha \nu + (1-\alpha) \bar{\nu}$.
 Then by Lemma \ref{lemtaildensity}, $\nu << \mu$  with a density  $f:=\frac{d \nu}{ d \mu} \in \mathcal{F}_\infty$.
  The latter is a density thus $\mu[f]=1$ and by tail-triviality of $\mu$ one also has $\mu[f]=f$ ($\mu$-a.s.). Hence $f=1$ ($\mu$-a.s.) and $\nu=\bar{\nu}=\mu$, which is thus an extreme element of $\mathcal{G}(\gamma)$:
$\; \; \; \mu(B)=0 \; {\rm or} \; 1,\; \forall B \in \mathcal{F}_\infty  \; \Longrightarrow \;  \mu \in {\rm ex} \mathcal{G}(\gamma).$

Let us prove now that any $\mu \in \mathcal{G}(\gamma)$ is uniquely determined within $\mathcal{G}(\gamma)$ by its restriction to $\mathcal{F}_\infty$. Consider $\mu,\nu \in \mathcal{G}(\gamma)$ such that $\mu(B)=\nu(B)$ for all $B \in \mathcal{F}_\infty$ and the convex combination $\bar{\mu}=\frac{1}{2} \mu + \frac{1}{2} \nu$, which is also in $\mathcal{G}(\gamma)$. By Lemma \ref{lemtaildensity}  one can write $\mu= f \cdot \bar{\mu}$ and $\nu= g \cdot \bar{\mu}$ with $f,g \in \mathcal{F}_\infty$. If $\nu=\mu$ on $\mathcal{F}_\infty$, then one also has  $\bar{\mu}=\nu=\mu$ on $\mathcal{F}_\infty$, thus $f=g=1$ and $\nu=\mu$. In particular, distinct extreme elements are mutually singular because they are trivial on the tail-$\sigma$-algebra: There exists then $B \in \mathcal{F}_\infty$ such that $\mu(B)=1$ and $\nu(B)=0$.\\

Hence, we have now a characterization of extremality in terms of
tail-triviality. To see how this leads to a unique simplicial
decompostion, we shall also use the following characterization of
extremal DLR measures, which can be derived from (\ref{tailext})
using standard arguments:
\begin{equation}\label{tailext1}
{\rm ex}\mathcal{G}(\gamma)= \Big \{ \mu \in \mathcal{G}(\gamma): \mu[A|\mathcal{F}_\infty]=\mu(A), \; \mu{\rm -a.s.}, \forall A \in \mathcal{F} \Big \}.
\end{equation}

{\bf 2.} To get the unique extreme decomposition of $\mu \in
\mathcal{G}(\gamma)$, we follow the spirit of the proof of Dynkin
\cite{Dy} as worked out in detail by Georgii \cite{Ge} in our
particular DLR case. We shall mention at the end of this chapter
other frameworks where such a decomposition holds.

We start thus from a specification $\gamma$ for which there exists $\mu \in \mathcal{G}(\gamma)$.
By consistency, for any $n \in \mathbb{N}$, the kernel $\gamma_{\Lambda_n}$ is  a regular version of conditional probability of $\mu$ given $\mathcal{F}_{\Lambda_n^c}$:
\begin{equation}\label{DLRcubes}
\gamma_{\Lambda_n}(A | \cdot)=\mu[A|\mathcal{F}_{\Lambda_n^c}](\cdot), \;  \mu-{\rm a.s.},\;\forall A \in \mathcal{F}.
\end{equation}
The existence of such a regular version of conditional probabilities is insured by Theorem \ref{regular},
and we shall sometimes denote formally  $\Omega_\mu$ the set of full $\mu$-measure set on which (\ref{DLRcubes})
 holds for all $A$ and for all $n$.
Remark that the uniformity in $n$ implies that  $\Omega_\mu$ is a tail event.\\

 The {\em backward martingale theorem} \cite{Will} ensures then  that the following almost-sure limit
\begin{equation}\label{BMT}
\lim_{n \to \infty} \; \mu[C|\mathcal{F}_{\Lambda_n^c}](\cdot) = \mu[C | \mathcal{F}_\infty](\cdot),\; \forall C \in \mathcal{C}
\end{equation}
exists also on a full measure set (that can be assumed  to be the same $\Omega_\mu$),
defining a regular version of the conditional probability w.r.t. the tail $\sigma$-algebra $\mathcal{F}_\infty$.
Remark that the regularity of such  versions is encoded in the order of the locutions ''$\mu$-a.-s.'' and ''$\forall A \in \mathcal{F}$''.

Our strategy is now to combine (\ref{DLRcubes}) and  (\ref{BMT})
with the core property (Definition \ref{Core}) to introduce
appropriate objects for a decomposition on the countable family of
cylinders first, to derive some (tail-) measurability properties
using countability, and thereafter to extend in a  standard way
the latter objects onto probability kernels that are in some sense
extreme "$\mu$-almost surely". In fact, the starting point of the
decomposition is one of the defining properties of  versions of
conditional probability with respect to the tail $\sigma$-algebra:
\begin{equation}\label{startdec0}
\forall \mu \in \mathcal{M}_1^+(\Omega),\; \mu(\cdot)=\int_\Omega \mu[\cdot | \mathcal{F}_\infty](\omega) \; d \mu(\omega).
\end{equation}

In some informal sense, the regular versions of $\mu[\cdot | \mathcal{F}_\infty]$ are the prototypes of extremal
measures entering in the decomposition of any $\mu \in \mathcal{G}(\gamma)$. To formalize this using consistency and the
 backward martingale theorem, we need to carefully define appropriate asymptotics of the specification $\gamma$, that will
 be probability kernels with an asymptotic properness leading to tail triviality, that eventually leads to a concentration
 on the extreme elements of $\mathcal{G}(\gamma)$.\\

{\em Step 1: $\mu$-asymptotics of the specification}\\

 We use first the core property (\ref{coreprop}) to $\mu$-almost-surely extend  our asymptotic
 kernels from a definition on the cylinders. Indeed, combining (\ref{DLRcubes}) and  (\ref{BMT}), one gets that
$$
\forall \omega \in \Omega_\mu, \; \lim_{n \to \infty} \;\gamma_{\Lambda_n}(C | \omega)  = \mu[C | \mathcal{F}_\infty](\omega), \; \forall C \in \mathcal{C}.
$$
By the core property, there exists then,  for all $\omega \in \Omega_\mu$,   $\pi^\omega \in \mathcal{M}_1^+(\Omega)$ s.t.
\begin{equation}\label{defkern1}
\forall C \in \mathcal{C}, \pi^\omega(C)= \lim_{n \to \infty} \;\gamma_{\Lambda_n}(C | \omega)= \mu[C | \mathcal{F}_\infty](\omega).
\end{equation}

Extending this construction to any $\omega \in \Omega_\mu^c$ by
requiring $\pi^\omega$ to be any arbitrary elements of
$\mathcal{M}_1^+(\Omega)$, one gets a probability kernel with some
nice specific properties, as seen in the following

\begin{lemma}\label{pipointkern}
Let $\mu_0 \in \mathcal{M}_1^+(\Omega)$, $\mu \in \mathcal{G}(\gamma)$, $\Omega_\mu \in \mathcal{F}_\infty$ as above for which (\ref{DLRcubes}) and  (\ref{BMT}) hold, and define $\; \; \pi^\cdot : \mathcal{F} \times \Omega \longrightarrow [0,1];\;(A,\omega) \longmapsto \pi^\omega(A)$ as follows:
\begin{itemize}
\item $\forall \omega \in \Omega_\mu$, $\pi^\omega$ is the unique element of $\mathcal{M}_1^+(\Omega)$ such that
$$
\forall C \in \mathcal{C},\; \pi^\omega(C)=\lim_n \gamma_{\Lambda_n}(C | \omega)=\mu[C | \mathcal{F}_\infty ](\omega).
$$
\item $\forall \omega \in \Omega_\mu^c$, $\pi^\omega$ is chosen to be the arbitrary probability measure $\mu_0$.
\end{itemize}
Then $\pi^\cdot$ is a probability kernel from $(\Omega,\mathcal{F}_\infty)$ to $(\Omega,\mathcal{F})$ such that
\begin{enumerate}
\item $\; \; \;\; \; \;\; \; \; \; \; \; \; \; \; \;\; \; \; \;\; \; \; \;\; \; \; \;\mu{\rm -a.s.},  \;\mu[A| \mathcal{F}_\infty]=\pi^\cdot(A),  \; \forall A \in \mathcal{F}.$
\item
\begin{equation}\label{pipointmu2}
\big \{ \pi^\cdot \in \mathcal{G}(\gamma) \big \} := \Big \{
\omega \in \Omega: \pi^\omega \in \mathcal{G}(\gamma) \Big \} \;
\in \; \mathcal{F}_\infty.
\end{equation}
and $\pi^\cdot$ is $\mu$-a.s. consistent with $\gamma$ in the sense that:
\begin{equation}\label{pipointmu3}
\mu\Big[\pi^\cdot \in \mathcal{G}(\gamma)\Big]=1.
\end{equation}
\end{enumerate}
\end{lemma}

{\bf Proof:} To prove that $\pi^\cdot$ is a probability kernel
from $(\Omega,\mathcal{F}_\infty)$ to $(\Omega,\mathcal{F})$, we
need to prove (see Definition \ref{Probabker}) that for all
$\omega \in \Omega$, $\pi^\omega( \cdot )$ is a probability
measure on $(\Omega,\mathcal{F})$, and that for all $A \in
\mathcal{F}$, $\pi^\cdot(A)$ is $\mathcal{F}_\infty$-measurable.
The first item is true by construction, thanks to the core
property. To prove the second one, denote $ \mathcal{D}=\big \{ A
\in \mathcal{F}: \pi^\cdot(A) \in \mathcal{F}_\infty \big \}.$ By
construction, $\mathcal{D}$ contains the set $\mathcal{C}$ of
cylinders. The latter is a $\pi$-system (i.e. a family of sets
stable by finite intersections, that generates the
$\sigma$-algebra $\mathcal{F}$) whereas $\mathcal{D}$ is a Dynkin
system (a family of subsets of $\Omega$ containing $\Omega$,
stable by subtractions of subsets and under  monotone limit of
sets) contained in $\mathcal{F}$. Then we use
\begin{lemma}[Dynkin lemma \cite{Will}]\label{dynkin}
Any Dynkin system which contains a $\pi$-system contains the $\sigma$-algebra generated by this $\pi$-system.
\end{lemma}
Thus the property characterizing $\mathcal{D}$ extends to the whole $\sigma$-algebra $\mathcal{F}$, because the former is
generated by the $\pi$-system of the cylinder, and thus: $\forall A \in \mathcal{F}, \; \pi^\cdot(A) \in \mathcal{F}_\infty$.\\

By construction, $\pi^\omega(C)=\mu[C|\mathcal{F}_\infty]$ is true for all cylinders $C$ for $\mu$-almost
 every $\omega$ by (\ref{defkern1}). For these $\omega$`s, $\mu[\cdot | \mathcal{F}_\infty](\omega)$ and
  $\pi^\omega$ are two probability measures that coincide on a $\pi$-system, which coincide then on the $\sigma$-algebra
  generated by this $\pi$-system, which is $\mathcal{F}$ itself here. This proves item 1. of the lemma.\\

Fix now a cylinder $C \in \mathcal{C}$ and focus first on the
(random) measures $\pi^\cdot(C)$. By construction, it inherits
first of all of the properness property and in particular for
$\mu$-almost every $\omega \in \Omega$,
$$
\forall C \in \mathcal{C}_\infty,\;\;\pi^\omega(C)=\mathbf{1}_C(\omega).
$$
It also inherits from consistency: For $\mu$-almost every $\omega$,
$$
\forall C \in \mathcal{C}, \pi^\omega \gamma_\Lambda(C)=\mu \big[\mu[C|\mathcal{F}_{\Lambda^c}]|\mathcal{F}_\infty \big]=\mu[C | \mathcal{F}_\infty]=\pi^\omega(C)
$$
Using Dynkin's lemma and standard extension techniques, this
implies that \be \mu{\rm -a.s.},\; \left\{
\begin{array}{lll}
\pi^\omega \gamma_\Lambda(A)=\pi^\omega(A), \; \; \; \; \forall \Lambda \in \mathcal{S}, \; \forall A \in \mathcal{F}.\\
\\
\pi^\omega(B)=\mathbf{1}_B(\omega),  \; \; \; \; \; \;  \forall B \in \mathcal{F}_\infty.
\end{array} \right.
\ee
In particular, one gets (\ref{pipointmu2}) and (\ref{pipointmu3}) and Lemma \ref{pipointkern} is proved.\\

{\em Step 2: Concentration on the extremal measures}\\

Consider now, for $\mu \in \mathcal{G}(\gamma)$,    $\pi^\cdot$ as
a random measure on the probability space $(\Omega,
\mathcal{F},\mu)$, taking values in the set
$\mathcal{M}_1^+(\Omega)$. By the former lemma, it
$\mu$-concentrates on $\mathcal{G}(\gamma)$ and has some specific
tail-measurable properties that are useful to relate it to
extremal measures. One gets then the  starting point of the
decomposition by rewriting (\ref{startdec0}) for $\mu \in
\mathcal{G}(\gamma)$: \be \label{startdec1} \forall A \in
\mathcal{F}, \; \mu(A)=\int_\Omega \pi^\omega(A) \;
d\mu(\omega)=\mu \big[ \pi^\cdot(A) \big]. \ee

Denote formally  $\alpha_\mu$ the law of $\pi^\cdot$ as a random variable on  the probability space $(\Omega, \mathcal{F},\mu)$. Writing $M_0=\pi^\cdot(\Omega)$, one can rewrite  formally (\ref{startdec1}) in the form

\be \label{startdec2}
\mu(\cdot) = \int_{\Omega}  \pi^\cdot d \mu= \int_{M_0} \nu \; \mu[\pi^\cdot \in d \nu]= \int_{M_0} \nu \; \alpha_\mu[d \nu].
\ee

Before focusing more properly
 on a rigorous definition of the weights $\alpha_\mu$ that leads to the correct decomposition,
 we first establish an extra important consequence of the previous lemma:
 The limiting procedure used to define $\pi^\cdot$ on the space $\Omega_\mu$ of full $\mu$-measure allows interesting
 probabilistic  properties of the measure $\pi^\omega$ for such typical $\omega$: The measure $\pi^\omega$ is an extreme element
 of $\mathcal{G}(\gamma)$ and the above integral reduces to the set $M_0={\rm ex}\mathcal{G}(\gamma)$.

\begin{lemma}\label{pipointexGgam}
$$
\Big \{ \pi^\cdot \in \; {\rm ex}\mathcal{G}(\gamma) \Big\} \in
\mathcal {F}_\infty\; \; {\rm and}, \;  \forall \mu
\in\mathcal{G}(\gamma), \; \mu\Big[\pi^\cdot \in {\rm
ex}\mathcal{G}(\gamma)\Big] = 1.
$$
\end{lemma}
{\bf Proof:} Firstly, recall that the consistency has a consequence on the expected value of $\pi^\cdot(A)$ as a
random variable (with values in $[0,1]$) on $(\Omega,\mathcal{F},\mu)$, as seen in (\ref{startdec1}): $\forall A \in \mathcal{F}$,
$$
\mu \big[ \pi^\cdot(A) \big]=\int_\Omega \pi^\omega(A) \; d \mu(\omega) = \int_\Omega \mu[A | \mathcal{F}_\infty](\omega) \;d \mu(\omega) = \mu \big[ \mu[A | \mathcal{F}_\infty] \big] = \mu(A)
$$
so that the expected value of $\pi^\cdot(A)$ is $\mu(A)$, while its variance under  $\mu \in \mathcal{G}(\gamma)$ is:
\begin{eqnarray}\label{var2}
\mathbb{E}_\mu\Big[\big(\pi^\cdot(A) - \mu(A) \big)^2 \Big]&=&\mu\Big[ \big(\pi^\cdot(A)\big)^2 - 2 \mu(A) \pi^\cdot(A) + \big( \mu(A) \big)^2 \Big]\\
&=& \mu\Big[ \big(\pi^\cdot(A)\big)^2 \Big] - 2 \mu(A)
\mu\big[\pi^\cdot(A) \big] + \big( \mu(A) \big)^2 \nonumber
\end{eqnarray}
in such a way that, when $\mu \in \mathcal{G}(\gamma)$, we can define it to be
\begin{equation} \label{var}
\sigma^2_A(\mu):=\mu \Big[\big(\pi^\cdot(A) \big)^2\Big] -
\Big[\mu \big(\pi^\cdot(A) \big) \Big]^2 =\mu
\Big[\big(\pi^\cdot(A) \big)^2 \Big] - \big(\mu(A)^2 \big)
\end{equation} or
\be \label{evalvar}
\sigma^2_A(\mu)=e_{\big(\pi^\cdot(A)\big)^2}(\mu)-\big(e_{\pi^\cdot(A)}(\mu)\big)^2
\ee that could in particular be used to get (tail) measurability.
For technical reasons, we extend the definition of this map
$\sigma_A^2$ on the whole space $\mathcal{M}_1^+(\Omega)$
using\footnote{Instead of the more usual variance (\ref{var2}). The
two expressions coincide on $\mathcal{G}(\gamma)$} the same
expression (\ref{evalvar}).

By its definition via (\ref{BMT}), $\pi^\cdot(A)$ is a version of $\mu[A \mid \mathcal{F}_\infty]$ for all
$A \in \mathcal{F}$ and
 by (\ref{tailext1}) $\mu$ is  extreme if and only if $\mu[A | \mathcal{F}_\infty] = \mu(A)$, $\mu$-a.s., $\; \forall A \in \mathcal{F}$. This implies that the starting DLR measure $\mu \in \mathcal{G}(\gamma)$ will be extreme  iff  $\;  \forall A \in \mathcal{F}, \pi^\cdot(A)= \mu(A), \;  \mu{\rm -a.s.}$
in such a way that \be \label{tailext2} {\rm
ex}\mathcal{G}(\gamma)=\Big \{ \mu \in
\mathcal{G}(\gamma):\pi^\cdot(A)=\mu(A) \; \; \mu{\rm -a.s.}, \;
\forall A \in \mathcal{F} \Big \} \ee Using now Dynkin`s lemma for
the $\mathcal{D}$-system $\big \{ A \in \mathcal{F}:
\pi^\cdot(A)=\mu(A) \; \; \mu{\rm -a.s.} \big \}$ that contains
$\mathcal{C}$, one gets
\be \label{dec} \Big \{ \pi^\cdot \in {\rm ex} \mathcal{G}(\gamma)
\Big \} = \Big \{ \pi^\cdot \in \mathcal{G}(\gamma) \Big \} \cap
\bigcap_{C \in \mathcal{C}} \Big\{ \pi^\cdot : \pi^\cdot(C)=\mu(C)
\Big\} \ee which in particular insures the
$\mathcal{F}_\infty$-measurability of  $ \big \{ \pi^\cdot \in
{\rm ex} \mathcal{G}(\gamma) \big \}$ from (\ref{pipointmu2}), and
from this of $\pi^\cdot(C)$ for any $C \in \mathcal{C}$. Hence,
extremal measures are the $\mu \in \mathcal{G}(\gamma)$ that
satisfy
$$
\forall C \in \mathcal{C},\; \pi^\cdot(C)=\mu(C),\; \mu-{\rm a.s.}
$$
i.e. that as a random variable, for all $C \in \mathcal{C}$, $\pi^\cdot(C)$ would be a.s. equals to its $\mu$-expectation, and,
 as in many cases in such situations, it implies that its variance (\ref{var}) should be $\mu$-a.s. zero
$$
\mu \in {\rm ex}\mathcal{G}(\gamma) \Longleftrightarrow \mu \in \mathcal{G}(\gamma) \; {\rm and} \; \forall C \in \mathcal{C},\;
  \mu \Big[\sigma^2_{C} \big(\pi^\cdot(C)\big)=0 \Big]=1
$$

This proves that $\pi^\cdot$ is itself extreme $\mu$-almost surely,
 as a consequence of the tail measurability of $\pi^\cdot : A \longmapsto \pi^\cdot(A)$.
 Indeed, one then has $\mu[(\pi^\cdot(A))^2 | \mathcal{F}_\infty]=(\pi^\cdot(A))^2$ $\mu$-a.s., for all $A \in \mathcal{F}$, and  in     particular, for all $C \in \mathcal{C}$,
\be \label{varzeromu}
\mu \Big[ \mu\big[(\pi^\cdot(C))^2 | \mathcal{F}_\infty\big](\omega) - \big(\pi^\omega(C)\big)^2 \Big]=0
\ee
which implies
$$
\mu \Big[ \sigma^2_{\pi^\cdot} (\pi^\cdot \big(C)\big)^2 \Big]=\mu \Big[ \pi^\omega[(\pi^\cdot(C))^2 ]- (\pi^\omega(C))^2 \Big]=0.
$$
This proves $\mu \big[ \pi^\cdot \in {\rm ex}\mathcal{G}(\gamma) \big]=1$ and the lemma using (\ref{dec}).\\

{\em Step 3: Extreme decomposition and its uniqueness}\\

To properly get the decomposition using the concentration of the
asymptotic kernels on the extreme DLR measures, we use the
tail-measurability of the previous lemma together with the very
definition of the conditional expectation in an extended version
of (\ref{startdec0}): \be \label{startdec3} \forall A \in
\mathcal{F}, \; \forall B \in \mathcal{F}_\infty, \int_B \mu[A
\mid \mathcal{F}_\infty](\cdot) \; d \mu(\cdot) = \int_B
\mathbf{1}_A(\cdot) \; d \mu(\cdot). \ee By Lemma
\ref{pipointexGgam}, $B=\big\{ \pi^\cdot \in {\rm
ex}\mathcal{G}(\gamma) \big\} \in \mathcal{F}_\infty$ so in
particular, one has for all $A \in \mathcal{F}$, \be
\label{startdec4} \int_{ \{ \pi^\cdot \in {\rm
ex}\mathcal{G}(\gamma) \}}\mu[A \mid \mathcal{F}_\infty](\cdot) \;
d \mu(\cdot) = \int_{\{ \pi^\cdot \in {\rm ex}\mathcal{G}(\gamma)
\}} \mathbf{1}_A(\cdot) \; d \mu(\cdot)=\mu\big (A \cap \{
\pi^\cdot \in {\rm ex}\mathcal{G}(\gamma) \}\big ) \ee and the
latter is exactly $\mu(A)$ by concentration of $\pi^\cdot$ on
${\rm ex}\mathcal{G}(\gamma)$ for $\mu \in \mathcal{G}(\gamma)$.
Thus one can rewrite (\ref{startdec1}) as \be \label{startdec5}
\mu(\cdot)=\mu\big(\cdot \cap  \{ \pi^\cdot \in {\rm
ex}\mathcal{G}(\gamma) \}\big) = \int_{ \{ \pi^\cdot \in {\rm
ex}\mathcal{G}(\gamma) \}}\mu[ \cdot \mid
\mathcal{F}_\infty](\omega) \; d \mu(\omega) = \int_{\{ \pi^\cdot
\in {\rm ex}\mathcal{G}(\gamma) \}} \pi^\omega(\cdot) \; d
\mu(\omega). \ee Consider $\alpha_\mu$ as in the definition
(\ref{weights}) of the theorem and extend it into a probability
measure $\alpha_\mu \in \mathcal{M}_1^+\big({\rm ex}
\mathcal{G}(\gamma),e({\rm ex} \mathcal{G}(\gamma)\big)$  defined
to be the law of $\pi^\cdot$ as a random extremal DLR measure,
i.e. by
$$
\alpha_\mu(M):=\mu(\pi^\cdot \in M), \; \;
\forall M \in e\big({\rm ex} \mathcal{G}(\gamma)\big).
$$
It is indeed a probability measure because for
 all $M \in e\big({\rm ex} \mathcal{G}(\gamma)\big)$ one has $\{ \pi^\cdot \in M\} \in \mathcal{F}_\infty \subset \mathcal{F}$,
 and because $\alpha_\mu({\rm ex} \mathcal{G}(\gamma)\big)=1$ by step 2.
 above\footnote{Rigorously speaking, one should use the expression $(\ref{evalvar})$ in terms of the evaluation
 maps to prove the measurability on $\big({\rm ex} \mathcal{G}(\gamma),e\big({\rm ex} \mathcal{G}(\gamma)\big)\big)$,
 but this is standard, see e.g. \cite{Ge}.}.
  Thus we identify  $M_0$ by ${\rm ex} \mathcal{G}(\gamma)$ in (\ref{startdec2}) and rewrite (\ref{startdec5})
  under the form
\be \label{startdec6} \mu=\int_{{\rm ex} \mathcal{G}(\gamma)} \;
\nu \; \alpha_\mu(d \nu). \ee Uniqueness of the representation
follows by the uniqueness of the extension in the core property
and from the uniqueness of the representation of a probability
measure via its action on measurable functions, and Theorem
\ref{pipoint} is proved.
\subsection{Selections by boundary conditions}
In statistical physics, Gibbs measures are often considered by
taking the infinite-volume limit  of finite volume specifications
with prescribed boundary conditions. It is not rigorously true for
general DLR measures, but a corollary of the simplicial
decomposition indeed indicates that it is true for extremal
measures\footnote{The converse statement is not true. There exists
non extremal measures that are such weak limits, see e.g. the
3-states Potts models with well-chosen external field. }. The direct
description of non-extremal ones is more peculiar, but of course can
be done using this decomposition, see \cite{Ge,VEFS} for a more
general description. Stronger results are also true in the
quasilocal context, but the latter is not necessary for what
follows; it is important for us while we shall consider
non-quasilocal measures within  the still active Dobrushin program
of restoration of Gibbsianness. Extreme points have thus the nice
general extra property to get identified with some particular
sequence of measures with boundary conditions: An infinite-volume
extremal measure specified by $\gamma$ can be selected by a
sequences of finite volume measures with boundary conditions that
are typical for it:

\begin{theorem}\label{selextbc}\cite{Ge}
Let $\gamma$ be a specification such that there exists $\mu \in {\rm ex}\mathcal{G}(\gamma)$. Then, for any sequence of cubes $(\Lambda_n) \in \mathcal{S}$, for any $f \in \mathcal{F}$ bounded, the following convergence holds:
\begin{eqnarray}
\gamma_{\Lambda_n} f (\cdot)&\;
\mathop{\longrightarrow}\limits_{n\to\infty}\;& \mu[f], \; \; \; \;
\mu{\rm -a.s.}\; \\ \label{selcbc1}
 \gamma_{\Lambda_n}(\cdot | \omega) &\stackrel{W}{
\mathop{\longrightarrow}\limits_{n\to\infty}}& \; \mu(\cdot), \; \;
{\rm for} \; \; \mu{\rm -a.e.}(\omega). \label{selcbc2}
\end{eqnarray}
\end{theorem}
In case of phase transitions, it provides a more explicit description of extremal measures:
\begin{theorem}\label{selextbc2}\cite{Ge}
Let $\gamma$ be a specification such that there exists $\mu \neq \nu\in {\rm ex}\mathcal{G}(\gamma)$. Consider $f \in \mathcal{F}$ bounded such that $\mu[f] \neq \nu[f]$. Then the tail-measurable sets
\begin{eqnarray*}
B_\mu^f &=& \Big \{ \omega: \gamma_{\Lambda_n} f(\omega)\; \;
\mathop{\longrightarrow}\limits_{n\to\infty}\; \; \mu[f] \Big\}\\
B_\nu^f &=& \Big \{ \omega: \gamma_{\Lambda_n} f(\omega)\; \;
\mathop{\longrightarrow}\limits_{n\to\infty}\; \; \nu[f] \Big\}
\end{eqnarray*}
are such that $\mu\big(B_\mu^f\big)=\nu\big(B_\nu^f\big)=1$ and $\nu\big(B_\mu^f\big)=\mu\big(B_\nu^f\big)=0.$
\end{theorem}

{\bf Basic example: 2d Ising model at low temperature}

 Although the simplicial representation is a very satisfactory result from
  a theoretical point of view, it is far from being an easy task to characterize and describe the extreme DLR
  measures for a given specification. We shall mention a few known examples later on in this course,
  but we also stress here that the description is still mostly incomplete from the mathematical point of view.
   The most complete results concern the standard Ising model on $\mathbb{Z}^2$. For this model, we have seen
   in Theorem \ref{thmIsing} that the set of DLR measures is the convex set $\big[\mu_\beta^-,\mu^+_\beta\big]$
   where the extremal measures are characterized by an opposite magnetization $\pm m_\beta \in [0,1]$, defined
   to be $m_\beta = \mathbb{E}_{\mu^+_\beta}[\sigma_0]=-\mathbb{E}_{\mu^-_\beta}[\sigma_0]$, with $m_\beta \neq 0$
   at low enough temperature. In such a case, one can write for any $\mu \in \mathcal{G}(\gamma)$, using (\ref{startdec6})
   and the definition (\ref{weightsnu}) of the weights,
$$
\mu=\alpha_\mu(\{\mu_\beta^+\}) \cdot \mu_\beta^+ +  \alpha_\mu(\{\mu_\beta^-\})  \cdot \mu_\beta^-
$$
where the weights can be shown to satisfy

$$
\alpha_\mu(\mu_\beta^\pm)= \mu \Big[ \big \{ \omega \in \Omega: \lim_n \gamma_{\Lambda_n} (C | \pm) = \mu_\beta^\pm(C) \;{\rm for \; any \; cylinder} \; C \big \} \Big]=\mu(B_{m_\beta^\pm})
$$
with the sets $B_m$ defined in (\ref{tailm}) for $m \in [0,1]$. For
the Gibbs measure with free boundary conditions\footnote{I.e.
without any boundary condition, see next section.}, one e.g.
recovers
$$
\mu^f=\frac{1}{2} \cdot \mu_\beta^+ + \frac{1}{2} \cdot \mu_\beta^-.
$$



\subsection{Ergodic vs. extremal DLR measures}

Thus, Theorem \ref{pipoint} tells us that for a given
specification, any  DLR measure is uniquely determined in terms of
the extremal ones, those that are trivial on the tail
$\sigma$-algebra $\mathcal{F}_\infty$,  i.e. for which global
macroscopic observables do not fluctuate. This is the reason why
they are related to macroscopic states of our system, as we shall
discuss soon. Nevertheless, one is also often interested by
translation-invariant quantities  and it appears that replacing
the tail events by translation-invariant ones  in the previous
decomposition leads to the famous ergodic decomposition of
translation-invariant probability measures. Let us describe
briefly the proof of \cite{Ge} that  gets the ergodic
decomposition as a corollary of the previous theorem, using a
particular specification related to a spatial  average operator.

Hence, as usual in ergodic theory, we focus now on the set
$\mathcal{M}^+_{1,{\rm inv}}(\Omega)$ of {\em
translation-invariant} probability measures  and on the
$\sigma$-algebra $\mathcal{F}_{\rm{inv}}$ of translation-invariant
events. Introducing a particular specification $\tilde{\gamma}$
defined for all $A \in \mathcal{F}, \omega \in \Omega$ and
$\Lambda \in \mathcal{S}$ by
$$
\tilde{\gamma}_\Lambda(A|\omega) = \frac{1}{|\Lambda|} \; \sum_{i \in \Lambda} \mathbf{1}_A(\tau_i \omega)
$$
one easily gets that the $\tilde{\gamma}$-invariant sets are
exactly the translation-invariant ones, i.e.
$$
\mathcal{F}_{\tilde{\gamma}}=\mathcal{F}_{\rm{inv}}.
$$
By an adaptation of  the proof of the extreme decomposition (\ref{startdec5}), in the more general framework of \cite{Dy},
one gets the following theorem, proved in Chapter 14 of \cite{Ge}:
\begin{theorem}\label{ergdec}\cite{Ge,VEFS}
The set $\mathcal{M}_{1,{\rm inv}}^+(\Omega)$ is  a convex subset of $\mathcal{M}_{1}^+(\Omega)$ such that:
\begin{enumerate}
\item Its extreme elements  are the probability measures  that are
trivial on the translation-invariant $\sigma$-algebra
$\mathcal{F}_{\rm{inv}}$, i.e. the {\em ergodic probability
measures} on $(\Omega,\mathcal{F})$:
\begin{equation}\label{invext}
{\rm erg}(\Omega)= \Big \{ \mu \in \mathcal{M}_{1,{\rm inv}}^+(\Omega): \mu(A)=0 \; {\rm or} \; 1, \; \forall A \in    \mathcal{F}_{\rm inv} \Big \}.
\end{equation}
Distinct ergodic measures $\mu,\nu$ are mutually singular: $\exists A \in  \mathcal{F}_{\rm inv}$, $\mu(A)=1$ and $\nu(A)=0$,
 and more  generally, each $\mu \in {\rm erg}(\Omega)$ is uniquely determined within the ergodic measures by its restriction
 on  $\mathcal{F}_{\rm inv}$.
\item $\mathcal{M}_{1,{\rm inv}}^+(\Omega)$ is a Choquet simplex: Any $\mu \in  \mathcal{M}_{1,{\rm inv}}^+(\Omega)$ can
be written in a unique way as
$$
\mu = \int_{{\rm erg}(\Omega)} \; \nu \cdot \alpha_\mu(d \nu)
$$
where $\alpha_\mu \in \mathcal{M}_1^+\Big({\rm erg}(\Omega),e({\rm erg}(\Omega))\Big)$ is defined
for all $M \in e({\rm erg}\mathcal{M})$ by
\be \label{weightserg}
\alpha_\mu (M) = \mu \Big[ \big \{ \omega \in \Omega: \exists \nu \in M,\; \lim_n \tilde{\gamma}_{\Lambda_n} (C | \omega) = \nu(C) \;{\rm for \; any \; cylinder} \; C \big \} \Big].
\ee
\end{enumerate}
\end{theorem}
In this theorem, translation-invariance and the ergodic theorem
play the roles respectively devoted to tail-triviality and the
backward martingale limit theorem  in  Theorem \ref{pipoint}. Many
other similarities exist, described e.g. in  Chapter 14 of
\cite{Ge} or in \cite{VEFS}.

\begin{remark}[Physical Phases] {\em Once we agree to describe the true physical phases of the system by {\em some}
 random field of $\mathcal{M}_{1,{\rm inv}}^+(\Omega)$, we have now three different mathematical manners to characterize the
 macroscopic ''states'' of the systems we want to model. Starting from a specification describing an equilibrium
 at finite volume\footnote{It corresponds to the finite volume {\em Boltzmann-Gibbs weights}, see next chapter.},
 one can consider first the set of DLR measures $\mathcal{G}(\gamma)$ as good candidates to play the same role at infinite volume,
  leaving moreover  the door open to the modelization of phase transitions. Requiring then that macroscopic, i.e. tail-measurable,
   observables should not fluctuate, one can then restrict the macroscopic description to the DLR measures  that are trivial
    on the tail $\sigma$-algebra, i.e. to ${\rm ex} \mathcal{G}(\gamma)$, and use thereafter the extreme
     decomposition to get a more general description that incorporates uncertainty of the experiment. For other purposes,
     one can also be interested in translation-invariant objects, in particular when the underlying system
      is translation-invariant, and chose then to focus either on the
       set $\mathcal{G}_{\rm{inv}}(\gamma)=\mathcal{M}_{1,{\rm inv}}^+(\Omega) \cap \mathcal{G}(\gamma)$ and
       on its extreme elements  ${\rm ex}\mathcal{G}_{\rm{inv}}(\gamma)$, or either on the translation-invariant extremal
       measures, i.e the translation-invariant elements of ${\rm ex}\mathcal{G}(\gamma)$.
       These approaches are far from being equivalent: The latter form a rather small
       and sometimes empty set, whereas the former consists of ergodic measures,
       that are in particular extreme because $\mathcal{F}_{\rm inv} \subset \mathcal{F}_\infty$.
        These ergodic measures are thus often chosen to be the {\em physical phases} that represent
        macroscopically the equilibrium state of the underlying interacting particle system. If one
        does not focus on translation-invariance, the structure of extremal states could be very rich and far from being equivalent to ergodicity, see various examples
        of Ising models at higher dimensions \cite{DOB0}, or antiferromagnetic \cite{VEFS}, or  on trees \cite{BLG,hig,jof,Pres2}.}
\end{remark}

\begin{remark}[de Finetti`s theorem and exchangeability] {\em Let us also mention now that within a slightly different framework, a
theorem similar to Theorem \ref{pipoint} is equivalent to the de Finetti`s theorem for exchangeable
 measures\footnote{See e.g. \cite{Diac}. It can also be used to complete the description of mean-field models.}.
 Instead of working on a lattice, consider $S$ to be the set of non-negative integers $\mathbb{N}$, and focus on
 permutations and its associate $\sigma$-algebra of symmetric events $\mathcal{I}=\cap_{n \in \mathbb{N}}$ already
 introduced in Section 2.2. This symmetric $\sigma$-algebra plays the role devoted to $\mathcal{F}_{\infty}$ in the
 extreme decomposition, and defining a family of proper probability kernels $\gamma=(\gamma_{n \in \mathbb{N}})_{n \in \mathbb{N}}$
 from $(\Omega,\mathcal{I})$ to $(\Omega,\mathcal{F})$ by,
$$
\forall n \in \mathbb{N}, \forall A \in \mathcal{F}, \forall \omega \in \Omega,\; \gamma_n(A | \omega) = \frac{1}{n!} \sum_{\tau \in I_n} \mathbf{1}_A(\tau\omega)
$$
one first gets that the set $\mathcal{G}(\gamma)$ of $\gamma$-invariant probability measures is exactly the set of exchangeable
measures. Using a corollary of Theorem \ref{pipoint}, Its   extreme points have then to be trivial on the symmetric
events $\mathcal{I}$, and it corresponds to  the product measures of the form $\lambda^{\otimes \mathbb{N}}$
with $\lambda \in \mathcal{M}_1^+(E, \mathcal{E})$. Proceeding in a similar way as in the extreme decomposition
 of Theorem \ref{pipoint}, one gets that all permutation-invariant (or exchangeable) measures are uniquely determined as
 convex combinations of product measures: This is exactly de Finetti`s theorem \cite{Diac}.  A more refined analysis of this
  analogy and of the $\sigma$-algebra압 $\mathcal{I}$, $\mathcal{F}_{\rm inv}$ and $\mathcal{F}_\infty$ also lead to
  related 0-1 laws, see \cite{BR,BR2,Dy,Ge}.}
\end{remark}

\chapter{Quasilocal and Gibbs measures}

\section{Quasilocality for measures and specifications}

\subsection{Essential continuity of conditional probabilities}

The link between continuity and quasilocality described through
the product topology enables to generalize the Markov property to
probability measures whose conditional expectations of local
functions depend only weakly of spins arbitrarily far away from
their support. This leads to the concept of {\em quasilocal
measures}, which, in addition to provide a good framework to get
the existence of specified measures, is also closely related to
the notion of Gibbs measures. It corresponds to the  concept of
{\em Feller kernels} in the standard theory of stochastic
processes.

\begin{definition}
A specification $\gamma$ is said to be \emph{quasilocal}
when
$$
\forall \Lambda \in \mathcal{S}, \; \; f \in
\mathcal{F}_{\rm{loc}} \; \Longrightarrow \; \gamma_\Lambda f \in
\mathcal{F}_{\rm{qloc}}. $$
 A measure is said to be {\em quasilocal} if there exists a quasilocal specification
   $\gamma$ such that $\mu \in \mathcal{G}(\gamma)$.
\end{definition}

Recall that for any  $\Lambda \in \mathcal{S}$,  $\gamma_\Lambda f$
is defined by:
$$
\forall \omega \in \Omega,\; \gamma_\Lambda f(\omega)=\int_\Omega
f(\sigma) \gamma_\Lambda(d \sigma \mid \omega).
$$

Thus, $\gamma_\Lambda$ being properly speaking a kernel from
$(\Omega_{\Lambda^c},\mathcal{F}_{\Lambda^c})$ to $(\Omega_{\Lambda}, \mathcal{F}_\Lambda)$, the function
$\gamma_\Lambda f$ is $\mathcal{F}_{\Lambda^c}$-measurable and
continuity as to be understood here as {\em continuity w.r.t. the
boundary condition} $\omega$ (depending only on
$\omega_{\Lambda^c}$): If $\gamma$ is a quasilocal,
then for any $f \in \mathcal{F}_{\rm loc}$, for any $\Lambda \in
\mathcal{S}$,
$$
\lim_{\Lambda' \uparrow \mathcal{S}} \;
\sup_{\sigma_{\Lambda'^c}=\omega_{\Lambda'^c}} \Big|
\gamma_\Lambda f(\omega) - \gamma_\Lambda f(\sigma) \Big| = 0.
$$

An  important consequence on the conditional probabilities of
local or quasilocal functions w.r.t. a quasilocal measure is the
following
\begin{proposition}[Essential continuity of conditional
probabilities]

Consider $\gamma$ to be a quasilocal specification on
$(\Omega,\mathcal{F})$ and $\mu \in \mathcal{G}(\gamma)$. Then,
for all $f \in \mathcal{F}_{\rm{qloc}}$, $\Lambda \in \mathcal{S}$
and $\omega \in \Omega$, there always exists a version of the
conditional probability $\mu[f \mid
\mathcal{F}_{\Lambda^c}](\cdot)$ that is continuous at $\omega$.
\end{proposition}
Indeed, $\mu \in \mathcal{G}(\gamma)$ implies that for all $f \in
\mathcal{F}_{\rm{loc}}$, one has $\mu[f \mid
\mathcal{F}_{\Lambda^c}](\cdot) = \gamma_\Lambda f (\cdot), \; \;
\mu$-a.s., that  has thus to be $\mu$-a.s. continuous. In
particular, for a given quasilocal measure $\mu$, it is not possible
to change a version of conditional probabilities to make it
discontinuous: Take e.g. $f(\sigma)=\sigma_0$, one should have for
all $\Lambda \in \mathcal{S}$ and for all $\omega \in \Omega$,

\begin{equation}\label{essecont}
\lim_{\Lambda' \uparrow \mathcal{S}} \; \sup_{\omega^1, \omega^2
\in \Omega} \Big | \mu[\sigma_0 \mid \omega_{\Lambda' \setminus
\Lambda} \omega^1_{\Lambda'^c}] -\mu[\sigma_0 \mid
\omega_{\Lambda' \setminus \Lambda} \omega^2_{\Lambda'^c}] \Big
|=0
\end{equation}
because the former conditionings are in open neighborhoods of
$\omega$, and open neighborhoods are automatically of positive
$\mu$-measures here \cite{F}. This also express an almost sure
asymptotic weak dependence in the conditioning that can be seen as
an asymptotic extension of Markov properties, as suggested by the
denomination {\em almost Markovian} chosen by Sullivan in
\cite{Su}.

The failure of this essential continuity (\ref{essecont}) will be
very important in the last chapter when dealing with
transformations of Gibbs or quasilocal measures and to detect
non-quasilocality via the following sufficient condition, which we
call  {\em essential discontinuity} although it is a bit stronger
in the following  formulation than the usual general meaning (see
\cite{F}).

\begin{proposition}[Essential discontinuity]
A probability measure $\mu \in \mathcal{M}_1^+(\Omega)$ is {\em
essentially discontinuous} at $\omega$ if there exists $\Lambda
\in \mathcal{S}$, $f \in \mathcal{F}_{\rm{loc}}$, $\delta >0$ and
$\mathcal{N}^1_\Lambda(\omega),\mathcal{N}^2_\Lambda(\omega)$ in
a neighborhood  $\mathcal{N}_\Lambda(\omega)$ such that
$$
\forall \omega^1 \in \mathcal{N}^1_\Lambda(\omega), \omega^2 \in
\mathcal{N}^2_\Lambda(\omega),\; \Big | \mu[ f \mid
\mathcal{F}_{\Lambda^c}](\omega^1) - \mu[f \mid
\mathcal{F}_{\Lambda^c}](\omega^2) \Big |
> \delta
 $$
or equivalently
\begin{equation}\label{essdisc}
\lim_{\Lambda' \uparrow \infty} \sup_{\omega^1,\omega^2 \in
\Omega} \Big | \mu[ f \mid
\mathcal{F}_{\Lambda^c}](\omega_{\Lambda'} \omega^1_{\Lambda'^c})
- \mu[f \mid \mathcal{F}_{\Lambda^c}](\omega_{\Lambda'}
\omega_{\Lambda'^c}^2) \Big | > \delta.
\end{equation}
\end{proposition}

\subsection{Existence results in the quasilocal framework}

In our finite spin-state settings, where compactness holds,
quasilocality insures thus the existence of a measure in
$\mathcal{G}(\gamma)$. The following proposition additionally
indicates how one can naturally construct such objects as the
limit of large but finite systems with some specified boundary
conditions, that have  then to be typical for the measure
constructed. The set of DLR measures is then a closed convex
subset, and this explains the usual introduction of Gibbs measures
as weak limits of finite-volume probability measures with boundary
conditions. For any specification $\gamma$ and sequence
$\big(\nu_n\big)_{n \in \mathbb{N}} \in \mathcal{M}_1^+(\Omega)$,
we recall that $\nu_n \gamma_n$ denotes the probability measure
acting on bounded $f \in \mathcal{F}$ via:
$$
\nu_n \gamma_{\Lambda_n}[f]=\int_\Omega \gamma_{\Lambda_n} f
(\omega) \nu_n(d \omega),\; \forall n \in \mathbb{N}.
$$

\begin{proposition}\cite{Ge}
Let $\Omega$ be a compact metric space and
$\gamma=(\gamma_\Lambda)_{\Lambda \in \mathcal{S}}$ a {\em quasilocal}
specification on it. Then, for any  sequences of cubes
$(\Lambda_n)_{n \in \mathbb{N}}$ and any arbitrary sequence
$(\nu_n)_{n \in \mathbb{N}}$ on $\mathcal{M}_1^+(\Omega)$, the
weak limit
$$
\mu :=\lim_{n \to \infty} \; \nu_n \gamma_n
$$
exists in  $\mathcal{M}_1^+(\Omega)$ and $\mu \in
\mathcal{G}(\gamma)$. In particular, $\mathcal{G}(\gamma)$ is a
non empty convex subset of $\mathcal{M}_1^+(\Omega)$.
\end{proposition}

 We introduce  now the main example of quasilocal measures,
 that are nothing but (infinite-volume) Gibbs measures, and explain why a converse
 statement telling that {\em most} quasilocal measures are Gibbs
 is also true.

\section{infinite-volume Gibbs measures}
\subsection{Equilibrium states at finite volume}

We recall here briefly some elementary physical concepts that led
Boltzmann and Gibbs to settle down their prescription  for
equilibrium states. This is very simplified, and probably too simple
in a physical point of view (the notion of entropy being for example
far from being so simple, see e.g. \cite{LebMaes,Lan}) but we state
it in order to formally justify the notion of equilibrium states
that we will develop in Chapter 4  within the so-called {\em
variational principle}. Hence, our aim is to provide a probabilistic
translation of the second law of thermodynamics that claims:

\begin{center}
{\em Equilibrium at a fixed value of energy maximizes entropy}
\end{center}
or, in an equivalent statement, {\em Equilibrium minimizes free energy}.\\

For a modelization at finite volume $\Lambda$, the microscopic
states are the collections $\sigma_\Lambda \in \Omega_\Lambda$ of
random variables $(\sigma_i)_{i \in \Lambda}$ and the macroscopic
states are their possible distributions $\mu_\Lambda \in
\mathcal{M}_1^+(\Omega_\Lambda,\mathcal{F}_\Lambda)$. The energy of
a configuration is represented by an {\em Hamiltonian} at finite
volume $H_\Lambda(\sigma_\Lambda)$ and thus the energy of a
macroscopic "state" $\mu_\Lambda$ is  represented by the average of
the Hamiltonian, i.e. \be \label{FVEnergy}
\mathbb{E}_\mu[H_\Lambda]:= \sum_{\sigma_\Lambda \in \Omega_\Lambda}
H_\Lambda(\sigma_\Lambda) \mu_\Lambda(d \sigma_\Lambda). \ee In one
of its original interpretations, the entropy of a system is supposed
to evaluate its degree of disorder. Translated into a probabilistic
framework and quoting Khinchin in \cite{Kh}, "it seems highly
desirable to introduce a quantity which in a reasonable way measures
the amount of uncertainty associated with a given probability
measure, that would be minimal for complete uncertainty,  positive
in other cases, maximal for the one with equally likely outcomes
(uniform distribution), and that would have some nice monotone
properties when the knowledge of the system increases". Following
these ideas, one could show that such a function of a measure should
involve the function $f(x)=x \ln(x)$ and the standard definition of
the {\em entropy of a (finite volume) probability measure}
$\mu_\Lambda$ is indeed given by: \be \label{FVEntropy}
\mathcal{H}_\Lambda(\mu)=-\sum_{\sigma_\Lambda \in \Omega_\Lambda}
\mu(\sigma_\Lambda) \ln \mu_\Lambda (\sigma_\Lambda). \ee In
classical thermodynamics, the free energy "$F$" of a system is
usually defined through the second law in the form "$F=U-TS$" where
$U$ is the (internal) energy, $T=\frac{1}{\beta}$ the temperature
and $S$ the entropy. To pursue the analogy, let us define the {\em
free energy at inverse temperature $\beta > 0$ of a (finite-volume)
probability measure} $\mu_\Lambda$ to be \be \label{FVFree}
F^\beta_\Lambda(\mu)=\mathbb{E}_\mu[H_\Lambda] - \frac{1}{\beta} \;
\mathcal{H}_\Lambda(\mu). \ee There exists two simple ways to see
which probability measures could reasonably be considered as
equilibrium states, following the two different statements of the
second law of thermodynamics. In the first formulation in terms of
maximization of entropy, it is an elementary exercise using
Lagrange`s multipliers \cite{Ja} to show that a probability measure
$\mu_\Lambda$ having  the given energy (\ref{FVEnergy}) and
maximizing its entropy (\ref{FVEntropy}) should give weights  of the
form given by the famous {\em Boltzmann-Gibbs weights} \cite{Bo,Gi}
to each configuration: \be \label{BG}
\nu^\beta_\Lambda[\sigma_\Lambda] = \frac{1}{Z_\Lambda^\beta} \cdot
e^{- \beta H_\Lambda(\sigma_\Lambda)} \ee where the normalization is
the partition function $Z_\Lambda^\beta = \sum_{\sigma_\Lambda \in
\Omega_\Lambda} \; e^{- \beta H_\Lambda(\sigma_\Lambda)}$, which
will be related soon to our free energy.

To establish (\ref{BG}) and illustrate the second law of thermodynamics in its second formulation, we introduce
another important concept, {\em the relative entropy of two
probability measures}. For simplicity, we consider the case of
$\mu_\Lambda,\nu_\Lambda$ being two non-null probability measures
on the finite volume configuration space
$(\Omega_\Lambda,\mathcal{F}_\Lambda)$, in the sense that any
configuration has a positive probability; the relative entropy of
$\mu_\Lambda$ with respect to $\nu_\Lambda$ is then defined to be
$$
\mathcal{H}_\Lambda(\mu \mid \nu)=\sum_{\sigma_\Lambda \in
\Omega_\Lambda} \mu_\Lambda (\sigma_\Lambda) \cdot \ln \frac{\mu_\Lambda
(\sigma_\Lambda)}{ \nu_\Lambda (\sigma_\Lambda)}.
$$

This function has among others the nice property to be non-negative for any
probability measures on $\Omega_\Lambda$ and to be zero if and
only the two measures coincide:
\begin{eqnarray*}
\mathcal{H}_\Lambda(\mu \mid \nu) &\geq& 0.\\
\mathcal{H}_\Lambda(\mu \mid \nu) &=& 0 \;\;  {\rm iff} \; \;
\mu_\Lambda=\nu_\Lambda.
\end{eqnarray*}
Observing that  $\mathcal{H}_\Lambda(\mu \mid \nu^\beta) =
F_\Lambda^\beta(\mu) + \frac{1}{\beta} \ln Z_\Lambda^\beta$, one
concludes that free energy is indeed minimal when
$\mu_\Lambda=\nu_\Lambda^\beta$ is given by the Boltzmann-Gibbs
weights (\ref{BG}). This minimal value of the free energy is then
$F_\Lambda^\beta (\nu^\beta) = - \frac{1}{\beta}  \ln
Z_\Lambda^\beta$, recovering thus the other form of the second law
of thermodynamics. These  justify the following introduction of  Gibbs specifications.

\subsection{Gibbs specifications and infinite-volume Gibbs measures}

\begin{definition}[Potential] A {\em potential} is a family $\Phi = (\Phi_{A})_{A \in \mathcal{S}}$
of functions  $$\Phi_A \colon \Omega \longrightarrow \mathbb R$$
indexed by the finite subsets of $S$, such that $\forall A \in
\mathcal{S}$, $\Phi_{A}$ is $\mathcal{F}_{A}$-measurable.
\end{definition}

Our infinite-volume formalism incorporates the finite-volume one by
considering {\em free (or empty) boundary conditions}, to extend
Hamiltonians from $\Omega_\Lambda$ to $\Omega$ in a well defined way
by considering finite sums in the following\footnote{This is true
only when $\Phi_A$ is bounded for all $A \in \mathcal{S}$. In a more
general framework involving ''hard-core exclusion'', $\Phi$ is
allowed to be $\infty$ and the formalism has been adapted, see e.g.
\cite{DOBPe,Lebo}.}

\begin{definition}[Hamiltonian with free boundary condition] Consider a potential $\Phi$. For all $\Lambda \in
\mathcal{S}$, the  {\em Hamiltonian at finite volume
$\Lambda$ with free boundary condition} associated with $\Phi$ is
 the well defined and $\mathcal{F}_\Lambda$-measurable map
\begin{eqnarray*}
\mathbf{H}_{\Lambda}^{\Phi,f} &:& \Omega \longrightarrow \mathbb R \\
& & \omega \longmapsto
\mathbf{H}_{\Lambda}^{\Phi,f}(\omega):=\sum_{A \in \mathcal{S},A
\subset \Lambda} \Phi_{A} (\omega).
\end{eqnarray*}
\end{definition}

Nevertheless, the sums involved in the Hamiltonians are not finite
in general and one should focus first on convergence properties of
potentials before introducing infinite-volume Hamiltonians with
prescribed boundary conditions. In the following definition,
convergence of series will be considered in the sense of the
convergence along nets already defined: A series $\sum_{\Lambda
\in \mathcal{S}} F_\Lambda$ converges iff the net
$\Big(\sum_{\Lambda \in \Delta} F_\Lambda \Big)_{\Delta \in
\mathcal{S}} $ converges to a finite limit as $\Delta \uparrow S$
in the sense of Definition \ref{net}. We shall illustrate this
convergence in some examples soon.

\begin{definition}[Convergence of potentials]
A potential $\Phi$ is said to be
\begin{enumerate}
\item \emph{Nearest neighbor} iff for all $\omega \in
\Omega$, $\Phi_A(\omega)=0$  unless $A=\{i\}$ or $A$ is a pair
$\langle ij \rangle$ of nearest neighbors.

\item \emph{Finite-range} iff there exists a {\em range} $R \in \mathbb{N}^*$ such that, for all $\omega$,
$\Phi_A(\omega)=0$  and for all $A$
such that $|A|  > R$, where $|A|=\sup_{i,j \in
A} d(i,j)$ is the diameter of A.

\item \emph{(Point-wise) convergent} at $\omega \in \Omega$ if,  for all $\Lambda \in \mathcal{S}$, the
Hamiltonian
\be \label{Hamilt}
\mathbf{H}_{\Lambda}^{\Phi} (\omega) :=
\sum_{A \in \mathcal{S},A \cap \Lambda \neq \emptyset} \Phi_{A}
(\omega)
\ee
exists, {\em convergent} when the convergence holds for all $\omega \in \Omega$ and {\em almost-surely convergent}
 when there exists $\mu \in
\mathcal{M}_1^+(\Omega)$ such that $\Phi$ is convergent at
$\mu$-a.e. $\omega \in \Omega$.

 \item \emph{Uniformly convergent} when the series defining
 (\ref{Hamilt})
are uniformly convergent in $\omega \in \Omega$, or equivalently
when
\begin{equation}\label{Unifcgt}
\lim_{\Delta \uparrow \mathcal{S}}\sup_{\omega \in
  \Omega} \Big | \sum_{A \in \mathcal{S},A \cap \Lambda \neq \emptyset,A
  \cap \Delta^{c} \neq \emptyset}\Phi_A(\omega) \Big | =0.
\end{equation}
\item \emph{Uniformly absolutely convergent (UAC)} when
\begin{equation}\label{UAC}
\forall i \in S, \;
\sum_{A \in \mathcal{S},A \ni i} \sup_{\omega \in \Omega} |
\Phi_{A} (\omega) | < + \infty.
\end{equation}
\end{enumerate}
\end{definition}

Nearest neighbor and finite-range potentials are UAC and obviously

\begin{lemma} $\Phi$ UAC $\Longrightarrow \; \Phi$ uniformly convergent
$\Longrightarrow \; \Phi$ convergent.
\end{lemma}

A potential that is UAC satisfies also for any $\Lambda \in
\mathcal{S}$, $\sum_{A \in \mathcal{S},A \cap \Lambda \neq
\emptyset} \sup_{\omega \in \Omega} | \Phi_A(\omega) | < + \infty$
which in particular implies uniform convergence (it corresponds to
{\em normal} convergence of series).

\begin{remark} {\em If we do not make precise the way these infinite sums are done, the
sum  $\mathbf{H}_{\Lambda}^{\Phi}$ in (\ref{Hamilt}) could be
ill-defined. Consider the {\em pair} (but not n.n.) potential
$\Phi$ defined for all $\omega \in \{-1,+1\}^{\mathbb Z}$ by
$\Phi_A(\omega)=\frac{1}{\mid i-j \mid} \; \omega_i\omega_j$ if
 $A=\{i,j\}$ and $\Phi_A = 0$ when $A$ is not a pair\footnote{This potential corresponds to the so-called Coulomb interactions, see e.g. \cite{Simon}.}.
Let $\Lambda \in \mathcal{S}, \omega \in \Omega$, write $S=B^+ \cup
B^-$,  with $B^\pm = B^\pm(\omega)=\{i \in \mathbb Z,\; \omega_i=\pm
1\}$ to get for $A \in \mathcal{S}$
\begin{eqnarray*}
 \sum_{A \cap \Lambda \neq \emptyset} \Phi_{A}
 (\omega)= \sum_{A \cap \Lambda \neq \emptyset,A
 \subset B^+}
 \Phi_{A} (\omega)+ \sum_{A \cap \Lambda \neq
 \emptyset, A
 \subset B^-}
 \Phi_{A} (\omega)
+ \sum_{A \cap \Lambda \neq
 \emptyset,A \cap B^+ \neq \emptyset,A \cap B^- \neq \emptyset} \Phi_{A} (\omega).
\end{eqnarray*}
But,
\begin{displaymath}
 \sum_{A \cap \Lambda \neq \emptyset, A
 \subset B^+} \Phi_{A}(\omega)\; = \sum_{A \cap \Lambda \neq \emptyset,A
 \subset B^-} \Phi_{A}(\omega)=\sum_{i \in \Lambda,j \in \mathbb{Z}} \frac{1}{\mid
 i-j \mid}
\end{displaymath}
are non-convergent series, whereas the series $\sum_{A \cap \Lambda \neq
 \emptyset \\ A \cap B^+ \neq \emptyset,A \cap B^- \neq \emptyset}
 \Phi_{A} (\omega)$
can be convergent for some $\omega$'s. Thus, the series could be
non-convergent whereas $\mathbf{H}_{\Lambda}^{\Phi}$ is
well-defined if we use nets as above.}
\end{remark}

\medskip

{\bf Examples of potentials}
\begin{enumerate}
\item{{\em Ising potentials:}} Recall that the single-spin
state-space is $E=\{-1,+1\}$ with a priori  measure
$\rho_{0}=\frac{1}{2}\delta_{-1}+\frac{1}{2}\delta_{+1}$. The (n.n.)
Ising potential  is $\Phi=(\Phi_A)_{A \in \mathcal{S}}$ defined by
\begin{displaymath}
\Phi_A(\omega)=\left\{
\begin{array}{lllll}
\; -J(i,j) \cdot \omega_i \cdot \omega_j \; \; & &\textrm{if} \; A=\{i,j\}\\
\;-h(i)\cdot \omega_i \; & &\textrm{iff} \; A=\{i\}\\
\; 0 \; & &\textrm{otherwise}
\end{array} \right.
\end{displaymath}

where $J : S \times S \longrightarrow {\mathbb R}$ is called the
coupling function and $h : S \longrightarrow {\mathbb R}$ the
external magnetic field. In the standard Ising model, $J(i,j)=0$ unless $i,j$ are n.n. and when both $J$ and $h$ are constant, we
call it {\em
  homogeneous} Ising model, {\em inhomogeneous} otherwise. It
is {\em ferromagnetic} when $J \geq 0$ and {\em
anti-ferromagnetic} otherwise. Less standard non-n.n. Ising models
are also sometimes considered. One studies e.g. {\em long-range
Ising models} when $J(i,j)=\frac{1}{|i-j|^r}$, well defined for
$r\in ]1,2]$, see next sections for a few results, but also
so-called {\em Kac-Ising potentials} which have
 a long but finite range. Their origin comes from a description of the van der Waals theory of
 liquid-vapor transitions initiated by \cite{Kac} and we follow here the terminology of \cite{AMasi}.
 The starting point is a smooth non-negative function supported by the unit ball $J(\cdot)$ and normalized as a
  probability kernel (i.e. $||J||_1=1$). The Kac interaction allows ranges from 1 to $+ \infty$ via a parameter
   $\gamma>0$ and coupling constants
$$
J_\gamma(i,j)= \gamma^d \cdot J(\gamma |i-j|), \; \; \forall i,j \in S.
$$
In the original van der Waals theory,
one is particularly interested in small $\gamma$ for which the model presents a long range interaction
 (of order $\gamma^{-1}$), small coupling constants (of order $\gamma^d$) and a total strength at each site
 of constant order 1, and in performing thereafter the limit $\gamma$ goes to 0 to approach mean-field models.
 It is the n.n. Ising model when $\gamma=1$, and thus Kac models allow an interplay between n.n. and mean-field models.
\item{{\em A (uniformly) convergent potential that is not UAC:}}

This example is due to Sullivan \cite{Su}. Consider
$\Omega=\{-1,+1\}^{\mathbb Z}$ and define a potential $\Phi$ that is
non-null only for the finite sets of adjacent sequences in $\mathbb
Z$ on which the spins are all $+1$, more precisely  such that for
all $A \in \mathcal{S}$, for all $\omega \in \Omega$,
$\Phi_A(\omega)=\frac{(-1)^n}{n^2} \; \textrm{ iff } \;
\omega_i=+1,\; \forall i \in A=\{k,\cdots ,k+n-1\},  k  \in \mathbb
Z,n \in \mathbb{N}^{*}$, and $\Phi_A=0$ otherwise. To prove that
$\Phi$ is a convergent potential, we  prove that  the series
$\mathbf{H}^{\Phi}_{\Lambda}(\omega)=\sum_{A \cap \Lambda \neq
  \emptyset,A \in \mathcal{S}}\Phi_A(\omega)$ are convergent for
$\Lambda=\{0\}$, the extension to all finite subsets $\Lambda$
being then straightforward. This amounts to  prove that, for all
$\omega \in \Omega$, the sequence of general term
$U_n(\omega)=\sum_{A \ni 0,A \cap \Lambda_{n}^c \neq \emptyset}
\Phi_A(\omega)$ converges to zero when n goes to infinity. Here it
becomes
\begin{displaymath}
\sum_{A \ni 0,A \cap \Lambda_{n}^c \neq
  \emptyset}
  \Phi_A(\omega)= \sum_{k>n}\sum_{A \ni
  0,|A|=k}\frac{(-1)^k}{k^2} \cdot \prod_{i \in
  A}\mathbf{1}_{\{\omega_{i}=+1\}}(\omega).
\end{displaymath}
and one has for $n$ large enough, uniformly in $\omega$,
\begin{displaymath}
0 \; \leq \big| U_n(\omega) \big| = \Big| \sum_{A \ni 0,A \cap
\Delta_{n}^c \neq
  \emptyset} \Phi_A(\omega) \Big| \leq \Big| U_n(+) \Big|=\Big |
  \sum_{k>n}(k+1)\frac{(-1)^k}{k^2} \Big |
\end{displaymath}
where the term on the right is the tail of a convergent alternating
series, which is convergent. This potential is thus uniformly convergent,
but it is not uniformly absolutely convergent (UAC):
\begin{eqnarray*}
 \sum_{A \ni 0,A \subset
  \Lambda_{n}}\sup_{\omega \in \Omega} | \Phi_A(\omega) |
  =\sum_{n=0}^{2k}(k+1) \Big | \frac{(-1)^{k}}{k^2} \Big | =\sum_{k=0}^{2n}\frac{k+1}{k^2}
\end{eqnarray*}
and the latter  is a non-convergent series.

\item {\em A.s. convergent potential}: Potentials associated to
renormalized measures to define them as weakly Gibbsian measures,
see Chapter 5 or \cite{MRM}.

\item {\em (Relatively) uniformly convergent potential}: Sullivan
\cite{Su} has introduced a notion of convergence slightly weaker
than uniform convergence, which can be  associated to any
quasilocal specification and which is translation-invariant when
the specification is. See next section and Remark \ref{rksull}.

\item {\em UAC}: We describe next section how to define a UAC potential from a quasilocal specification,
following a general construction of Kozlov \cite{Ko}.
\end{enumerate}

Before introducing Gibbs measures properly speaking, we give a
general  definition of particular potentials that will be used to
build a convergent potential associated to a quasilocal
specification in the forthcoming Theorem \ref{Gibbsrepthm}, and
later on in Chapter 5 to establish thermodynamic properties in the
generalized Gibbsian framework.

\begin{definition}[Vacuum potential] Let $\Phi$ be a potential and
denote by $\mathbf{+}$ a particular\footnote{This could be any configuration, we denote it by "+" only
 by analogy with the Ising model.}
 configuration of $\Omega$. We say that $\Phi$ is
a {\em vacuum
  potential} with vacuum state $\mathbf{+} \in \Omega$ iff
$\Phi_A(\omega)=0$ whenever $\omega_i=+$ for some $i \in A \in
\mathcal{S}$.
\end{definition}
For such a potential, also called {\em lattice gas potential}, the
Hamiltonian with free boundary conditions can be seen as an
Hamiltonian with the vacuum state as boundary condition, when
Hamiltonians with boundary condition are defined by the following
\begin{definition}[Hamiltonian at volume $\Lambda$ with boundary condition $\omega$]
If  $\Phi$ is a convergent potential, the
{\em Hamiltonian at volume $\Lambda \in \mathcal{S}$ with boundary condition $\omega \in \Omega$} is defined for all $\sigma \in \Omega$ by
\be \label{Hamiltbc}
\mathbf{H}_{\Lambda}^{\Phi,\omega}
(\sigma)=\mathbf{H}_{\Lambda}^{\Phi} (\sigma \mid \omega)
:=\mathbf{H}_{\Lambda}^{\Phi}
(\sigma_{\Lambda}\omega_{\Lambda^{c}})=\sum_{A \in \mathcal{S},A \cap \Lambda \neq \emptyset} \phi_A(\sigma_{\Lambda}\omega_{\Lambda^{c}}).
\ee
\end{definition}

A convergent potential is regular enough to define this
Hamiltonian with boundary conditions, but it will not be enough to
define  Gibbs measures with the right expected properties, for
which UAC is usually required. At finite volume $\Lambda$, the
Hamiltonian with free boundary conditions of a configuration
$\sigma$ is seen as the energy of the system contained in
$\Lambda$ when it is in the configuration $\sigma$, and a UAC
convergence means  that a change of a configuration in a finite
part of the infinite system produces always a finite change of the
total energy. Requiring for a potential to be UAC will be enough
to define a Gibbsian specification associated with this potential,
and then to provide a ''reasonable'' modelling of the physical
properties of the system\footnote{See also a general discussion
about Banach spaces of interactions in \cite{VEFS}.}. This
requirement actually seems to be too strong, and this possibly too
strong requirement causes troubles in the analysis of some
renormalization group transformations, leading to generalized
Gibbs measures described in Chapter 5.

We are now ready to introduce Gibbs specifications and Gibbs measures,
defined from a UAC potential. First we introduce the following
normalization, central in statistical physics and related to the
free energy of the system as we shall see in Chapter 4.

\begin{definition}[Partition function]\label{partitionZ}
Let $\Phi$ be a convergent potential, $\omega \in \Omega$, $\beta
> 0$ and $\Lambda \in \mathcal{S}$. We call \emph{partition
  function} at temperature $\beta^{-1}$, at volume $\Lambda$, with
potential $\Phi$ and boundary condition $\omega$, the
$\mathcal{F}_{\Lambda^c}$-measurable function
\begin{displaymath}
\mathbf{Z}_{\Lambda}^{\beta\Phi}(\omega)=\int_{\Omega} e^{-\beta
\mathbf{H}_{\Lambda}^{\Phi}(\sigma)} \rho_{\Lambda} \otimes
\delta^{\otimes
\Lambda^{c}}_{\omega_{\Lambda^{c}}}(d\sigma)=\int_{\Omega}
e^{-\beta
\mathbf{H}_{\Lambda}^{\Phi}(\sigma)}\kappa_{\Lambda}(d\sigma)=\int_{\Omega_{\Lambda}}e^{-\beta
\mathbf{H}_{\Lambda}^{\Phi}(\sigma |\omega)}
\rho_{\Lambda}(d\sigma_{\Lambda})
\end{displaymath}
where  $\kappa_{\Lambda}=\rho_{\Lambda} \otimes \delta^{\otimes
  \Lambda^{c}}_{\omega_{\Lambda^{c}}} \in \mathcal{M}_1^+(\Omega)$,
 and $\delta_x$ is the Dirac measure on $x \in E^\Lambda$.
\end{definition}
When free boundary conditions are considered, the partition function is denoted $Z_\Lambda^{\beta \Phi,\rm{f}}$.

\begin{definition}[Gibbs distribution at finite volume $\Lambda$]
Let $\Phi$ be a UAC potential. For $\Lambda \in \mathcal{S}$, we
call the \emph{Gibbs distribution at finite volume $\Lambda$},
with potential $\Phi$, at temperature $\beta^{-1}$ and with
boundary condition $\omega \in \Omega$, the probability measure
$\gamma_{\Lambda}^{\beta \Phi}(\cdot |\omega)$ on
$(\Omega,\mathcal{F})$ defined by:
$$
\forall A \in \mathcal{F}, \gamma_{\Lambda}^{\beta \Phi}(A |
\omega)=\frac{1}{\mathbf{Z}_{\Lambda}^{\beta \Phi}(\omega)}
\int_{\Omega} \mathbf{1}_{A}(\sigma)  e^{-\beta
\mathbf{H}_{\Lambda}^{\Phi}(\sigma)} \;  \kappa_{\Lambda}(d\sigma)
$$
where $\kappa_{\Lambda}=\kappa_\Lambda^\omega$ still denotes the
product measure $\rho_{\Lambda}\otimes \delta^{\otimes
\Lambda^{c}}_{\omega_{\Lambda^{c}}}$ on $(\Omega,\mathcal{F})$.
\end{definition}
In order to underline the role of the
boundary condition $\omega$, one also writes
\begin{displaymath}
\gamma_{\Lambda}^{\beta \Phi}(A |
\omega)=\frac{1}{\mathbf{Z}_{\Lambda}^{\beta \Phi}(\omega)}
\int_{\Omega_{\Lambda}}\mathbf{1}_{A}(\sigma_{\Lambda}\omega_{\Lambda^c}) e^{-\beta
\mathbf{H}_{\Lambda}^{\Phi}(\sigma|\omega)} \;\rho_{\Lambda}(d\sigma_{\Lambda})
\end{displaymath}

\begin{theorem}[Gibbs specification]\label{Gibbsspeql}
Let $\Phi$ be a UAC potential and $\beta >0$. The family of
kernels $\gamma^{\beta \Phi}=(\gamma_\Lambda^{\beta
\Phi})_{\Lambda \in \mathcal{S}}$ is a specification, called {\em
Gibbs specification with (UAC) potential $\Phi$, at inverse
temperature $\beta >0$}.
\end{theorem}

{\bf Proof}:  It is straightforward to prove that for a UAC
potential, the Hamiltonian with boundary condition is bounded, and
thus the partition function exists as a function of the boundary
condition. Define now for all $A \in \mathcal{F}$, $\Lambda \in
\mathcal{S}$ and $\sigma\in \Omega$, the density-type function \be
\label{denspe}
f_{\Lambda}(\sigma)=\frac{1}{\mathbf{Z}_{\Lambda}^{\beta
\Phi}(\sigma)}\cdot e^{-\beta
\mathbf{H}_{\Lambda}^{\Phi}(\sigma)}. \ee

There is no boundary condition $\omega$ involved in this function,
although it incorporates the partition function. This  being
$\mathcal{F}_{\Lambda^c}$-measurable, one  nevertheless recovers
$$ \forall \sigma,\omega \in \Omega, \;
f_\Lambda(\sigma_\Lambda
\omega_{\Lambda^c})=\frac{1}{Z_\Lambda^{\beta \Phi}(\omega)} \cdot
e^{-\beta \mathbf{H}_{\Lambda}^{\Phi}(\sigma \mid \omega)}.
$$ We have also $ 0 < | f_\Lambda
| \leq 1$ and for all $A \in \mathcal{F}$ and $\omega \in
\Omega$,
\be \label{fdenscons}
\gamma_{\Lambda}^{\beta \Phi}(A
|\omega)=\int_{A} f_{\Lambda}(\sigma) \; \kappa_{\Lambda}^\omega(d\sigma)
=\int_{\Omega_{\Lambda}} \mathbf{1}_{A}(\sigma_{\Lambda}\omega_{\Lambda^c}) \cdot f_{\Lambda}(\sigma_{\Lambda}\omega_{\Lambda^c}) \; \rho_{\Lambda}(d\sigma_{\Lambda}).
\ee
and it is straightforward to check that $\gamma_\Lambda^{\beta
\Phi}$ is a probability kernel satisfying properties  1. in
 Definition \ref{specification} of a specification. Properness is also directly verified: Let $B \in \mathcal{F}_{\Lambda^{c}}$.  $\forall \sigma,\omega \in
\Omega, \;
 \mathbf{1}_{B}(\sigma_{\Lambda}\omega_{\Lambda^c})$ is independent of
 $\sigma$ and $\mathbf{1}_{B}(\sigma_{\Lambda}\omega_{\Lambda^c})=\mathbf{1}_{B}(\omega_{\Lambda}\omega_{\Lambda^c})=\mathbf{1}_{B}(\omega)$.
Therefore, for all $\omega \in \Omega$ and $B \in
\mathcal{F}_{\Lambda^c}$,
\begin{eqnarray*}
& & \gamma_{\Lambda}^{\beta\Phi}(B|\omega)=\frac{1}{\mathbf{Z}_{\Lambda}^{\beta\Phi}(\omega)}
\int_{\Omega_{\Lambda}}\mathbf{1}_{B}(\sigma_{\Lambda}\omega_{\Lambda^c}) e^{-\beta
\mathbf{H}_{\Lambda}^{\Phi}(\sigma|\omega)}
\rho_{\Lambda}(d\sigma_{\Lambda})\\
&=&\frac{1}{\mathbf{Z}_{\Lambda}^{\beta \Phi}(\omega)}
\int_{\Omega_{\Lambda}}
\mathbf{1}_{B}(\omega_{\Lambda}\omega_{\Lambda^c}) e^{-\beta
\mathbf{H}_{\Lambda}^{\Phi}(\sigma|\omega)} \rho_{\Lambda}(d\sigma_{\Lambda})
=\frac{\mathbf{1}_{B}(\omega)}{\mathbf{Z}_{\Lambda}^{\beta
\Phi}(\omega)} \int_{\Omega_{\Lambda}}e^{-\beta
\mathbf{H}_{\Lambda}^{\Phi}(\sigma|\omega)}  \rho_{\Lambda}(d\sigma_{\Lambda})=\mathbf{1}_{B}(\omega).
\end{eqnarray*}

To prove consistency (\ref{DLRcons}), we assume
without any loss of generality that
$\beta=1$ and consider $\Lambda \subset \Lambda' \in \mathcal{S}$, $A \in \mathcal{F}$ and $\omega
 \in \Omega$. To prove that $\gamma_{\Lambda'}(A|\omega)=\gamma_{\Lambda'}\gamma_{\Lambda}(A|\omega)$, we
write
\begin{displaymath}
\gamma_{\Lambda'}(A|\omega)=\int_{\Omega_{\Lambda'}}\mathbf{1}_{A}(\tau_{\Lambda'}\omega_{\Lambda'^c})f_{\Lambda'}
(\tau_{\Lambda'}\omega_{\Lambda'^c})\rho_{\Lambda'}(d\tau_{\Lambda'})
\end{displaymath}
and $\; \; \; \; \; \;\gamma_{\Lambda'}\gamma_{\Lambda}(A|\omega)=\int_{\Omega}\gamma_{\Lambda}(A|\tau)\; \gamma_{\Lambda'}(d\tau|\omega)$
\hspace{2cm}
$$
=\int_{\Omega_{\Lambda'}}\Big(\int_{\Omega_{\Lambda}}\mathbf{1}_{A}(\sigma_{\Lambda}\tau_{\Lambda'\backslash\Lambda}
\omega_{\Lambda'^c})f_{\Lambda}(\sigma_{\Lambda}\tau_{\Lambda'\backslash\Lambda}\omega_{\Lambda'^c})
d\sigma_{\Lambda}\Big) \cdot f_{\Lambda'}(\tau_{\Lambda'}\omega_{\Lambda'^c}) \; d\tau_{\Lambda'}
$$
where we have written $d\sigma_{\Lambda}$ instead of
$\rho_{\Lambda}(d\sigma_{\Lambda})$. To prove that the family is invariant under expectations
w.r.t. the conditioning in intermediate regions, we modify the
latter expression in order to extract what is really needed in
terms of the functions $f_\Lambda$. Write
\be \label{prt}
\gamma_{\Lambda'}\gamma_{\Lambda}(A|\omega)=\int_{\Omega_{\Lambda'\backslash\Lambda}}g_{\Lambda,\Lambda'}(\tau_{\Lambda'\backslash\Lambda})d\tau_{\Lambda'\backslash\Lambda}
\ee
where, using Fubini`s theorem and trivial changes of variables, one has
\begin{eqnarray*}
g_{\Lambda,\Lambda'}(\tau_{\Lambda'\backslash\Lambda})
&=&\int_{\Omega_{\Lambda}}f_{\Lambda'}(\tau_{\Lambda}\tau_{\Lambda'\backslash\Lambda}\omega_{\Lambda'^c})  \cdot \Big(\int_{\Omega_{\Lambda}}
\mathbf{1}_A(\sigma_{\Lambda}\tau_{\Lambda'\backslash\Lambda}\omega_{\Lambda'^c})f_{\Lambda}(\sigma_{\Lambda}
\tau_{\Lambda'\backslash\Lambda}\omega_{\Lambda'^c})d\sigma_{\Lambda}\Big) \; d\tau_{\Lambda}\nonumber\\
&=&\int_{\Omega_{\Lambda}}\mathbf{1}_A(\tau_{\Lambda}\tau_{\Lambda'\backslash\Lambda}\omega_{\Lambda'^c}) \cdot
f_{\Lambda}(\tau_{\Lambda'}\omega_{\Lambda'^c}) \cdot \Big(\int_{\Omega_{\Lambda}}f_{\Lambda'}(\sigma_{\Lambda}
\tau_{\Lambda'\backslash\Lambda}\omega_{\Lambda'^c})d\sigma_{\Lambda}\Big) \; d\tau_{\Lambda}
\end{eqnarray*}

To get consistency at the level of the density functions
$(f_\Lambda)_{\Lambda \in \mathcal{S}}$, one would like to get rid
of the last integral in this expression. This is provided by the
following lemma, which indicates under which conditions on the
family $\big(f_\Lambda\big)_{\Lambda \in \mathcal{S}}$ of densities
one recovers consistency at the level of specifications. It is in
fact crucial to check one of its items to get consistency for the
Gibbs kernels\footnote{This property of specifications will be very
useful to play on the conditioning for boundary conditions that
coincide outside some finite sets, in particular to get the Kozlov's
potential next section. It corresponds to the {\em key
bar-displacement} property of \cite{F}, where densities of
specifications are explicitly introduced.}.

\begin{lemma}[Consistency for densities]\label{keylemma}
Let $(f_{\Lambda})_{\Lambda \in \mathcal{S}}$ be a family of
(strictly) positive measurable functions $f_{\Lambda}$ such that $\forall \Lambda \in \mathcal{S}$, $\forall
\omega \in \Omega$,
$\int_{\Omega_{\Lambda}}f_{\Lambda}(\sigma_{\Lambda}\omega_{\Lambda^{c}})\rho_{\Lambda}(d\sigma_{\Lambda})=1$.
The following statements are equivalent:
\begin{enumerate}
\item $\forall \Lambda \subset \Lambda'  \in \mathcal{S},\;
\forall \omega, \omega' \in \Omega \; \textrm{s.t.} \;
\omega_{\Lambda^c}=\omega'_{\Lambda^c}$,
\begin{equation}\label{keybar}
\frac{f_{\Lambda'}(\omega')}{f_{\Lambda'}(\omega)}=\frac{f_{\Lambda}(\omega')}{f_{\Lambda}(\omega)}.
\end{equation}
\item $\forall \Lambda \subset \Lambda'  \in \mathcal{S},\;
\forall \omega \in \Omega$,
\begin{equation}\label{keybar2}
f_{\Lambda'}(\omega)=f_{\Lambda}(\omega) \cdot \int_{\Omega_{\Lambda}}f_{\Lambda'}(\sigma_{\Lambda}\omega_{\Lambda^c})
\rho_{\Lambda}(d\sigma_{\Lambda}).
\end{equation}
\end{enumerate}
\end{lemma}

{\bf Proof }: Let us prove that 1.
$\Longrightarrow$ 2., writing $d\sigma_{\Lambda}$ instead of
$\rho_{\Lambda}(d\sigma_{\Lambda})$. Assume (\ref{keybar}) holds for $\Lambda \subset
\Lambda' \in \mathcal{S}$ and let $\omega \in \Omega$. Then
\begin{eqnarray*}
f_{\Lambda}(\omega)\int_{\Omega_{\Lambda}}f_{\Lambda'}(\sigma_{\Lambda}\omega_{\Lambda^c})d\sigma_{\Lambda}
&=&\int_{\Omega_{\Lambda}}f_{\Lambda}(\omega)f_{\Lambda'}(\sigma_{\Lambda}\omega_{\Lambda^c})d\sigma_{\Lambda}
=\int_{\Omega_{\Lambda}}f_{\Lambda'}(\omega)f_{\Lambda}(\sigma_{\Lambda}\omega_{\Lambda^c})d\sigma_{\Lambda}\\
&=&f_{\Lambda'}(\omega) \int_{\Omega_{\Lambda}}f_{\Lambda}(\sigma_{\Lambda}\omega_{\Lambda^c})d\sigma_{\Lambda}=f_{\Lambda'}(\omega)
\end{eqnarray*}
because
$\int_{\Omega_{\Lambda}}f_{\Lambda}(\sigma_{\Lambda}\xi_{\Lambda^c})d\sigma_{\Lambda}=1$. Thus $1.\Longrightarrow 2.$

Consider now $\omega,\omega',\Lambda,\Lambda'$ as above, with
$\omega_{\Lambda^c}=\omega'_{\Lambda^c}$. Using
\begin{displaymath}
f_{\Lambda'}(\omega)=f_{\Lambda}(\omega) \cdot \int_{\Omega_{\Lambda}}  f_{\Lambda'}(\sigma_{\Lambda}\omega_{\Lambda^c}) \; d\sigma_{\Lambda}, \; \; \;
f_{\Lambda'}(\omega')=f_{\Lambda}(\omega') \cdot \int_{\Omega_{\Lambda}}f_{\Lambda'}(\sigma_{\Lambda}\omega'_{\Lambda^c}) \; d\sigma_{\Lambda}
\end{displaymath}
with  $\omega_{\Lambda^c}=\omega'_{\Lambda^c}$, we get
\begin{displaymath}
\int_{\Omega_{\Lambda}}f_{\Lambda'}(\sigma_{\Lambda}\omega_{\Lambda^c}) \; d\sigma_{\Lambda}=\int_{\Omega_{\Lambda}}f_{\Lambda'}(\sigma_{\Lambda}\omega'_{\Lambda^c}) \; d\sigma_{\Lambda}
\end{displaymath}
and then
\begin{displaymath}
f_{\Lambda'}(\omega') \cdot f_{\Lambda}(\omega) \cdot \int_{\Omega_{\Lambda}}f_{\Lambda'}(\sigma_{\Lambda}\omega_{\Lambda^c}) \;d\sigma_{\Lambda}=
f_{\Lambda}(\omega') \cdot \Big(\int_{\Omega_{\Lambda}}f_{\Lambda'}(\sigma_{\Lambda}\omega_{\Lambda^c}) \; d\sigma_{\Lambda}\Big)\cdot f_{\Lambda'}(\omega)
\end{displaymath}
and we conclude the proof of the lemma by non-nullness of our Gibbs weights.\\

To prove Theorem \ref{Gibbsspeql}, we check now that item 1. is true when $f_{\Lambda}$ is given by (\ref{denspe}).
Consider  $\Lambda \subset \Lambda' \in \mathcal{S}$, $\omega$ and $\omega'$ s.t.
$\omega_{\Lambda^c}=\omega'_{\Lambda^c}$. By definition
$$
\frac{f_{\Lambda'}(\omega')}{f_{\Lambda'}(\omega)}=\Big(\frac{\mathbf{Z}_{\Lambda'}(\omega')}{\mathbf{Z}_{\Lambda'}(\omega)}\Big)^{-1}
\cdot \frac{\exp(-\sum_{A \cap \Lambda' \neq
\emptyset}\Phi_A(\omega'))}{\exp(-\sum_{A \cap \Lambda' \neq
\emptyset} \Phi_A(\omega))}.
$$
But, by $\mathcal{F}_A$-measurability of $\Phi_A$, for $A \subset
\Lambda$, $\omega'_{\Lambda^c}=\omega_{\Lambda^c}$ implies
$\Phi_A(\omega)=\Phi_A(\omega')$, and thus
\begin{eqnarray*}
\frac{\exp(-\sum_{A \cap \Lambda' \neq
\emptyset}\Phi_A(\omega'))}{\exp(-\sum_{A \cap \Lambda' \neq
\emptyset} \Phi_A(\omega))}=e^{-\sum_{A \cap \Lambda' \neq
\emptyset}(\Phi_A(\omega')- \Phi_A(\omega))} =e^{-\sum_{A \cap
\Lambda \neq \emptyset}(\Phi_A(\omega')- \Phi_A(\omega))}.
\end{eqnarray*}
 The ratio of the partition functions is also the same for such $\omega$ 's and
 $\omega'$ 's
and eventually (\ref{keybar}) holds. This implies  Item 2. of the lemma: For all $\Lambda
\subset \Lambda' \in \mathcal{S}$, for all $\tau,\omega \in
\Omega$,
\begin{displaymath}
f_{\Lambda}(\tau_{\Lambda'}\omega_{\Lambda'^c}) \cdot \Big(\int_{\Omega_{\Lambda}}f_{\Lambda'}(\sigma_{\Lambda}
\tau_{\Lambda'\backslash\Lambda}\omega_{\Lambda'^c})d\sigma_{\Lambda}\Big)=f_{\Lambda'}(\tau_{\Lambda'}\omega_{\Lambda'^c})
\end{displaymath}
and thus (\ref{prt}) holds with
\begin{displaymath}
g_{\Lambda,\Lambda'}(\tau_{\Lambda'\backslash\Lambda})=\int_{\Omega_{\Lambda}}\mathbf{1}_A(\tau_{\Lambda'}\omega_{\Lambda'^c}) \cdot f_{\Lambda'}(\tau_{\Lambda'}\omega_{\Lambda'^{c}})d\tau_{\Lambda}
\end{displaymath}
yielding consistency
\begin{eqnarray*}
\gamma_{\Lambda'}\gamma_{\Lambda}(A|\omega)&=&
\int_{\Omega_{\Lambda'\backslash\Lambda}}  \Big(\int_{\Omega_{\Lambda}}\mathbf{1}_A(\tau_{\Lambda'}\omega_{\Lambda'^c}\Big) \cdot f_{\Lambda'}(\tau_{\Lambda'}\omega_{\Lambda'^c})d\tau_{\Lambda})d\tau_{\Lambda'\backslash\Lambda}\\
&=&\int_{\Omega_{\Lambda'}}\mathbf{1}_A(\tau_{\Lambda'}\omega_{\Lambda'^c})f_{\Lambda'}(\tau_{\Lambda'}\omega_{\Lambda'^c})d\tau_{\Lambda'}
=\gamma_{\Lambda'}(A|\omega).
\end{eqnarray*}

The relationships between potentials and Gibbs specifications is
not one-to-one: local changes in the potentials can be made
without affecting the kernels $\gamma^{\beta \Phi}$, yielding
equivalent descriptions of measures. This leads to the concept of
physical equivalence\footnote{Other equivalence classes and spaces
of potentials exist, see \cite{VEFS,Is}. In particular, to get the
following equivalence it is crucial to focus on UAC potentials.}.

\begin{definition}[Physical equivalence]
Two potentials $\Phi$ and $\Phi'$ are {\em physically equivalent}
if the Gibbs kernels $\gamma_\Lambda^{\beta \Phi}$ and
$\gamma_\Lambda^{\beta \Phi'}$ are the same for all $\Lambda \in
\mathcal{S}$.
\end{definition}

\begin{definition}[Gibbs measures]\label{Gibbsmeasure}
A probability measure $\mu \in \mathcal{M}_1^+(\Omega)$ is said to
be {\em a Gibbs measure} if there exists a UAC potential $\Phi$
and $\beta >0$ such that $\mu \in \mathcal{G}(\gamma^{\beta
\Phi})$. We often say that $\mu$ is a Gibbs measure for the UAC
potential $\Phi$.
\end{definition}

{\bf Examples of Gibbs measures and phase transitions:}
\begin{enumerate}
\item{{\bf One-dimensional homogeneous Ising models}}
\begin{enumerate}
\item{{\em Ferromagnetic n.n.}}: We have already described in the
previous chapter how ergodic Markov chains could be described as
Gibbs measures for the homogeneous Ising model. The converse is also
possible, and is indeed achieved in a general framework via the
introduction of stochastic matrices defined in terms of the n.n.
potential. The one-dimensional  n.n. Ising model with coupling $J>0$
and external magnetic field $h$ is described and analyzed in this
way in \cite{Ge} using a matricial formalism that leads  in
particular to the well known absence of phase transitions in one
dimension. When $h>0$, uniqueness is proved with a  Gibbs measures
$\mu_\beta^+$ of positive magnetization, that converges weakly to
the Dirac measure at the all $+$ configuration ($\delta_+$) when the
temperature goes to zero, while  opposite  measures $\mu_\beta^-$
and $\delta_-$ are reached when $h<0$. In absence of magnetic field
$h=0$, the unique Gibbs measure is the neutral convex combination
$\mu_\beta=\frac{1}{2} \mu_\beta^+ + \frac{1}{2} \mu_\beta^-$, which
weakly converges to the convex combination $\frac{1}{2} \delta^+ +
\frac{1}{2} \delta^-$, exhibiting a so called {\em asymptotic loss
of tail-triviality} responsible of the phase transition observed in
dimension 2. \item{{\em Anti-ferromagnetic n.n.}}: The same
matricial formalism is also used to deal with the anti-ferromagnetic
case $J<0$. At high magnetic field $|h|>2$, the $+$-phase is the
unique one and  weakly converges, when the temperature goes to zero,
to $\delta_+$ and the opposite situation occurs when $|h| <2$. The
boundary case $|h|=2$ also leads to uniqueness, leading to a unique
phase $\mu_F$ described in terms of Fibonacci's numbers reflecting
an highly non-trivial phenomenon: As claimed in \cite{Ge}, in spite
of the existence of infinitely many ground states\footnote{These are
minimizers of the Hamiltonian, useful to describe the phases at zero
temperature and by extension to low temperatures within the
Pirogov-Sinai theory, see \cite{PS}.} there is no asymptotic loss of
tail   triviality. This loss occurs when the magnetic field is
lower, $|h|<2$, where one gets as a unique Gibbs measure a convex
combination $\mu=\frac{1}{2} \mu^\pm_\beta + \frac{1}{2} \mu^{\mp}$
of two symmetric measures whose typical configurations have either
mostly pluses on a the (say) odd sublattice and mostly minuses on
the even one. When the temperature goes to zero, $\mu$ weakly
converges to a similar Dirac measures $\frac{1}{2} \delta_\pm +
\frac{1}{2} \delta_\mp$, see again a detailed analysis in \cite{Ge}.
\item{{\em Long-range one dimensional Ising models}}: The potential
has already been introduced in the beginning of this section. When
the polynomial decay $r=1$, it is not UAC, but a formalism that
corresponds to so-called Coulomb interactions can be developed
within the weaker notion of uniform convergence, see\cite{Simon} and
also the previous chapter. For $r>1$, this has been studied by e.g.
\cite{Dys,spi2} and it leads in particular to phase transitions in
one dimension when $1<r<2$. Phase transition also occurs in the case
$r=2$, yielding a particular decay of correlations known as a
Thouless effect, see \cite{FS}.
\end{enumerate}
\item{{\bf 3d Ising models}}: We shall be laconic to describe this
very important example of mathematical statistical mechanics:
Theorem \ref{thmIsing} is not valid in dimension $d \geq 3$ and they
do exist non-translation-invariant extreme Gibbs measures. This has
been achieved by Dobrushin in \cite{DOB01}, with a shorter proof in
\cite{Bei}, and these non-translation-invariant so-called {\em
Dobrushin states} are related to the stability of an interface
between a $+$-like phase and a $-$-like one, and to each interface
corresponds an extremal Gibbs measure, in addition to the usual $+$-
and $-$-phases. This example is very relevant for comparing the
notions of ergodic and extremal Gibbs measures discussed at the end
of the previous chapter. \item{{\bf Ising models on Cayley trees}}:
This example is also very interesting from the ergodic vs. extreme
point of view and there also exists an (uncountable) infinite number
of extremal Gibbs measures at low temperature, and depending on the
temperature there could exist two or three translation-invariant
ones, see all the work done in \cite{BLG,hig,jof,Pres2} and a whole
chapter in \cite{Ge}.
\item{{\bf Kac-models}}: A careful adaptation of the Peierls
argument allows to establish the occurrence of a phase transition
for this model at low temperature and long enough range, see
\cite{BZ,CP} in dimension $d \geq 2$.
\begin{theorem}\label{KacIsing}
For $d \geq 2$, for any $\beta >1$, there exists
$\gamma=\gamma(\beta)$ such that for all $\gamma < \gamma(\beta)$,
there exists at least two distinct DLR measures $\mu_\gamma^- \neq
\mu^+_\gamma$.
\end{theorem}
\end{enumerate}

Now that our central objects are properly defined,
we can prove a previous claim providing Gibbs measures as the
main example of quasilocal measures. It is the easiest part of the
link between these two notions, a partial converse statement will
be established next Section.

\begin{theorem}\label{Gibbsqloc}
Let $\beta >0$ and $\Phi$ be a UAC potential. Then the Gibbs
specification  $\gamma^{\beta \Phi}$ is quasilocal. Thus, any
Gibbs measure is also quasilocal.
\end{theorem}

{\bf Proof:} If $\Phi$ be UAC
potential it implies in particular that:
\begin{equation}\label{uac2}
\sum_{A \in \mathcal{S},A \cap \Lambda \neq \emptyset} \sup_{\omega \in \Omega} | \Phi_A(\omega)| < + \infty
\end{equation}
which in turns implies that  for all $\Lambda \in \mathcal{S}$, $H_\Lambda^\Phi$ is a
quasilocal function. Indeed, if $\mathcal{S} \ni \Lambda' \supset
\Lambda$ and consider two configurations $\sigma$ and $\omega$
such that $\sigma_{\Lambda'}=\omega_{\Lambda'}$, we have
$$
\big | \mathbf{H}_{\Lambda}^{\Phi}(\omega)-
\mathbf{H}_{\Lambda}^{\Phi}(\sigma) \big |
\leq  2 \sum_{A \in \mathcal{S},A \cap \Lambda \neq \emptyset,A
\cap \Lambda'^c \neq \emptyset} \sup_{\omega \in \Omega} | \Phi_A(\omega)|
$$
and the latter converges to zero as a consequence of (\ref{uac2}). Thus, one gets the quasilocality of the Hamiltonians:
\begin{displaymath}
\lim_{\Lambda' \uparrow \mathcal{S}} \sup_{\sigma,\omega \in
\Omega,\sigma_{\Lambda'}=\omega_{\Lambda'}}\mid
\mathbf{H}_{\Lambda}^{\Phi}(\omega)-\mathbf{H}_{\Lambda}^{\Phi}(\sigma)
\mid \; = \; 0.
\end{displaymath}
Quasilocality of Gibbs specifications follows.

\begin{remark}[Uniform convergence and quasilocality] {\em Requiring for a potential to be uniformly absolutely
convergent is actually too strong a requirement for merely proving
the quasilocality of the Gibbs specification. Uniform convergence
is actually enough to prove the quasilocality of the Hamiltonian.
In such case, one has
\begin{displaymath}
\sup_{\sigma,\omega \in
  \Omega,\sigma_{\Lambda'}=\omega_{\Lambda'}} \big |\mathbf{H}_{\Lambda}^{\Phi}(\omega)- \mathbf{H}_{\Lambda}^{\Phi}(\sigma) \big |
   \; \leq \; 2 \sup_{\omega \in \Omega} \Big | \sum_{A \in \mathcal{S},A \cap \Lambda
\neq \emptyset,A \cap \Lambda'^{c} \neq
   \emptyset}\Phi_A(\omega) \Big |
\end{displaymath}
and
\begin{displaymath}
\lim_{\Lambda' \uparrow \mathcal{S}}\sup_{\sigma,\omega \in
  \Omega} \Big | \sum_{A \in \mathcal{S},A \cap \Lambda \neq \emptyset,A
  \cap \Lambda'^{c} \neq \emptyset}\Phi_A(\omega) \Big | =0
\end{displaymath}
means the uniform convergence of this
potential. Thus, when the potential is uniformly convergent, the
Hamiltonian is a well-defined quasilocal function and so is the
specification.}
\end{remark}
\begin{remark}[Non-Gibbsianness and essential discontinuity]{\em Let $\mu$ be a Gibbs measure: By theorem \ref{Gibbsqloc}, there exists a quasilocal specification $\gamma$ s.t. $\mu \in \mathcal{G}(\gamma)$ and
$$
\forall A \in \mathcal{F},\; \mu[A
\mid\mathcal{F}_{\Lambda^{c}}](\cdot)=\mathbb
E_{\mu}[\mathbf{1}_{A}\mid
\mathcal{F}_{\Lambda^{c}}](\cdot)=\gamma_{\Lambda}(A \mid\cdot) \;
\mu\textrm{-a.s.}
$$
so that there exists always \emph{one
continuous version}, as a function of
  the boundary condition $\omega$, of the conditional probabilities
  of $\mu$ with respect to the $\sigma$-algebra generated by the
  outside of finite sets. This will be used in Chapter 5 to detect
  non-Gibbsianness by proving the existence of special configurations
  that are point of essential discontinuities, for which there exists conditional
  expectations of local functions that have no continuous version.}
\end{remark}

A Gibbs specification is quasilocal but the converse is not
true in general. However, {\em most} of the quasilocal
specifications are Gibbsian, and we make this now.

\subsection{Gibbs representation theorem}

In this section, we want to characterize a Gibbs measure at the level of specifications: Let $\mu$ be a DLR measure, i.e. such that there
is a specification $\gamma$ with $\mu \in \mathcal{G}(\gamma)$.
To characterize $\mu$ as a Gibbs measure, one should manage to express the
weights of configurations in an exponential form and in some sense
every configuration should receive a non-zero weight. One says
that the specification has to be {\em non-null} in the following sense:

\begin{definition}[Uniform non-nullness]
A specification $\gamma$ is said to be  {\em uniformly non-null} iff $\forall \Lambda \in
\mathcal{S}, \; \exists \; \alpha_{\Lambda},\; \beta_{\Lambda}$
with $0<  \alpha_{\Lambda}  \leq \beta_{\Lambda}  <
\infty$ s.t.
\begin{equation}\label{uninnonnull}
0 < \alpha_{\Lambda} \cdot  \rho(A) \; \leq \; \gamma_{\Lambda}(A \mid
\omega) \; \leq \beta_{\Lambda} \; \rho(A), \; \forall \omega \in
\Omega, \forall A \in \mathcal{F}.
\end{equation}
\end{definition}

If $\gamma$ is quasilocal, non-nullness, in the sense that
$\rho(A)>0 \; \Longrightarrow \;\gamma_{\Lambda}(A \mid \omega)>
\;0$, for all $\omega \in \Omega$, is equivalent to uniform
non-nullness \cite{Ge}. It is also-called the {\em finite energy
condition} in percolation circles. A measure $\mu$ is then said to
be
 (uniformly) non-null if there exists a (uniformly) non-null specification $\gamma$ such that $\mu \in
\mathcal{G}(\gamma)$.\\

 We are now ready to give a partial converse statement of Theorem \ref{Gibbsqloc}. For the purpose of this theorem, the inverse temperature $\beta$ has been incorporated in the potential.
\begin{theorem}[Gibbs representation theorem \cite{Ko,F}]\label{Gibbsrepthm} Let $\mu$ be a quasilocal and uniformly non-null
probability measure on $(\Omega,\mathcal{F})$. Then $\mu$ is a
Gibbs measure, i.e. there exists a UAC potential $\Psi$ such that $\mu
\in\mathcal{G}(\gamma^{\Psi})$.
\end{theorem}

\textbf{Proof}:  Let $\mu$ non-null and quasilocal: There exists a
quasilocal specification $\gamma$ such that $\mu \in
\mathcal{G}(\gamma)$  and (\ref{uninnonnull}) holds, and let us
try to guess which necessary property  a potential should have
such for  $\gamma=\gamma^{\beta \Phi}$ to hold, by considering
such a Gibbs specification $\gamma$ first. Among all the
physically equivalent potentials that define this specification,
let us also assume for the moment that a vacuum potential $\Phi^+$
exists, with a vacuum state denoted by +. The vacuum property and
its link with  free boundary conditions  will be very useful to
relate the specification and the potential. Indeed, considering
$\Lambda \in \mathcal{S}$, then  one obviously has
$H_\Lambda^{\Phi^+}(+ | +)=0$ by the vacuum property and thus
$$
\gamma_\Lambda(+|+)=\frac{1}{Z_\Lambda (+)}
$$
so for  any other configuration $\sigma \in \Omega$
$$
\gamma_\Lambda(\sigma|+)=\gamma_\Lambda(+|+) \; e^{-
H_\Lambda^{\Phi^+}(\sigma | +)}
$$
To exploit consistency via Lemma (\ref{keylemma}), we introduce
the density $f_\Lambda(\sigma):=\gamma_\Lambda(\sigma | \sigma)$.
 Consistency implies that  it satisfies the conditions of Lemma \ref{keylemma},
 and we shall use indifferently both expressions, in terms of $f$ or in
  terms\footnote{The formulation in terms of the density $f_\Lambda$
  is handy to use consistency via Lemma \ref{keylemma}, while the expression
  in terms of $\gamma$ is more familiar. Fern{\'a}ndez \cite{F} has introduced
  densities for specifications and has expressed Lemma \ref{keylemma} in terms of $\gamma$ directly.} of $\gamma$.
Then, by non-nullness and the defining equation (\ref{denspe}),
$$
H_\Lambda^{\Phi^+}(\sigma | +)=-\ln
\frac{\gamma_\Lambda(\sigma|+)}{\gamma_\Lambda(+|+)}=\ln \frac{f_\Lambda(+)}{f_\Lambda(\sigma_\Lambda +_{\Lambda^c})}.
$$
It then possible to derive a vacuum
potential from this Hamiltonian, mainly because conditioning
prescribing a vacuum boundary condition  is equivalent to consider the
Hamiltonian with free boundary condition, for which the use of an
inversion formula from Moebius is direct. Indeed, the vacuum
condition yields
\begin{equation}\label{freebc}
H_\Lambda^{\Phi^+}(\sigma | +)=
 \sum_{A \cap \Lambda \neq \emptyset} \Phi^+_A(\sigma_\Lambda
+_{\Lambda^c})=\sum_{A \subset \Lambda} \Phi^+_A (\sigma) + \sum_{A
\cap \Lambda \neq \emptyset,A \cap \Lambda^c \neq \emptyset}
\Phi^+_A(+_\Lambda +_{\Lambda^c})
\end{equation}
where the last sum is null by the vacuum property. Thus
$$
\forall \Lambda \in \mathcal{S},\; \forall \sigma \in \Omega,\;
H_\Lambda^{\Phi^+} (\sigma \mid +)=H_\Lambda^{\Phi^+,f}(\sigma):=
\sum_{A \subset \Lambda} \Phi^+_A(\sigma)).
$$
In particular one gets directly the single-site potentials:
\begin{equation}\label{singlesite}
\forall i \in S, \forall \sigma \in \Omega,\; \Phi^+_{\{i\}}
(\sigma)= H_{\{i\}}^{\Phi^+,f}(\sigma) = - \ln
\frac{\gamma_{\{i\}}(\sigma | +)}{\gamma_{\{i\}}(+ | +)}.
\end{equation}

To get an insight of the mechanism of the Moebius  inversion
formula, which will enables us to rewrite $\Phi^+$ from $\gamma$,
let us use (\ref{freebc}) to derive the potential for finite
regions consisting of two and three sites, for a fixed $\sigma$
that we forget in the notation. For $\Lambda=\{i,j\}$, write

$$
H_{\{i,j\}}^{\Phi^+,f} = \Phi^+_{\{i\}} + \Phi^+_{\{j\}} +
\Phi^+_{\{i,j\}}
$$
so, using the single-site expression (\ref{singlesite}), one gets  for all
$\sigma \in \Omega$
\begin{eqnarray*}
\Phi^+_{\{i,j\}}&=& H_{\{i,j\}}^{\Phi^+,f} - H_{\{i\}}^{\Phi^+,f}
- H_{\{j\}}^{\Phi^+,f}\\
&=& - \ln\frac{\gamma_{\{i,j\}}(\sigma | +)}{\gamma_{\{i,j\}}(+ |
+)} + \ln \frac{\gamma_{\{i\}}(\sigma | +)}{\gamma_{\{i\}}(+ | +)}
+ \ln \frac{\gamma_{\{j\}}(\sigma | +)}{\gamma_{\{j\}}(+ | +)}.
\end{eqnarray*}
For $\Lambda = \{i,j,k\}$, write similarly, thanks to the vacuum
condition,
\begin{eqnarray*}
H_{\{i,j,k\}}^{\Phi^+,f} &=& \Phi^+_{\{i\}} + \Phi^+_{\{j\}} +
\Phi^+_{\{k\}} + \Phi^+_{\{i,j\}} + \Phi^+_{\{i,k\}} +
\Phi^+_{\{j,k\}} +  \Phi^+_{\{i,j,k\}}\\
&=&\Phi^+_{\{i\}} + \Phi^+_{\{j\}} + \Phi^+_{\{k\}} +
H_{\{i,j\}}^{\Phi^+,f} - H_{\{i\}}^{\Phi^+,f} -
H_{\{j\}}^{\Phi^+,f}\\ &+& H_{\{i,k\}}^{\Phi^+,f} -
H_{\{i\}}^{\Phi^+,f} - H_{\{k\}}^{\Phi^+,f} +
H_{\{j,k\}}^{\Phi^+,f} - H_{\{j\}}^{\Phi^+,f} -
H_{\{k\}}^{\Phi^+,f} + \Phi^+_{\{i,j,k\}}\\
&=& - H_{\{i\}}^{\Phi^+,f} - H_{\{j\}}^{\Phi^+,f}
 - H_{\{k\}}^{\Phi^+,f} + H_{\{i,j\}}^{\Phi^+,f} +
 H_{\{i,k\}}^{\Phi^+,f} +  H_{\{j,k\}}^{\Phi^+,f} + \Phi^+_{\{i,j,k\}}
\end{eqnarray*}
and thus
$$
\Phi^+_{\{i,j,k\}}= H_{\{i,j,k\}}^{\Phi^+,f} -
H_{\{i,j\}}^{\Phi^+,f} - H_{\{i,k\}}^{\Phi^+,f}  -
H_{\{j,k\}}^{\Phi^+,f} + H_{\{i\}}^{\Phi^+,f}
+H_{\{j\}}^{\Phi^+,f}+ H_{\{k\}}^{\Phi^+,f}.
$$
Proceeding by induction, one could reconstruct the potential in
this way. It is actually formally proved using the following
formula, proved in this way e.g. in \cite{F}.

\begin{proposition}[Moebius "inclusion-exclusion" inversion
formula] Let $\mathcal{S}$ be a countable set of finite sets and
$H=(H_\Lambda)_{\Lambda \in \mathcal{S}}$ and $\Phi=(\Phi_A)_{A
\in \mathcal{S}}$ be set functions from $\mathcal{S}$ to
$\mathbb{R}$. Then

\begin{equation}\label{Moeb}
\forall \Lambda \in \mathcal{S}, H_\Lambda=\sum_{A \subset
\Lambda} \Phi_A \; \Longleftrightarrow \; \forall A \in
\mathcal{S}, \; \Phi_A= \sum_{B \subset A} (-1)^{|A \setminus B |}
H_B.
\end{equation}

\end{proposition}
We use it and propose then the
\begin{definition}A vacuum potential for a given specification $\gamma$ is the
potential defined for all $\sigma \in \Omega$ by
$\Phi^+_\emptyset=0$ and
\begin{equation}\label{vacuumpot}
\forall A \in \mathcal{S}, \; \Phi^+_A(\sigma) = - \sum_{B \subset A}
(-1)^{|A \setminus B |} \ln \frac{\gamma_B(\sigma | +)}{\gamma_B(+
| +)}=\sum_{B \subset A} (-1)^{|A \setminus B |} \ln \frac{f_B(+)}{f_B(\sigma_B +_{B^c})}.
\end{equation}
\end{definition}

\begin{lemma}[Convergence and consistency of the vacuum potential]\label{vacuumcons}
Let $\gamma$ be any quasilocal and non-null specification.
 Then  $\Phi^+=(\Phi^+_A)_{A \in \mathcal{S}}$ defined by (\ref{vacuumpot})
  is a vacuum potential with vacuum state $+ \in \Omega$, for any reference configuration $+ \in \Omega$. It is moreover convergent and its corresponding Gibbs specification $\gamma^{\Phi}$ coincides with $\gamma$.
\end{lemma}

It is obviously a potential. We consider any reference configuration  $+ \in \Omega$ and prove first that $\Phi^+$ satisfies the vacuum condition. It will be crucial to get consistency.
Consider $A \in \mathcal{F}$ and  $\sigma \in \Omega$ such
that there exists $i \in A$ where $\sigma_i=+_i$. One has
$$
\Phi^+_A(\sigma)=-\sum_{B \subset A} (-1)^{|A \setminus B|}\ln
\frac{\gamma_B(\sigma | +)}{\gamma_B(+ | +)}=\sum_{B \subset A}
(-1)^{|A \setminus B|} H_B^{\Phi^+,f} (\sigma).
$$
where by the Moebius formula, one has for all $B \in
\mathcal{S}$,

\begin{equation}\label{Hfree}
 H_B^{\Phi^+,f} (\sigma)=\sum_{A
\subset B} \Phi^+_A(\sigma)=-\ln \frac{\gamma_B(\sigma |
+)}{\gamma_B(+ | +)}=\ln \frac{f_B(+)}{f_\Lambda(B_\Lambda+_{\Lambda^c})}.
\end{equation}
Using  Equation (\ref{keybar}), one first gets
\begin{equation}\label{consistenciii}
\forall i \in B \subset A, \; H_B^{\Phi^+,f} (\sigma)= H_{B
\setminus i}^{\Phi^+,f} (\sigma).
\end{equation}
Indeed, by consistency property of the specification, one can
rewrites
$$
\frac{f_B(+)}{f_B(\sigma_B+_{B^c})}=\frac{f_B(+_{B \setminus i} +_i)}{f_{B}(\sigma_{B \setminus i}
+_i)} =\frac{f_{B \setminus i}(+)}{f_{B \setminus i}(\sigma_{B
\setminus i} +_i))}
$$
to eventually get (\ref{consistenciii}). Now define for any site $i \in S$ a partition of
$\mathcal{S}$ by $\mathcal{S}=\mathcal{S}_{A,i} \cup
\mathcal{S}_{A,i}^{c}$ with
$\mathcal{S}_{A,i}=\{V \in \mathcal{S}, V \subset A, V \ni i\}.$
An obvious bijection from $\mathcal{S}_{A,i}$ to
$\mathcal{S}_{A,i}^{c}$ links $B \in \mathcal{S}_{A,i}$ to
$B \backslash i \; \in \mathcal{S}_{A,i}^{c}$, so one gets
\begin{eqnarray*}
\Phi_A^+(\omega)&=& \sum_{B \in \mathcal{S}_{A,i}} (-1)^{\mid A
\backslash
  B \mid}H_B^{\Phi^+}(\omega) +\sum_{B \in \mathcal{S}_{A,i}^{c}} (-1)^{\mid A \backslash
  B \mid}H_B{\Phi^+}(\omega)\\
&=&\sum_{B \in \mathcal{S}_{A,i}} \Big[(-1)^{\mid A \backslash
  B \mid}H_B{\Phi^+}(\omega) +(-1)^{\mid A \backslash
 \{B \backslash i\}\mid}H_{B \backslash i}{\Phi^+}(\omega)\Big]\\
&=&\sum_{B \in \mathcal{S}_{A,i}}(-1)^{\mid A \backslash
  B \mid} \Big[H_B{\Phi^+}(\omega) - H_{B \backslash i}{\Phi^+}(\omega)\Big]\\
&=&0
\end{eqnarray*}
and $\Phi^+$ is indeed a vacuum potential associated to the
specification $\gamma$. This potential need not be convergent or
consistent with $\gamma$ in general, and we verify it now in this
non-null quasilocal case, that is we first need to prove that we
can always define
$$
\forall \sigma,\omega \in  \Omega, \; H_\Lambda^{\Phi^+} (\sigma |
\omega)=\sum_{A \cap \Lambda \neq \emptyset, A \in \mathcal{S}}
\Phi^+_A(\sigma_\Lambda \omega_{\Lambda^c})
$$
i.e.  to extend  the definition of the Hamiltonian with free b.c.
(\ref{Hfree}) to an Hamiltonian with any $\omega \in \Omega$ as a
boundary condition. It amounts to proving the convergence of the
potential, i.e. that for all $\sigma \in \Omega$,
$$
H_\Lambda^{\Phi^+} (\sigma):=\sum_{A \cap \Lambda \neq \emptyset,
A \in \mathcal{S}} \Phi^+_A(\sigma) \; < \; + \infty
$$
in the sense that the limit as $\Delta \uparrow \mathcal{S}$ of
the net $\Big(\sum_{A \cap \Lambda \neq \emptyset, A \subset \Delta}
\Phi^+_A(\sigma) \Big)_{\Delta \in \mathcal{S}}$ is finite.
Recall that we have already been able to define the Hamiltonian
with free b.c. as
$$
H_\Lambda^{\Phi^+,f}(\sigma)=\sum_{A \subset \Lambda} \Phi^+_A(\sigma) = \ln
\frac{f_\Lambda(+)}{f(\sigma_\Lambda +_{\Lambda^c})}.
$$

To prove  now that it is a convergent potential using the
quasilocality of the function $\omega \longmapsto
f_{\Lambda}(\sigma_{\Lambda}\omega_{\Lambda^c})$, we re-write
\begin{displaymath}
\sum_{A \cap \Lambda \neq \emptyset,A \subset
  \Delta}\Phi_A(\sigma)=\sum_{A \subset \Delta}\Phi_A(\sigma)-\sum_{A
  \subset \Delta \cap \Lambda^c}\Phi_A(\sigma).
\end{displaymath}
Using twice Moebius inversion formula (\ref{Moeb}), one obtains
\begin{displaymath}
\sum_{A \subset
  \Delta}\Phi_A(\sigma)=\ln{\frac{f_{\Delta}(+)}{f_{\Delta}(\sigma_{\Delta}+_{\Delta^c})}} \; \textrm{and} \; \sum_{A
  \subset \Delta \cap \Lambda^c}\Phi_A(\sigma)
  =\ln{\frac{f_{\Delta \cap \Lambda^c}(+)}{f_{\Delta \cap \Lambda^c}(\sigma_{\Delta \cap \Lambda^c}+_{\Delta^c \cup \Lambda})}}.
\end{displaymath}
By consistency (Lemma \ref{keylemma}), we get
$\ln{\frac{f_{\Delta}(+)}{f_{\Delta}(\sigma_{\Delta \cap
        \Lambda^c}+_{\Delta^c \cup \Lambda})}}$ for the second term, because,  on $\Delta \cap \Lambda^c$, the two involved configurations coincide
        and
      eventually
\begin{displaymath}
\sum_{A \cap \Lambda \neq \emptyset,A \subset
  \Delta}\Phi_A(\sigma)=\ln{\frac{f_{\Delta}(\sigma_{\Delta
  \cap \Lambda^c}+_{\Delta^c \cup \Lambda})}
{f_{\Delta}(\sigma_{\Delta}+_{\Delta^c})}}
\end{displaymath}
and using again Lemma (\ref{keylemma}), with the sets $(\Delta \cap
\Lambda,\Delta)$, and the configurations $(\sigma_{\Delta
  \cap \Lambda^c}+_{\Delta^c \cup \Lambda}, \sigma_{\Delta}+_{\Delta^c})$ which agree outside $\Delta \cap
  \Lambda$, we get
\begin{displaymath}
\sum_{A \cap \Lambda \neq \emptyset,A \subset
  \Delta}\Phi_A(\sigma)=\ln{\frac{f_{\Lambda \cap \Delta}(\sigma_{\Delta
  \cap  \Lambda^c}
  +_{\Delta^c \cup \Lambda})}{f_{\Lambda \cap \Delta}(\sigma_{\Delta}+_{\Delta^c})}}.
\end{displaymath}
Let $\Delta \uparrow S$ in the sense defined. For $\Delta
\supset \Lambda$, one gets
\begin{displaymath}
\sum_{A \cap \Lambda \neq \emptyset,A \subset
\Delta}\Phi_A(\sigma)=\ln{\frac{f_{\Lambda}(+_{\Lambda}\sigma_{\Delta
  \backslash \Lambda}
  +_{\Delta^c})}{f_{\Lambda}(\sigma_{\Delta}+_{\Delta^c})}}.
\end{displaymath}
Thus quasilocality implies that the potential $\Phi$ is
convergent and that
$$
\forall \sigma \in \Omega, \; H_\Lambda^{\Phi^+}(\sigma)=\sum_{A
\cap \Lambda \neq
  \emptyset, A \in \mathcal{S}}\Phi_A(\sigma)=-\ln \frac{\gamma_\Lambda(\sigma | \sigma)}{\gamma_\Lambda(+ | \sigma)}=\ln{\frac{f_{\Lambda}(+_{\Lambda}\sigma_{\Lambda^c})}{f_{\Lambda}(\sigma)}} < + \infty.
$$
Hence, every quasilocal and non-null specification $\gamma$ is consistent
with the convergent vacuum potential $\Phi^+$ whose the Hamiltonian with boundary
condition $\omega \in \Omega$ is defined for all $\sigma \in \Omega$
\be \label{Hambc}
\forall \omega \in \Omega, \; H_\Lambda^{\Phi^+}(\sigma | \omega)=- \ln \frac{\gamma_\Lambda (\sigma | \omega)}{\gamma_\Lambda (+ | \omega)}=\ln \frac{f_\Lambda(+_\Lambda \omega_{\Lambda^c})}{f_\Lambda(\sigma_\Lambda \omega_{\Lambda^c})} < + \infty.
\ee
This proves Lemma \ref{vacuumcons}: Any quasilocal and non-null specification is consistent with a convergent potential.\\

Unfortunately, this vacuum potential is not UAC in the sense of
(\ref{UAC}). To gain summability and absoluteness,
 Kozlov \cite{Ko} introduced a particular re-summation procedure by telescoping the terms of the Hamiltonian with
 free boundary conditions in large enough annuli to recover absoluteness, but carefully keeping consistency,
  to eventually get a potential $\Psi$ such that, for all $\sigma \in \Omega$
\begin{equation}\label{conskoz}
\sum_{A \subset \Lambda} \Psi_A(\sigma)=H_\Lambda^{\Phi^+,f}(\sigma)=\sum_{A \subset \Lambda} \Phi^+_A(\sigma) = \ln
\frac{f_\Lambda(+)}{f(\sigma_\Lambda +_{\Lambda^c})}
\end{equation}
with the extra summability property
\begin{equation}\label{uackoz}
\forall i \in S,\; \sum_{A \in \mathcal{S},A \ni i} \sup_\omega \Big | \Psi_A(\omega) \Big | < + \infty.
\end{equation}
We shall describe it formally following the pedagogical exposition of Fern{\'a}ndez \cite{F},
and describe a bit more explicitly the telescoping at the end of this proof.

In our  settings with a finite single-spin state space, non-nullness and quasilocality can be

reduced to site-characterizations that are very useful to get the
stronger summability around sites (\ref{uackoz}). Introduce, for
any site $i \in S$ and any cube $\Lambda_n$, the quantities
$$
m_i:=\inf_\omega f_{\{i\}} (\omega) = \inf_{\omega \in \Omega}\gamma_{\{i\}} (\omega | \omega)
$$
and
$$
g_i(n) = \sup_\omega \big| f_{\{i\}} (\omega_{\Lambda_n} +_{\Lambda_n^c}) - f_{\{i\}} (\omega) \big|
$$
By non-nullness, one has $m_i >0$ for all $i \in S$ and quasilocality reads
$$
\forall i \in S, \; g_i(n) \;
\mathop{\longrightarrow}\limits_{n\to\infty}\; 0.
$$

Starting from the expression of the Hamiltonians in terms of the
vacuum potential, which itself is expressed as the logarithm of
ratios of densities, Kozlov used the inequality
$$
\Big | \ln \frac{a}{b} \Big | \leq \frac{|a-b|}{\rm{min}(a,b)},\;
\forall a,b >0
$$
to get that for all $i \in S$
\be \label{upbound}
\sup_{\omega} \Big | \sum_{A \subset \Lambda_n, A \ni i} \Phi^+(\sigma) \Big | \leq \frac{g_i(n)}{m_i}
\ee
and in particular that
\be \label{upbound2}
\sup_\omega \Big | \sum_{A \subset \Lambda_n, A \ni i} \Phi^+(\sigma) - \sum_{A \subset \Lambda_{n-1}, A \ni i} \Phi^+(\sigma)\Big |\leq \frac{g_i(n) + g_i(n-1)}{m_i}.
\ee

Kozlov used then these bounds to reduce the lack of absolute
convergence by grouping terms of the vacuum potential within
intermediate annuli chosen large enough to exploit quasilocality.
To do so, the telescoping has to integrate larger boxes, i.e.
along subsequences of cubes $\Lambda_{n_k^i}, k \geq 1$ in
(\ref{upbound2})  chosen such that for any $i \in S$,
$$
\sum_{k \geq 1} g_i(n^i_k) <  \infty
$$
which is always possible because for any $i \in S$, the sequence $\big(g_i(n)\big)_{n \in \mathbb{N}}$ converges to zero.

For any $i \in S$, we consider then the subsequence of cubes $\Lambda_{n_k^i}$, centered in $i$, of radius $n_k^i$ such that
 the annuli $\Lambda_{n^i_k} \setminus \Lambda_{n_{k-1}^i}$ is thus large enough, as we shall see.
 This size will allow the use of the bounds (\ref{upbound}) for any $i \in S$ and to get the right summability properties,
 the telescoping is done by following the bonds along these cubes, adding at each steps the terms of the vacuum potential
 that correspond to bonds of the annulus, and that were not in the previous cubes. Define then, for any $i \in S$, any $k \geq 1$
$$
S_k^i= \big \{ B \subset \Lambda_{n_k}^i: B \ni i \big \} \setminus S_{k-1}^i
$$
with  $S^i_0=\{i\}$, and introduce the potential\footnote{This potential is not yet the Kozlov potential, so we write it $\tilde{\Psi}$, because the resummation uses several times terms involving sites $i$.} $\tilde{\Psi}$ defined by

\begin{displaymath}
\tilde{\Psi}_A(\sigma)=\left\{
\begin{array}{lll}
\;\tilde{\Psi}_A(\sigma)= \sum_{B \in S_k} \Phi_B^+(\sigma) \;\; \;  \;  \; \textrm{if} \;  \;  \; \; \; A=\Lambda_{n_k}^i \; \textrm{for some} \; k \geq 1, {\rm some} \; i \in S.\\
\\
\; 0 \; \; \; \; \; \; \; \;  \; \; \; \;  \;  \; \; \textrm{otherwise.}
\end{array} \right.
\end{displaymath}
By (\ref{upbound}), one has
$$
\sup_\omega \big| \tilde{\Psi}_A(\sigma) \big| = \sup_\omega \Big | \sum_{B \ni i, B \subset \Lambda_{n_k^i}} \tilde{\Phi}_B^+(\sigma) - \sum_{B \ni i, B \subset \Lambda_{n_{k-1}^i}} \tilde{\Phi}_B^+(\sigma) \Big | \leq \frac{g_i(n_k^i) + g_i(n_{k-1}^i)}{m_i}.
$$
in such a way that one has the right summability property at the site $i$:
$$
\sum_{A \ni i} \sup_\omega | \tilde{\Psi}_A^i (\sigma) |  \leq \frac{2}{m_i} \cdot \sum_{k \geq 1} g_i(n_k^i) < + \infty.
$$

Nevertheless, we need to do the telescoping more carefully to keep
the consistency, in general lost in the procedure above: For a
given $B \in \mathcal{S}$, the same vacuum interaction $\Phi_B$
could have been  used more than once. To avoid it, one has to find
a way of grouping terms of the vacuum interaction without using
the terms already used, i.e. one has to run the sequence of cubes
by using any finite set $B$ of bonds only once. To do so,
Fern{\'a}ndez \cite{F} proposed the following presentation of
Kozlov`s potential, now denoted by $\Psi$. The sites of the
lattice will be now lexicographically ordered and still
generically denoted by $i$. For any site, one replaces the
previous subsequence of cubes $\Lambda_{n_k^i}$ by rectangles
around it that do not incorporate $B$'s (or i's) already
considered. Hence, one defines for each $i=1,2, \dots$, a sequence
$(L_k^i)_{k \geq 1}$ defined for all $i,k \geq 1$ by
$$
L_k^i= \big \{ j \in S: i \leq j \leq r_k^i \big \}
$$
where the diameters $r_k^i$ are chosen such that $n_k^i={\rm
diam}(L_k^i)=r_k^i-i$ in order to keep the same large enough
sequence of annuli. These groups of bonds will be the only one
involved in the potential and to perform a correct re-summation
procedure,  one defines for any site $i \in S$, a family of {\em
disjoints} subsets of $\mathcal{S}$ containing $i$ by
$S_0^i=\{i\}$ and
$$
S_k^i= \big \{ B \subset L_k^i: B \ni i \big \} \setminus S_{k-1}^i
$$
in such a way that $B \in \cup_{j=1}^i \cup_{k \geq 1} S_j^i$ and  any $B$ containing $i$ is uniquely contained in  one of them.
By this procedure, any set of bonds B  is considered only once and we get the Kozlov potential $\Psi$ defined by
\begin{equation}\label{kozlovpot}
\Psi_A(\sigma)=\left\{
\begin{array}{lll}
\;\Psi_A(\sigma)= \sum_{B \in S_k^i} \Phi_B^+(\sigma) \; \textrm{if} \; A=L_{k}^{i} \; \textrm{for some} \; (i,k), i \in S, k \geq 1\\
\\
\; 0 \; \; \; \; \; \; \;\; \textrm{otherwise.}
\end{array} \right.
\end{equation}
yielding the following
\begin{lemma}
The Kozlov`s potential $\Psi$ defined by (\ref{kozlovpot}) is a UAC potential consistent with the non-null
and quasilocal specification $\gamma$, and thus any quasilocal measure $\mu \in \mathcal{G}(\gamma)$ is a Gibbs measure.
\end{lemma}

Consistency holds because the careful procedure yields the same
Hamiltonian with free boundary conditions for the vacuum and
Kozlov potentials, and convergence is due to the choice of the
subsequences $n_k^i$:
\begin{eqnarray*}
\forall i \in S, \sum_{A \in \mathcal{S},
 A \ni i} \sup_\omega \Big | \Psi_A(\omega) \Big | &\leq& \sum_{j=1}^i \sum_{k \geq 1} \sup_\omega \big | \Psi_{L_k^i}(\omega)
 \big |\leq \sum_{j=1}^i \frac{2}{m_i} \; \sum_{k \geq 1} g_i(n_k^i) < \infty.
\end{eqnarray*}
This proves the lemma and  the Gibbs
representation theorem \ref{Gibbsrepthm}.

\begin{remark}[Telescoping procedure] \label{telescop}
{\em To get an idea of the type of telescoping that has to be done,
we informally detail it starting from the expression (\ref{Hfree})
of the Hamiltonian with free boundary condition, in order to see how consistency
is important to get it. It is also the way the procedure is done in an adaptation to generalized Gibbs
measures in Chapter 5 to get {\em weakly Gibbsian measures}, see \cite{Maes1,MRM,Maes3,Maes4}, and we
also use a similar procedure  in \cite{KLNR} to get a variational principle for translation-invariant quasilocal measures.

Let us start from the Hamiltonian with free boundary condition for some $\Lambda$ containing the origin and
assume $\Lambda$ to be a cube $\Lambda_{n_l}$ of the subsequence already taken, and write $L_{n_l}$ for
the annulus $\Lambda_{n_l} \setminus \Lambda_{n_{l-1}}$ . One has by consistency and  (\ref{Hfree})
$$
H_{\Lambda_{n_l}}^{\Phi^+,f}(\sigma)= \ln
\frac{f_{\Lambda_{n_l}}(+)}{f(\sigma_{\Lambda_{n_l}} +_{\Lambda_{n_l}^c})}.
$$
The idea now is to telescope this term by incorporating terms
corresponding to so-called relative energies by flipping the spin
in the annulus only
$$
H_{\Lambda_{n_l}}^{\Phi^+,f}(\sigma)= \ln
\frac{f_{\Lambda_{n_l}}(+)}{f_{\Lambda_{n_l}}(+_{\Lambda_{n_{l-1}}}\sigma_{L_{n_l-1}} +_{\Lambda_{n_l}^c})} \cdot \frac{f_{\Lambda_{n_l}}(+_{\Lambda_{n_{l-1}}}\sigma_{L_{n_{l-1}}} +_{\Lambda_{n_l}^c})}{f_{\Lambda_{n_l}}(+_{\Lambda_{n_{l-2}}}\sigma_{L_{n_{l-2}}} +_{\Lambda_{n_{l-1}}^c})}.
$$
Doing it for any $k=1 \dots l$, one gets
$$
H_{\Lambda_{n_l}}^{\Phi^+,f}(\sigma)= \sum_{k=1}^l \ln
\frac{f_{\Lambda_{n_l}}(+_{\Lambda_{n_k}}\sigma_{L_{n_k+1}} +_{\Lambda_{n_{k+1}}^c})}
{f_{\Lambda_{n_l}}(+_{\Lambda_{n_{k-1}}}\sigma_{L_{n_k}} +_{\Lambda_{n_k}^c})} \cdot \frac{f_{\Lambda_{n_l}}(+_{\Lambda_{n_{k-1}}}\sigma_{L_{n_k}} +_{\Lambda_{n_k}^c})}{f_{\Lambda_{n_l}}(+_{\Lambda_{n_{k-2}}}\sigma_{L_{n_{k-1}}} +_{\Lambda_{n_{k-1}}^c})}.
$$
Now using consistency via our
key lemma \ref{keylemma}, one can replace the densities $f_\Lambda$  by densities corresponding
to the proper annulus in order to get proper measurability conditions for the potential; To see this, one rewrites thus
$$
H_{\Lambda_{n_l}}^{\Phi^+,f}(\sigma)= \sum_{k=1}^l \ln
\frac{f_{\Lambda_{n_k}}(+_{\Lambda_{n_k}}\sigma_{L_{n_k+1}} +_{\Lambda_{n_{k+1}}^c})}{f_{\Lambda_{n_k}}(+_{\Lambda_{n_{k-1}}}\sigma_{L_{n_k}} +_{\Lambda_{n_k}^c})} \cdot \frac{f_{\Lambda_{n_{k-1}}}(+_{\Lambda_{n_{k-1}}}\sigma_{L_{n_k}} +_{\Lambda_{n_k}^c})}{f_{\Lambda_{n_{k-1}}}(+_{\Lambda_{n_{k-2}}}\sigma_{L_{n_{k-1}}} +_{\Lambda_{n_{k-1}}^c})}
$$
to eventually get a potential of the
form $\Psi_{L_{n_k}}=\Phi_{\Lambda_{n_k}} - \Phi_{\Lambda_{n_{k-1}}}$ for which
consistency holds together with the required convergence property. We shall use such a  procedure,
 using the expression of (\ref{keylemma}) in terms of $\gamma$ instead of $f$ in Chapter 4.}
\end{remark}

\begin{remark}[Translation-invariance and Kozlov vs. Sullivan압 results]\label{rksull}
 {\em The procedure due to Kozlov to introduce its UAC potential does not yield a
 translation-invariant one, as explicitly seen in the site-dependent way of
  re-ordering terms of the Hamiltonian with free boundary condition. It is
  nevertheless possible to consider larger rectangles partitioning $\Lambda$
  similarly for all site considered but they have to be larger and require a condition a bit stronger
   than quasilocality. It is an open question whether this condition is technical or not,
    but Sullivan \cite{Su}  has observed that the vacuum potential,
     in addition to be convergent and translation-invariant, is {\em relatively uniformly convergent} in the sense that the series
\be \label{relatunifcon} \sum_{A \cap \Lambda \neq \emptyset,A \in
\mathcal{S}} \big | \Phi^+_A(\omega) -
\Phi^+_A(\omega'_{\Lambda}\omega_{\Lambda^+}) \big | < + \infty
\ee are uniformly convergent in $(\omega,\omega')$, and also that
this was enough to get quasilocality of the Gibbs specification.
For the vacuum potential itself and $\omega'=+$, this implies the
usual uniform convergence of the vacuum potential. Reciprocally,
any quasilocal specification has a relatively uniformly convergent
potential and we shall see that it is also interesting to derive
thermodynamical properties for these measures, although this
convergence is too weak to get all the flavor of the Gibbsian
theory (see again the discussion in \cite{VEFS}).}
\end{remark}

\chapter{Equilibrium approach}

We present now an alternative approach to  describe equilibrium states at
infinite-volume that has a more physical flavor and which will eventually appear to be partially equivalent to
the DLR construction presented above. Inspired by the second law of thermodynamics described at finite volume in
 the introduction of Chapter 3, this so-called {\em equilibrium approach to Gibbs measures} provides
  thermodynamic functions at infinite-volume and yields afterwards infinite-volume counterparts of the
   second law of thermodynamics in terms of zero relative entropy or in terms of minimization of free energy.
   This approach is restricted to a {\em translation-invariant} framework, mostly because it is mainly untractable
    otherwise\footnote{Except in a few situations when e.g. {\em periodic} boundary conditions are considered \cite{LPS,Rue}.},
     and we shall characterize the translation-invariant equilibrium states  of a given system in terms of
     {\em variational principles}, either specification-dependent or specification-independent depending on
      the choice made to  characterize of the second law, as we shall see. We thus  describe, in a rather
      general framework that will be useful for generalized Gibbs measures next chapter, how this approach is
      in some sense equivalent to the DLR approach  restricted to translation-invariant measures, and describe
       the general  proof given in  \cite{KLNR} of the recent result that a (specification-dependent) variational
       principle holds for translation-invariant quasilocal specifications in general.\\

In all this chapter, we also restrict the  infinite-volume limit procedure by mostly considering the
 limit $\Lambda \uparrow \mathcal{S}$ along  sequences of cubes $(\Lambda_n)_{n \in \mathbb{N}}$, or at least
 sequences s.t. the ratio (surface boundary)/(volume) $\frac{|\partial \Lambda |}{|\Lambda |} \to 0$, within
  the so-called {\em thermodynamic limit}, although the latter is a bit more general \cite{VEFS,Is}.
   For this reason, we focus on the $d$-dimensional regular lattice $S=\mathbb{Z}^d$, because this
   thermodynamic limit does not hold for cubes on trees\footnote{The ratio (surface boundary)/(volume)
    do not vanish in the limit,
    see also \cite{BFS}.}, which are the other lattices sometimes considered in these notes. We also incorporate the temperature
    in the potential when dealing with Gibbs specifications and measures, or equivalently here we assume  $\beta=1$. At the end
    of the chapter, we briefly mention related large deviation properties for the considered measures and introduce
    thereafter another way to consider Gibbs measures as equilibrium states of the system, defining them as  invariant
    measures of Markov processes on the configuration space, and illustrate this notion by the so-called {\em stochastic Ising models}, that will also be discussed in the generalized Gibbs framework in Chapter 5.

\section{Thermodynamic properties}
\subsection{Thermodynamic functions}
We have already introduced in Chapter 3 the \emph{relative
entropy} at finite volume $\Lambda \in \mathcal{S}$ of $\mu$
relative to $\nu$ for two translation-invariant measures $\mu,\nu
\in \mathcal{M}_{1,\rm{inv}}^+(\Omega)$, defined to be
\begin{equation}\label{eq:rex}
\mathcal{H}_{\Lambda}(\mu | \nu)=\int_{\Omega}  \Big(\frac{d\mu_{\Lambda}}{d\nu_{\Lambda}} \Big)
\cdot \log \Big(\frac{d\mu_{\Lambda}}{d\nu_{\Lambda}}\Big) \; d \nu
\end{equation}
when the projection $\mu_\Lambda$ of $\mu$ on
$(\Omega_\Lambda,\mathcal{F}_\Lambda)$ is absolutely continuous
w.r.t. the projection $\nu_\Lambda$ of $\nu$, and to be
$\mathcal{H}_{\Lambda}(\mu | \nu)= + \infty$ otherwise. To avoid
trivial divergences at infinite-volume, one  considers quantities
{\em per unit of volume}, writes for any $n \in \mathbb{N}$
\begin{equation}\label{finiteentrop}
h_n(\mu | \nu):=\frac{1}{|\la_n|}\sum_{\sigma_{\la_n} \in \Omega_{\Lambda_n}}
\mu(\sigma_{\la_n}) \cdot \log
\frac{\mu(\sigma_{\la_n})}{\nu(\sigma_{\la_n})}.
\end{equation}
 and introduces  the
\emph{relative entropy density} of $\mu$ relative to $\nu$ to be the
limit
\begin{equation}\label{eq:red}
h(\mu | \nu)\;=\;\lim_{n \to \infty} h_n(\mu | \nu)
\end{equation}
provided it exists.  The limit is known to exist for any arbitrary
$\mu\in \mathcal{M}_{1,\rm{inv}}^{+}(\Omega)$ when
$\nu\in\mathcal{M}_{1,\rm{inv}}^{+}(\Omega)$ is a Gibbs measure
(for a UAC potential) and, more generally, if $\nu$ is
asymptotically decoupled\footnote{These  are measures  introduced
by Pfister \cite{P} to state general large deviation principles,
see next chapter.}. We extend this result next section for
general {\em translation-invariant quasilocal measures} in Theorem
\ref{chp:thm7}, whose proof follows \cite{KLNR}. We also recall
(see e.g. \cite{Bil2}) that for $\mu \in
\mathcal{M}_{1,\rm{inv}}^+(\Omega)$, the entropy per unit of
volume \be \label{limKSn}
h_n(\mu)=-\frac{1}{|\la_n|}\sum_{\sigma_{\la_n}}
\mu(\sigma_{\la_n})\log \mu(\sigma_{\la_n}) \ee has a well-defined
limit \be \label{KS} h(\mu):=- \lim _{n \to \infty}
\frac{1}{|\la_n|}\sum_{\sigma_{\la_n}} \mu(\sigma_{\la_n})\log
\mu(\sigma_{\la_n}) \ee
 called the {\em
Kolmogorov-Sinai entropy} of $\mu$.\\

For $\nu \in
\mathcal{M}_{1,{\rm inv}}^+(\Omega)$ and $f \in \mathcal{F}$  bounded, the \emph{pressure}  for
$f$ relative to $\nu$ is defined as
\begin{equation} \label{Pressuregen}
p(f | \nu)= \lim_{n \to \infty}\,\frac{1}{|\Lambda_n |^{d}}\, \log\int\exp\Bigl(\sum_{x
\in \Lambda_n}\tau_xf \Bigr)\,d\nu
\end{equation}
whenever this limit exists.  This limit exists, for every quasilocal
function $f$, if  $\nu$ is Gibbsian \cite{VEFS,Ge} or asymptotically decoupled
\cite{P}.

When dealing with  a {\em translation-invariant}  potential
$\Phi$, the particular choice
$$
f=f_\Phi:=\sum_{A \ni 0} \frac{1}{|A|} \cdot \Phi_A
$$
and with the a priori product measure $\rho$ as reference measure,
it connects with the usual pressure in the case of a lattice gas
and is more  generally related  with the free energy of a system
obtained from a partition function with boundary condition
$\omega$ or with free boundary condition. At finite volume
$\Lambda_n$ they are respectively defined to be
$$
P^f_{\Lambda_n}(\Phi)=\frac{1}{|\Lambda_n|} \; \ln Z_{\Lambda_n}^{\Phi, f} \; \;\; {\rm and} \; \ \ P^\omega_{\Lambda_n}(\Phi)=\frac{1}{|\Lambda_n|} \; \ln Z_{\Lambda_n}^{\Phi} (\omega).
$$
When the limit  exists, it captures many information of the
particle system\footnote{This is probably the most important part
of mathematical statistical mechanics that we do not develop in
this course, see again  \cite{VEFS,Is,Ge}.} and for UAC potentials
it turns out to be independent of the boundary condition:
\begin{theorem}[Pressure of a UAC potential  \cite{Is}]\label{pressureUAC}
Let $\Phi$ be a U.A.C. translation-invariant potential. Then,  the
limits
\begin{equation} \label{pressUAC}
\lim_{n \to \infty} \; \frac{1}{|\Lambda_n|} \; \ln Z_{\Lambda_n}^{\Phi,\rm{f}} \; \; \; \; {\rm and} \ \ \ \
\lim_{n \to \infty} \; \frac{1}{|\Lambda_n|} \; \ln Z_{\Lambda_n}^{\Phi}(\omega)
\end{equation}
exist, coincide for all $\omega \in \Omega$ and define the {\em pressure of the potential} $\Phi$:
\be \label{PressureUAC2}
P(\Phi):=\lim_{n \to \infty} P^f_{\Lambda_n}(\Phi)=\lim_{n \to \infty} P^\omega_{\Lambda_n}(\Phi)
\ee
exists and is thus independent of the boundary condition $\omega \in \Omega$.
\end{theorem}

We recall briefly the philosophy of  the proof of Israel \cite{Is}, because we  extend it
in this chapter to deal with the vacuum potential of a translation-invariant  quasilocal
 measure, which is translation-invariant but not UAC, using its {\em relative uniform convergence}, introduced
  in Remark \ref{rksull}, following Sullivan \cite{Su}. The proof focuses first on  finite-range potentials
  and using the finiteness of the range, and thus the independence of  spins from sites far enough, one introduces a
  partition of the volume $\Lambda$ into cubes and corridors, the width of the latter being at least the range of the potential,
  to eventually get a factorization of the partition function up to some boundary terms that are negligible in the thermodynamic
  limit\footnote{This is the reason why this restriction on the ratio (surface boundary)/(volume) is made in this approach.}.
   This factorization leads to sub-additivity of the logarithm of  the partition function, which in turns implies the existence
   of the pressure. This result  is then extended to  general UAC potentials using the density of the finite-range potential in
   the Banach space of translation-invariant UAC potentials \cite{Is}, and thereafter the strong UAC  convergence to get the
   independence with the boundary condition.\\

 We
also introduce the {\em $\nu$-specific energy of a reference configuration}
$+ \in \Omega$:
\be \label{nuispecificenergy}
e_{\nu}^+:=-\lim_{\la \uparrow \Z^d} \frac{1}{|\la|}
 \log \nu(+_{\la})
\ee whenever it exists. We prove its existence for
translation-invariant quasilocal measures in this chapter in order
to get a general variational principle for translation-invariant
quasilocal measures. We emphasize the fact that it is not properly
speaking\footnote{For a vacuum potential with vacuum state $+$,
and a quasilocal specification, it even coincides with the
pressure, see the proof of of Lemma \ref{chp:lem2}.} the
infinite-volume counterpart of the energies of Chapter 3.

\subsection{Variational principles}
To translate into a proper mathematical framework the second law
of thermodynamics in the vein of the finite-volume description
given in Chapter 3, we  distinguish between thermodynamical
(specification-independent) and statistical mechanical
(specification-dependent) variational principles.
\begin{definition}[Thermodynamic variational principle]\label{chp:def3}
$\nu \in \mathcal{M}_{1,\rm{inv}}^+(\Omega)$ is said to satisfy a
(thermodynamical) variational principle if the  relative entropy $h(\mu | \nu)$ and
the pressure $p(f | \nu)$ exist for all $\mu \in
\mathcal{M}_{1,\rm{inv}}^+(\Omega)$ and all $f \in
\mathcal{F}_{{\rm qloc}}$, and are conjugate convex functions
in the sense that
\begin{equation}\label{eq:r2}
\forall f \in \mathcal{F}_{{\rm qloc}}, \; \; \; p(f | \nu) \;= \; \sup_{\mu \in \mathcal{M}_{1,{\rm
inv}}^{+}(\Omega)} \Bigl[\mu(f) - h(\mu | \nu) \Bigr].
\end{equation}
\begin{equation}\label{eq:r3}
\forall \mu \in \mathcal{M}_{1,\rm{inv}}^{+}(\Omega), \; \; \; h(\mu | \nu) \;= \;\sup_{f \in \mathcal{F}_{{\rm
qloc}}}\Bigl[\mu(f) - p(f | \nu)\Bigr].
\end{equation}
\end{definition}

This is the infinite-volume counterpart of the formulation of the
second law of thermodynamics in terms of the  minimization of free
energy, which coincides with the pressure here. Gibbs measures
satisfy this specification-independent principle and Pfister
\cite{P} has  extended its validity to the larger class of
asymptotically decoupled measures, described next chapter within
the generalized Gibbsian framework. These conjugate convex
functions are also very important to study large deviation
properties of DLR measures. We do not focus much on this type of
variational principle in these lectures, and prefer focusing on
the other formulation of the second law of thermodynamics, in
terms of zero relative entropy.

\begin{definition}[Variational principle relative to a specification]\label{chp:defr}
Consider a specification $\gamma$  and $\nu \in \mathcal{G}_{\rm
inv}(\gamma)$.  A variational principle occurs for
$(\nu,\gamma)$ iff
\begin{equation}\label{eq:r4}
\forall \mu\in\mathcal{M}_{1,\rm{inv}}^{+}(\Omega), \; h(\mu|\nu) =0 \ \Longleftrightarrow\ \mu\in \mathcal{G}_{\rm inv}(\gamma)\;.
\end{equation}
\end{definition}

Hence, this property  has to be related to the  formulation of the
second law of thermodynamics in terms of  zero relative entropy,
which at finite volume implies  equality of measures, as explained
in the beginning of Chapter 3. At infinite-volume nevertheless,
two different measures could have zero relative entropy, but when
the reference measure has some nice locality
properties\footnote{Like e.g. being Gibbs. To get a
counterexample, i.e. two measures having zero relative entropy
without having much in common, consider the voter model in
dimension $d=3$, with its Dirac invariant measure as a reference
measure.}, the other measure, although possibly different, should
share the same system of conditional probabilities, identifying
the corresponding measures as equilibrium states of the system.
This result is well known for Gibbs measures consistent with a
translation-invariant UAC potential \cite{Ge} and has thus been
extended recently to translation-invariant quasilocal DLR measures
\cite{KLNR}. We describe  now  this result under a more general
form that will be useful for its extension to non-Gibbsian and
non-quasilocal measures in Chapter 5.

\section{Topological criterion for variational principles}

In this section, we consider specifications in the general
framework of Chapter 2, non necessarily Gibbsian or quasilocal,
and use some specific concentration properties on some points of
(sometimes partial) continuity of the specification. We introduce
different  sets of points of continuity. First, the {\em set
$\Omega_\gamma$ of good configurations of $\gamma$} is the set of
points of continuity, i.e.
\begin{equation}\label{Omegagam}
\Omega_\gamma = \Big\{ \omega \in \Omega: \forall \Lambda \in \mathcal{S},\; \forall f \in \mathcal{F}_{\rm{loc}}, \;
\lim_{n \to \infty} \; \sup_{\sigma \in \Omega} \big| \gamma_\Lambda f(\omega_{\Lambda_n}
\sigma_{\Lambda_n^c}) - \gamma_\Lambda f(\omega) \big|=0 \Big\}.
\end{equation}

We also consider points of continuity in some specific direction, say  $+ \in \Omega$.
Introduce, for all $n \in \mathbb{N}$, the {\em truncated kernels} $\gamma_{\Lambda}^{n,+}$
defined for all $f \in \mathcal{F}_{\rm{loc}}$ and all $\Lambda \in \mathcal{F}$ by

\begin{equation}\label{trunckern}
\forall \omega \in \Omega,\; \gamma^{n,+}_\Lambda f (\omega)=
\gamma_\Lambda f(\omega_{\Lambda_n}\,+_{\Lambda_n^c}).
\end{equation}
The set of points\footnote{Remark that a function can be
continuous in any direction $+$ without being continuous, see
\cite{FLNR2}.} that are {\em continuous in the $+$-direction} is
then defined to be

\begin{equation}\label{Omegagam+}
\Omega_\gamma^+ = \Big\{ \omega \in \Omega: \forall \Lambda \in \mathcal{S},\;\forall f \in \mathcal{F}_{\rm{loc}}, \; \lim_{n \to \infty} \; \gamma_\Lambda^{n,+} f (\omega) = \gamma_\Lambda f (\omega) \Big\}.
\end{equation}

For  presumably\footnote{See the telescoping procedure next
section and a discussion in Chapter 5.} technical reasons due to
some telescoping procedure, we introduce also the set of
configurations $\sigma$  for which there is continuity at some
particular {\em concatenated configuration} $\sigma^+ \in \Omega$
defined for all $\sigma \in \Omega$ by
\begin{equation}\label{concat}
\sigma^+_i=\sigma_i \; \; \; {\rm if} \; \ \ i \geq 0, \; {\rm and} \ \sigma_i=+_i \; {\rm otherwise}.
\end{equation}
We denote this set by
\begin{equation}\label{Omegagam+inf}
\Omega_\gamma^{<0} = \Big\{ \sigma \in \Omega : \sigma^+ \in
\Omega_\gamma \Big\}.
\end{equation}
Due to the definition  of these sets by a limiting procedure,
involved in any continuity-type property, one can prove  that these sets are tail-measurable and such that
$$
\Omega_\gamma \subset \Omega_\gamma^+ \subset \Omega_\gamma^{<0} \in \mathcal{F}_\infty.
$$
\subsection{Second part of the variational principle: General criterion}

Getting consistency from zero relative entropy is seen as the
''easiest part'', usually called the {\em second part}. The result and its proof are standard, but we give a slightly
more general version of both, that will be useful also for
non-quasilocal measures in Chapter 5. What is actually really needed is  some  weaker {\em continuity in the
+-direction} for some reference configuration $+ \in \Omega$, as we see now.

\begin{theorem}\label{r.teo1}\cite{FLNR,FLNR2}
Let $\gamma$ be a  specification that is quasilocal in the
direction $+\in\Omega$, and $\nu\in\mathcal{G}_{\rm
inv}(\gamma)$. If $\mu\in\mathcal{M}_{1}^+(\Omega)$ is such that $h(\mu|\nu)=0$, then
\begin{equation}\label{eq:r.10}
\mu\in \mathcal{G}(\gamma) \; \Longleftrightarrow \;
\nu\Bigl[\, g_{\Lambda_n\setminus\Lambda} \cdot \Bigl(\gamma^{n,
+}_\Lambda f - \gamma_\Lambda f\Bigr) \Bigr]
\;\mathop{\longrightarrow}\limits_{n\to\infty}\;0
\end{equation}
for all $\Lambda\in\mathcal{S}$ and $f\in\mathcal{F}_{\rm loc}$,
 where $g_{\Lambda_n \setminus \Lambda}:=\frac{d\mu_{\Lambda_n\setminus\Lambda}}{d\nu_{\Lambda_n\setminus\Lambda}}$
 provided it exists.
\end{theorem}

Thus, to get information about consistency from zero relative entropy requires that the
concentration properties of the density of $\mu_{\Lambda_n \setminus \Lambda}$ w.r.t $\nu_{\Lambda_n \setminus \Lambda}$ to
beat asymptotic divergence due to the lack of continuity of $\gamma$. When $\gamma$ is quasilocal, this lack of continuity
 never exists ($\Omega_\gamma=\Omega$) and is thus always beaten, yielding  the standard proof of the second part of the
 variational principle for translation-invariant quasilocal specifications.\\

{\bf Proof of Theorem \ref{r.teo1}}:

 It comes from \cite{FLNR,FLNR2} and is an adaptation of the standard proof \cite{Ge,P} in  the quasilocal
 case, where the criterion (\ref{eq:r.10}) is trivially valid.

  When  $h(\mu|\nu)=0$ holds, the latter relative entropy is in particular well-defined and as a consequence, for $n$ sufficiently
large, the   $\mathcal{F}_{\Lambda_n}$-measurable density $g_{\Lambda_n}:=d\mu_{\Lambda_n}/d\nu_{\Lambda_n}$ exists.
Fix $f \in \mathcal{F}_{\rm{loc}},\; \Lambda\in\mathcal{S}$ and pick $n$ big enough to get both
$\Lambda_n\supset \Lambda$ and the existence of  $g_{\Lambda_n}$. To prove that $\mu \in \mathcal{G}(\gamma)$,
we prove that $\mu \gamma_\Lambda[f] = \mu[f]$, and approximate first $\gamma$ by the truncated kernel (\ref{trunckern}), writing
$$
\mu \gamma_\Lambda[f] = \mu \gamma^{n,+}_\Lambda[f] + A_n
$$
where $A_n= \mu \big[\gamma_\Lambda f - \gamma^{n,+}_\Lambda f \big]$ goes to zero when $n$ goes to
infinity by continuity in the $+$-direction. Now we can use consistency,
 the $\mathcal{F}_{\Lambda_n \setminus \Lambda}$-measurability of the truncated kernel
 $\gamma^{n,+}_\Lambda[f]$ and  the $\mathcal{F}_{\Delta}$-measurability of any density $g_{\Delta}$ to rewrite
$$
\mu  \gamma^{n,+}_\Lambda[f]=\mu\big[  \gamma^{n,+}_\Lambda f \big] = \nu \big [g_{\Lambda_n \setminus \Lambda}  \cdot \gamma^{n,+}_\Lambda f \big] = \nu \big [g_{\Lambda_n \setminus \Lambda}  \cdot \gamma_\Lambda f \big] + B_n
$$
where $$B_n=\nu \Big[g_{\Lambda_n \setminus \Lambda}  \cdot (
\gamma^{n,+}_\Lambda f -  \gamma_\Lambda f )\Big]$$ is controlled
by the criterion (\ref{eq:r.10}). Now, by
$\mathcal{F}_{\Lambda^c}$-measurability of $g_{\Lambda_n \setminus
\Lambda}$ and consistency of $\nu$, one rewrites
$$
\nu \big [g_{\Lambda_n \setminus \Lambda}  \cdot \gamma_\Lambda f \big]= \nu \Big[ \gamma_\Lambda (g_{\Lambda_n \setminus \Lambda}  \cdot f) \Big]= \nu \Big[g_{\Lambda_n \setminus \Lambda} \cdot f \Big]
$$
which has to be compared with
$$
\mu[f]=\nu \big[g_{\Lambda_n} \cdot f \big]= \nu \big[g_{\Lambda_n \setminus \Lambda} \cdot f \big] + C_n
$$
with
$$
C_n=\nu \Big[ (g_{\Lambda_n \setminus \Lambda} - g_{\Lambda_n}) \cdot f \Big].
$$
We can evaluate  $C_n$ using the following in equality due to
Csisz\'ar \cite{C}
\begin{eqnarray*}
|C_n| =\Bigl|\nu\Bigl[(g_{\Lambda_n}  -
g_{\Lambda_n\setminus\Lambda})\, f\Bigr]\Bigr|
\leq \Bigl(2 \cdot  \sup_{\omega \in \Omega} |f(\omega)| \cdot \big(H_{\Lambda_n}(\mu | \nu) -
H_{\Lambda_n \setminus\Lambda}(\mu | \nu)\big) \Big)^{1/2}
\end{eqnarray*}
and use that the hypothesis of zero relative entropy implies that
finite-volume  relative entropy cannot grow faster that the
volume, and thus that the same is true for densities. Thus one has
$\mu \gamma_\Lambda [f]=\mu[f]$ if and only if
$B_n\;\mathop{\longrightarrow}\limits_{n\to\infty}\;0$, which
proves the lemma.

This criterion is a way to express that getting zero relative
entropy is meaningful only when the measures share some locality
properties. Without properties of that type, things could be
different, as shown by an example of $\cite{Xu}$. This criterion
has been upgraded in \cite{VEV} via the following theorem, that
will is also useful for the generalized Gibbs measures.

\begin{theorem}\cite{VEV}
Let $\mu \in \mathcal{G}(\gamma)$ and $\nu \in
\mathcal{M}_1^+(\Omega)$ such that $h(\mu|\nu)=0$. If, for all
$\sigma \in \Omega$, for $\nu$-a.e. $\omega$,
$$
\mu(\sigma_\Lambda|\omega_{\Lambda_n \setminus \Lambda}) \;\mathop{\longrightarrow}\limits_{n\to\infty}\;
\gamma_\Lambda(\sigma | \omega)
$$
then $\nu \in \mathcal{G}(\gamma)$.
\end{theorem}
\subsection{First part of the variational principle: General criterion}

In the usual theory of Gibbs measures, this part of the
variational principle, i.e getting zero relative entropy from
consistency,  is known when the latter  holds with a
translation-invariant specification defined via  a {\em
translation-invariant UAC} potential,
 and goes via existence and boundary condition independence of pressure (see \cite{Ge}).
Since for a general translation-invariant quasilocal specification
$\gamma$ we cannot rely on the existence of such a
translation-invariant potential, we shall use  the weaker property
of relative uniform convergence of the (translation-invariant)
vacuum potential which can be associated to the quasilocal
$\gamma$, as discussed in Remark \ref{rksull}, is enough to obtain
zero relative entropy. This result is  a consequence of more
general results on generalized Gibbs
measures developed in \cite{KLNR} that will be useful in the next chapter. \\

We consider a translation-invariant specification $\gamma$ and  a
probability measure  $\nu \in \mathcal{G}_{\rm{inv}}(\gamma)$.
\begin{theorem}\label{chp:thm2}\cite{KLNR}
If $\mu \in \mathcal{M}_{1,\rm{inv}}^+(\Omega)$ is such that $\mu(\Omega_\gamma^{<0})=1$ and $e_{\nu}^+$ exists, then
\begin{enumerate}
\item $h(\mu|\nu)$ exists and is given by
\begin{equation}\label{chp:exist}
h(\mu|\nu)= - h(\mu) + e_{\nu}^+  - \int_{\Omega} \log
\frac{\gamma_0(\sigma^+ | \sigma^+)}{\gamma_0( + |
\sigma^+)}\mu(d\sigma).
\end{equation}
\item If moreover $\mu \in \mathcal{G}_{\rm{inv}}(\gamma)$  exists, then \be\label{tension} h(\mu |\nu)=\lim_{\la
\uparrow \Z^d} \frac{1}{|\la|} \log
\frac{\mu(+_{\la})}{\nu(+_{\la})}. \ee
\end{enumerate}
\end{theorem}

To prove this theorem, we establish a more general lemma that will help  to restore
the thermodynamic properties of generalized Gibbs   in  Chapter 5. It topologically
captures what is  needed for a specification to get zero relative entropy of its DLR measures.
The full proof is given in \cite{KLNR}.
\begin{lemma}\label{chp:lem1}\cite{KLNR}
If $\mu(\Omega_{\gamma}^{<0})=1$, then \ben \item Uniformly in
$\om \in \Omega$,
\begin{displaymath}
\lim_{n \to \infty} \frac{1}{|\la_n|} \int_{\Omega} \log
\frac{\gamma_{\la_n}(\si|\omega)}{\gamma_{\la_n}(+|\omega)} \mu(d
\sigma)=\int_{\Omega}\log
\frac{\gamma_0(\si^+|\si^+)}{\gamma_0(+|\si^+)} \mu(d\sigma).
\end{displaymath}
\item For $\nu\in\mathcal{G}(\gamma)$,
\begin{displaymath}
\lim_{n \to \infty} \frac{1}{|\la_n|} \int_{\Omega} \log
\frac{\nu(\si_{\la_n})}{\nu(+_{\la_n})} \mu(d \sigma)=
\int_{\Omega} \log \frac{\gamma_0(\si^+|\si^+)}{\gamma_0(+|\si^+)}
\mu(d\sigma).
\end{displaymath}
In particular, the limit depends only on the pair $(\gamma,\mu)$.
\een
\end{lemma}
{\bf Proof :}
\begin{enumerate}
\item The proof relies on the uniform convergence of the
translation-invariant vacuum potential with vacuum state +
established  by Sullivan and described in Remark \ref{rksull}, using
a particular case of the telescoping procedure described in Remark
\ref{telescop}. We assumed there quasilocality of the specification
but it is not difficult to extend its validity when only almost-sure
continuity in the $+$ direction holds, which is in particular
implied by the condition $\mu(\Omega_\gamma^{<0})=1$. This latter
condition is not optimal\footnote{It is indeed extended in
\cite{KLNR}.}, and comes from the telescoping procedure used here,
but we do not know if it is technical or not, and in particular if
it can be relaxed to a general almost-sure quasilocality property, a
very important property in the context of generalized Gibbs
measures. Following Sullivan \cite{Su2}, we define, for all $\sigma
\in \Omega$
\begin{equation} \label{D}
\; \; D(\si)= E^+_{\{ 0 \}} (\si|\si)=\log
\frac{\gamma_0(\si|\si)}{\gamma_0(+|\si)}.
\end{equation}
and consider an approximation  of $\sigma^+$ at finite volume $\la \in \s$
with boundary condition $\om$ by defining the {\em telescoping
configuration} at $i \in S$ $T_{\la}^{\om}[i,\si,+]$, defined for all $j \in S$ by:
\[
\big(T_{\la}^{\om}[x,\si,+]\big)_j= \left\{
\begin{array}{lll}
\omega_j \;  &\rm{if}& \; j \in \la^c\\
\sigma_j \;&\rm{if}& \; j \; \leq \; i, \; j \in \la\\
+ 1 \; &\rm{if}& \;  j\; > \; i,\; j \in \la.\\
\end{array}
\right.
\]
To perform the telescoping,
denote $\la_{\leq i} = \{ j\in\Lambda: j\leq i\}$, $\la_{< i} =
\la_{\leq i}\setminus \{ i\}$, $\la_{> i} = \Lambda\setminus
\la_{\leq i}$ and  let $\la = \{ i_1, \ldots i_N\}$ denote an
enumeration of $\la$ in lexicographic order. By consistency and Lemma \ref{keylemma}: \be\frac{\gamma_\la (\si|\omega )}{\gamma_\la
(+|\omega )} = \prod_{k=1}^N \frac{\gamma_\la (\si_{\la_{\leq
i_k}}+_{ \la_{> i_k}} |\omega )} {\gamma_\la (\si_{\la_{\leq
i_{k-1}}}+_{ \la_{> i_{k-1}}} |\omega )}
= \prod_{k=1}^N \frac{\gamma_{i_k} (\si_{i_k}| \si_{\la_{< i_k}}
+_{\la_{> i_k}}\omega_{\la^c})} {\gamma_{i_k} (+_{i_k}|
\si_{\la_{< i_k}} +_{\la_{> i_k}}\omega_{\la^c})}. \ee Taking the
logarithm yields for $\Lambda=\Lambda_n$
\begin{eqnarray*}
\int_{\Omega} \log \frac{\gamma_{\Lambda_n}(\sigma|\omega)}{\gamma_{\Lambda_n}(+|\omega)} \;  \mu(d \sigma)&=&\sum_{i
\in \la_n} \int_{\Omega} D(\tau_{-i} T_{\la_n}^{\om}[i,\tau_i
\si,+]) \;  \mu(d \sigma).
\end{eqnarray*}
in such a way that proving Item 1. amounts to prove that,
uniformly in $\om$, the r.h.s divided by the volume converges to
$\int_\Omega D(\sigma^+) \; \mu(d\sigma)$. This is obtained in
\cite{KLNR} by carefully counting the points of the  set $A_n$
where this telescoping configuration and $\sigma^+$ differ, see
that $|A_n | = \circ(| \Lambda_n |)$, to get, due to the
continuity of $D$ at the
 configuration $\sigma^+$,  that for all $\epsilon >0$,
 $$
 \frac{1}{| \la_n  |} \Bigl|\; \sum_{i \in \la_n}
\big[ D(\tau_{-i} T_{\la_n}^{\om}[i,\tau_x \si,+])-D(\sigma^+
\big]\; \Bigr|
 \leq  \epsilon + 2 \sup_{\omega} | D(\omega)| \cdot \frac{|A_n|}{| \la_n  |}$$

which is less than $2
\epsilon$ for $n$ big enough. So we obtain that
\begin{displaymath}
\frac{1}{| \la_n  |} \Bigl|\; \sum_{i \in \la_n} \big[ D(\tau_{-i}
T_{\la_n}^{\om}[i,\tau_x \si,+])-D(\sigma^+) \big] \; \Bigr|
\end{displaymath}
converges to zero on the set of $\Omega^{< 0}_{\gamma}$ of full
$\mu$-measure, uniformly in $\om$,
which implies statement 1 of the lemma by dominated convergence. \item By consistency $\nu \in \mathcal{G}(\gamma)$, one rewrites for all $\sigma \in \Omega$
$$
\nu(\sigma_{\Lambda_n})=\int_{\Omega} \gamma_{\la_n}(\sigma |
\om) \nu(d \om)
$$
so that
\[
F_{\la_n}(\mu,\nu):=\frac{1}{ | \la_n |} \int_{\Omega} \log
\frac{\nu(\sigma_{\la_n})}{\nu(+_{\la_n})} \mu(d \sigma)=\frac{1}{| \la_n  |}
\int_{\Omega} \log \frac{\int_{\Omega} \gamma_{\la_n}(\sigma |
\om) \nu(d \om)}{\int_{\Omega} \gamma_{\la_n}( + | \om) \nu(d
\om)} \mu(d \si).
\]
We use now the obvious bound
\begin{displaymath}
\inf_{\om \in \Omega}\frac{\gamma_{\la_n}(\sigma | \om)}{
\gamma_{\la_n}( + | \om) } \; \leq \; \frac{\int_{\Omega}
\gamma_{\la_n}(\sigma | \om) \nu(d \om)}{\int_{\Omega}
\gamma_{\la_n}( + | \om) \nu(d \om)} \; \leq \; \sup_{\om \in
\Omega} \frac{\gamma_{\la_n}(\sigma | \om)}{ \gamma_{\la_n}( + |
\om)}.
\end{displaymath}

to get for $\epsilon >0$  given that there exists
$\om=\om(n,\sigma, \epsilon)$, $\om'=\om'(n,\sigma, \epsilon)$
such that
\[
\int_{\Omega} \inf_{\om \in \Omega}\log
\frac{\gamma_{\la_n}(\sigma | \om)}{ \gamma_{\la_n}( + | \om) }
\mu(d \sigma) \geq \int_{\Omega} \log
\frac{\gamma_{\la_n}(\sigma|\om(n,\sigma,
\epsilon))}{\gamma_{\la_n}(+|\om(n,\sigma, \epsilon))} - \epsilon
\]
and
\[
\int_{\Omega} \sup_{\om \in \Omega}\log
\frac{\gamma_{\la_n}(\sigma | \om)}{ \gamma_{\la_n}( + | \om) }
\mu(d \sigma) \leq \int_{\Omega} \log
\frac{\gamma_{\la_n}(\sigma|\om'(n,\sigma,
\epsilon))}{\gamma_{\la_n}(+|\om'(n,\sigma, \epsilon))} +
\epsilon.
\]
Now use the first item of the lemma and choose $N$ such that for
all $n \geq N$,

\[
\sup_{\om} \Big| \frac{1}{|\la_n |}\int_{\Omega} \log
\frac{\gamma_{\la_n}(\sigma | \om) }{\gamma_{\la_n}( + | \om) }
\mu(d \sigma) - \int_{\Omega} D(\sigma^+) \mu(d \si) \Big| \leq
\epsilon
\]
to get, for $n \geq N$,
\begin{displaymath}
\int_{\Omega} D(\sigma^+) \mu(d \sigma) -2 \epsilon \; \leq \;
F_{\la_n} (\mu | \nu) \; \leq \; \int_{\Omega} D(\sigma^+)
\mu(d\si) + 2 \epsilon
\end{displaymath}
and eventually prove the lemma.

\end{enumerate}

{\bf Proof of Theorem \ref{chp:thm2}:}

\begin{enumerate}
\item
By a short computation at finite volume, rewrite \be
\label{dec1} h_n(\mu | \nu)=\, -h_n(\mu) -\,
\frac{1}{|\la_n|}\sum_{\sigma_{\la_n} \in \Omega_{\Lambda_n}} \mu(\sigma_{\la_n})\log
\frac{\nu(\sigma_{\la_n})}{\nu(+_{\la_n})}\, -
\,\frac{1}{|\la_n|}\log \nu(+_{\la_n}). \ee When $\mu(\Omega_\gamma^{<0})=1$
holds, the asymptotic behavior of the second term of the r.h.s. is
given by  Lemma \ref{chp:lem1} and under  the existence of $e_\nu^+$ one gets (\ref{chp:exist}).
\item For  $\mu \in \mathcal{G}_{\rm{inv}}(\gamma)$ such
that $\mu(\Omega_{\gamma}^{< 0})=1$,  rewrite now :
$$ h_n(\mu|\nu)=
\frac{1}{|\la_n|}\Big(\sum_{\sigma_{\la_n} \in \Omega_{\Lambda_n}}
\mu(\sigma_{\la_n})\log \frac{\mu(\sigma_{\la_n})}{\mu(+_{\la_n})}
- \sum_{\sigma_{\la_n} \in \Omega_{\Lambda_n}}
\mu(\sigma_{\la_n})\log \frac{\nu(\sigma_{\la_n})}{\nu(+_{\la_n})}
+ \log \frac{\mu(+_{\la_n})}{\nu(+_{\la_n})}\Big). $$ \normalsize
By Lemma \ref{chp:lem1}, in the limit $n\to\infty$, the first two
terms of the r.h.s. are functions of $(\gamma,\mu)$ rather than
functions of $\nu, \mu \in \mathcal{G}_{\rm{inv}}(\gamma)$ and
cancel out. Hence, the relative entropy exists if and only if the
fourth term converges. Using Item 1 (existence of relative
entropy), we obtain the existence of the limit in (\ref{tension})
and the equality
\[
h(\mu | \nu)=\lim_{n \to \infty} \frac{1}{| \, \la_n \, |} \log
\frac{\mu(+_{\la_n})}{\nu(+_{\la_n})}
\]
and in particular the r.h.s is a well-defined limit.
\end{enumerate}

\subsection{Application: VP for translation-invariant quasilocal measures}
\begin{theorem}\label{VPQl}\cite{KLNR}
Any $\mu \in \mathcal{M}^+_{1,{\rm
inv}}(\Omega)$ quasilocal satisfies the variational principle (\ref{eq:r4}).
\end{theorem}

{\bf Proof:} When $\gamma$ is quasilocal, one has $\Omega_\gamma=\Omega$ and the second part
is a direct consequence of Theorem \ref{r.teo1}, because the convergence required is true by
any convergence theorem (uniform or dominated). We recover thus the usual standard proof of
\begin{theorem}[2nd part of the VP for quasilocal measures]\label{2dpartql} Let $\gamma$ be a
 translation-invariant quasilocal specification and $\nu \in \mathcal{G}_{\rm
inv}(\gamma)$. Then for any  $\mu\in\mathcal{M}_{1,{\rm inv}}^+(\Omega)$,
$$
h(\mu|\nu)=0 \; \; \Longrightarrow \; \; \mu\in \mathcal{G}_{\rm inv}(\gamma).
$$
\end{theorem}

Thus, the  result of \cite{KLNR}, that extends the variational
principle from {\em translation-invariant Gibbs measures with a
translation-invariant UAC potential} to {\em translation-invariant
quasilocal measures}, for which the UAC potential derived from
Kozlov (and described in the previous chapter) is not necessarily
UAC, relies on the following
\begin{theorem}[1st part of the VP for quasilocal measures]\label{chp:thm7}
Let $\gamma$ be a translation-invariant quasilocal specification,
$\nu \in \mathcal{G}_{\rm{inv}}(\gamma)$ and $\mu \in
\mathcal{M}_{1,\rm{inv}}^+(\Omega)$. Then $h(\mu | \nu)$ exists
for all $\mu \in \mathcal{M}_{1,\rm{inv}}^+(\Omega)$ and
\[
\mu \in \mathcal{G}_{\rm{inv}}(\gamma) \, \Longrightarrow \,
h(\mu|\nu)=0.
\]
\end{theorem}

{\bf Proof:} We need the following lemma to use Theorem \ref{chp:thm2}:

\begin{lemma}\label{chp:lem2}
For  $\mu,\nu \in \mathcal{G}_{\rm{inv}}(\gamma)$ with $\gamma$
t.i. and quasilocal, $e^+_{\nu},e^+_{\mu}$ exist
and
\begin{displaymath}
\lim_{n \to \infty} \frac{1}{|\la_n|} \log
\frac{\mu(+_{\la_n})}{\nu(+_{\la_n})}=0.
\end{displaymath}
\end{lemma}

{\bf Proof:} It is a direct consequence of the  uniform convergence of the vacuum potential $\Phi^+$ that we now
 associate with $\gamma$ for the reference  configuration $+$, keeping the notation of Chapter 3.
 It is a consequence of another lemma from \cite{KLNR},  which is an adaptation of the
argument used by Israel \cite{Is} to prove existence and boundary
condition independence of the pressure for a UAC potential.

\begin{lemma}\label{chp:lem3}
The vacuum potential with vacuum state $+$ associated with the
quasilocal specification $\gamma$ is such that:
\begin{enumerate}
\item
\begin{equation}\label{chp:e4.3.1}
 \lim_{n \to \infty} \sup_{\omega,\eta,\sigma}
\frac{1}{|\la_n|} \Bigl|  H_{\la_n}^{\Phi^+}(\sigma|\eta) -
H_{\la_n}^{\Phi^+}(\sigma|\omega) \Big| =0.
\end{equation}
\item
\begin{equation}\label{chp:e4.3.2}
\lim_{n \to \infty} \sup_{\omega,\eta}\frac{1}{|\la_n|} \log
\frac{Z_{\la_n}^{\Phi^+}(\omega)}{Z_{\la_n}^{\Phi^+}(\eta)} =0.
\end{equation}
\end{enumerate}
\end{lemma}
{\bf Proof:} Clearly, (\ref{chp:e4.3.1}) implies
(\ref{chp:e4.3.2}): For all $n \in \mathbb{N}$,
\begin{displaymath}
\exp\Big\{-\sup_{\omega,\eta,\sigma} \Bigl|
H_{\la_n}^{\Phi^+}(\sigma|\eta) - H_{\la_n}^{\Phi^+}(\sigma|\omega)\Bigr|
\Big\} \leq
\sup_{\omega,\eta}\frac{Z_{\la_n}^{\Phi^+}(\omega)}{Z_{\la_n}^{\Phi^+}(\eta)} \leq
\exp\Big\{\sup_{\omega,\eta,\sigma} \Bigl|
H_{\la_n}^{\Phi^+}(\sigma|\eta) - H_{\la_n}^{\Phi^+}(\sigma|\omega)\Bigr|
\Big\}.
\end{displaymath}
Now one  proves (\ref{chp:e4.3.1}) by rewriting
\begin{displaymath}
H_{\la_n}^{\Phi^+}(\sigma|\eta) - H_{\la_n}^{\Phi^+}(\sigma|\omega)=\sum_{A \cap
\la_n \neq \emptyset, A \cap \la_n^c \neq \emptyset} \Bigl[
\Phi^+_A(\sigma_{\la_n}\eta_{\la_n^c})-\Phi^+_A(\sigma_{\la_n}\omega_{\la_n^c})
\Bigr].
\end{displaymath}
which goes uniformly to zero by  relative uniform convergence, see \cite{KLNR} for details.\\

 To derive Lemma \ref{chp:lem2} from Lemma \ref{chp:lem3}, we only
have to prove that for all $\nu \in
\mathcal{G}_{\rm{inv}}(\gamma)$, $e_{\nu}^+$ exists and is
independent of $\gamma$. For such a measure $\nu$, write, using  Lemma
\ref{chp:lem3}
\begin{displaymath}
\nu(+_{\la})=\int_{\Omega}
\frac{e^{-H^{\eta}_{\la_n}(+)}}{Z_{\la_n}(\eta)}\nu(d\eta)\cong \int_{\Omega}
\frac{e^{-H_{\la}^+(+)}}{Z_{\la}^+} \nu(d\eta)
\end{displaymath}
where $a_{\la} \cong  b_{\la}$ means $\lim_{\la} \frac{1}{|\la|} |
\log \frac{a_{\la}}{b_{\la}} |=0$. Since $\Phi^+$ is the vacuum
potential with vacuum state $+$, $H_{\la}^{+}(+_{\la})=0$ and
\begin{displaymath}
\nu(+_{\la})=\frac{1}{Z_{\la}^{\Phi^+}}(+)=\frac{1}{Z_{\la}^{\Phi^+,\rm{f}}}=\Big[\sum_{\sigma \in \Omega_{\la}} \exp (-\sum_{A \subset \la}
\Phi_A^+(\sigma)) \Big]^{-1}
\end{displaymath}
 Fix $R>0$ and define a finite-range potential $\Phi^R$ by putting  $\Phi^{(R)}_A(\sigma) :=  \Phi^+_A(\sigma)$ if
$|A| \leq R$  and $\Phi^{(R)}_A(\sigma):=0$ otherwise. To use the
existence of pressure $P(\Phi^{(R)})$ from  \cite{Is} for
finite-range (translation-invariant) potentials sketched in the
beginning of this chapter, we write
\begin{eqnarray*}
\log \frac{\sum_{\sigma} \exp{(-\sum_{A \subset \la}
\Phi_A^+(\sigma))}}{\sum_{\sigma} \exp{(-\sum_{A \subset \la}
\Phi_A^{(R)}(\sigma))}} \; & \leq &\; \sup_{\sigma} \Big|\sum_{A
\subset \la,|A| > R} \Phi_A^+(\sigma) \Big|
 \leq  \sup_{\sigma}  \sum_{i \in \la} \Big| \sum_{A \ni i,
|A|>R} \Phi_A^+(\sigma) \Big|\\
&\leq & \sum_{i \in \la} \sup_{\sigma} \Big|\sum_{A \ni i,
|A|>R} \Phi_A^+(\sigma) \Big|
= |\la| \cdot \sup_{\sigma} \Big|\sum_{A \ni 0, |A|>R}
\Phi_A^+(\sigma) \Big|.
\end{eqnarray*}
This is the tail of a uniformly convergent series (see Remark
\ref{rksull}) from which we conclude by that $\{P(\Phi^{(R)}), R>0
\}$ is a Cauchy net with limit
\begin{displaymath}
\lim_{R \to \infty} P\big(\Phi^{(R)}\big)=\lim_{\la \uparrow
\mathbb{Z}^d} \frac{1}{|\la|} \log
Z_{\la}^{\Phi^+,\rm{f}}=\lim_{\la \uparrow \mathbb{Z}^d}
\frac{1}{|\la|} \log Z_{\la}^{\Phi^+}(+)=-e_{\nu}^+
\end{displaymath}
which depends only on the vacuum potential (hence on the
specification $\gamma$). This proves that $e_{\nu}^+$ and
$e_{\mu}^+$ exist for all $\mu, \nu \in
\mathcal{G}_{\rm{inv}}(\gamma)$, and depends of $\gamma$ only.
Therefore,
\begin{displaymath}
\lim_{\la \uparrow \Z^d} \frac{1}{|\la|} \log
\frac{\mu(+_{\la})}{\nu(+_{\la})}=e_{\nu}^+-e_{\mu}^+=0.
\end{displaymath}

\br \label{rkpressure} \er
 In the standard theory of Gibbs measures, the existence of
$h(\mu | \nu)$ and the identity \ref{chp:exist} are obtained by
proving existence and boundary condition independence of the
pressure, see e.g. \cite{Ge} p 322  for a similar expression. This requires the existence of a UAC potential, which in
our case is replaced by regularity properties of the specification
and existence of the limit  defining $e_{\nu}^+$. The existence is
guaranteed e.g. for renormalization group transformations of Gibbs
measures, and for $\nu$ with positive correlations (by
subadditivity). Moreover, in the case of transformations of Gibbs
measures, convergence to zero of (\ref{tension}) is also easy to verify \cite{KLNR}.

\section{More on equilibrium: LDP and Stochastic Ising models}

\subsection{Large deviation properties}

The thermodynamical variational principle (\ref{chp:def3}), which
holds for translation-invariant Gibbs measures consistent with a
translation-invariant potential, as proved e.g. in \cite{Ge}, or in
\cite{P} for the more general class of asymptotic decoupled
measures, and in particular the convex conjugation of the relative
entropy and the pressure, is a first step for the statement of a
large deviation principle for such measures. Indeed, when relative
entropy is defined, it is often the so-called rate function for such
an LDP at the level of measures. We do not enter into the general
details of this important subject of probability theory in these
lectures,  see e.g. \cite{Com,LP,LPS,VEFS} for a precise formulation
of  large deviation principles and their relationships with the
notion of entropy in thermodynamics.

Roughly speaking, large deviations consist in an estimation of the
probability of rare events by estimating the usually very small
probabilities of large simultaneous fluctuations in a system
consisting of a large number of random variables, and claiming
that such a principle holds at the level of measures is a way of
expressing the fact that the probability that a typical
configuration for a measure $\mu \in \mathcal{M}_{1,\rm{inv}}^+$
looks  typical in $\Lambda$ for the measure $\nu \in
\mathcal{M}_{1,\rm{inv}}^+$  decays exponentially fast with the
volume with a rate equals to the relative entropy $h(\mu |\nu)$:
$$
{\rm Prob}_\nu \big[\omega_\Lambda \; {\rm typical \; for} \;
\mu_\Lambda \big] \; \approx \; e^{-|\Lambda| \; h(\mu |\nu)}.
$$

Thus, when such a principle holds for our Gibbs measures, the rate
function gets its minimum at zero when $\mu$ and $\nu$ are Gibbs
for the same specification, and thus the probability of getting
from $\mu$ a typical configuration for $\nu$ decays exponentially
with the {\em order of the surface only}, or at least at a
sub-volumic rate. This fact can be used to prove that a measures
is not Gibbs for the same potential as a reference measure and has
indeed been used to detect non-Gibbsianness in the projection of
the Ising model \cite{sch}.
\subsection{Stochastic Ising models}

Stochastic Ising models are particular types of Markov Processes on
the configuration space $(\Omega,\mathcal{F})$ when the single-site
state space is  $E=\{-1,+1\}$. They are widely described in
\cite{ligg} and allow, under mild conditions, to get Gibbs measures
as invariant reversible measures for these stochastic processes. Let
us focus on the most standard local stochastic dynamics, the
so-called {\em Glauber or spin-flip dynamics}, which corresponds to
usual birth and death processes in standard probability theory.
Starting from an a priori configuration, one would like to  change,
or flip, the configuration at a given site randomly depending on its
neighbors in order to get a suitable convergence to a typical
configuration of a given Gibbs measure.

To formalize this a bit, one denotes for any $i \in S$ and $\sigma
\in \Omega$ the flipped configuration $\sigma^i$ to be defined by
$\sigma^i_i=-\sigma_i$ and $\sigma^i_j=\sigma_j$ for all $i\neq
j$, and consider a collection of {\em spin-flip rates} $\big \{
c_i(\sigma), i \in S, \sigma \in \Omega \}$, assumed to be of
finite-range, strictly positive and translation-invariant, in
order to uniquely define a Feller process $(\eta_t)_{t \geq 0}$ on
$\Omega,\mathcal{F}$ with generator\footnote{For a  rigorous
description, consult \cite{ligg,EFHR}, one can consider the
closure of the generator.} $L$ defined on local functions $f \in
\mathcal{F}_{\rm{loc}}$ by
$$
Lf (\sigma)=\sum_{i \in S} c_i(\sigma) \big[f(\sigma^i) - f(\sigma) \big].
$$
Denoting by $S(t)$ the corresponding semi-group (see \cite{ligg})
and by $\mathbb{E}_\sigma$
the expectation under the corresponding path-space measure $\mathbb{P}_{\sigma}$ given the
initial configuration $\eta_0=\sigma$, one gets an action on functions $f \in \mathcal{F}_{\rm{loc}}$ with for all $t >0$,
$$
S(t) f (\sigma)= \mathbb{E}_\sigma [f(\eta_t)]
$$
and on measures $\nu$ to get a measure $\nu S(t)$ defined by its
expectations on local functions
$$
\int f d(\nu S(t)) =  \int S(t) f d \nu.
$$
The corresponding measure $\nu S(t)$ being thus the distribution
of the configuration at time $t$ if the initial distribution at
time zero is $\nu$. A probability measure $\mu\in
\mathcal{M}_{1,\rm{inv}}(\Omega)$ is then called {\em invariant}
for the process (or {\em for the dynamics}) with generator $L$ iff
$$
\int Lf d\mu =0, \; \forall f \in \mathcal{F}_{\rm{loc}}
$$
or equivalently iff $\mu S(t)=\mu$ for all time $t$; an invariant
measure $\mu$ is {\em reversible} when
$$
\int (Lf) g d\mu =0, \; \forall f,g \in \mathcal{F}_{\rm{loc}}.
$$
In words, a probability measure is invariant when the process
$(\eta_t)_t$ obtained by using $\mu$ as initial distribution is
stationary in time, so that its definition can be extended to
negative times, and is thereafter reversible  when the process
$(\eta_t)_t$ and $(\eta_{-t})_t$ have the same distribution. In
our case, reversibility is equivalent to a standard detailed
balance condition on the rates \be \label{rev} \frac{d \mu^i}{d
\mu}= \frac{c_i(\sigma^i)}{c_i(\sigma)} \ee where $\mu^i$ is the
image law of $\mu$ by the spin-flip at site $i$, $\sigma
\longmapsto \sigma^i$. Thus,  to get a dynamics evolving towards a
given measure $\mu$, it is enough to choose the rates according to
(\ref{rev}). This is the way Gibbs measures are obtained as
reversible invariant measures in the so-called {\em Glauber
dynamics at inverse temperature} $\beta$: Given a UAC potential
$\Phi$, one introduces the rates
$$
c_i(\sigma)=\exp{\Big\{\frac{\beta}{2} \sum_{A \ni i}
\big[\Phi_A(\sigma)-\Phi_A(\sigma^i)\big] \Big\}}
$$
in order that (\ref{rev}) holds for the Gibbs measures corresponding to the potential $\Phi$.
In such a case, the rates generates jump-processes on the configuration space with independent
Poisson clocks attached at each site that randomly produce spin-flips  according to
 the considered Gibbs measures, which are eventually the invariant reversible measures reached at equilibrium.\\

Similar stochastic Ising models can be introduced by changing the
spin-flip rules (in the so-called Metropolis-Hastings dynamics) or
by exchanging the spins between two sites (Kawasaki dynamics,
equivalent to an exclusion-process in the lattice gas settings with
single site state-space $E=\{0,1\}$), leading at equilibrium to the
same reversible Gibbs measures.

\chapter{Generalized Gibbs measures}
\section{Heuristics}

Let us consider physical systems with a large number of particles in
thermal equilibrium modelled by the Gibbsian formalism described in
Chapter 3 and consider more precisely the example of particles of
water. Although water is too complicated a system to be described
precisely by the Gibbsian formalism\footnote{In particular, the
solid phase has to be more carefully described as we do now, due to
the particular nature of crystals. One expects that a crystal breaks
the translation symmetry, so that translation-invariant Gibbs
measures can be ergodic but not extremal in the set of Gibbs
measures. Moreover, although there is a gas-liquid point, there is
not a second-order transition from solid to liquid or gas in great
generality. Our present heuristics have thus to be taken very
carefully while dealing with this solid phase.}, it allows to give a
qualitative picture of the phase transition phenomenon, in
accordance with the precise description of the phase diagram that
can be achieved within the Gibbsian formalism by the Pirogov-Sinai
theory \cite{PS}. The system could be in different states depending
on the temperature, and assume that these states are described by
extremal Gibbs measures $\mu_S$, $\mu_L$, $\mu_V$ at low temperature
or $\mu$ at high temperature. One observes the existence of a  {\em
critical temperature} $T_c$ that distinguishes a region of
temperature where the physical system can only be in a unique phase
and a lower dimensional manifold\footnote{The existence of these
manifolds comes from the {\em Gibbs phase rule}, see the
introduction of Wightman in \cite{Is}.} corresponds in the
$(P,T)$-plane where the system can coexist in two or three different
states, depending on the pressure for a given temperature, yielding
the following (qualitative) phase diagram.
\begin{description}
\item[$T \; > \; T_c$:] Uniqueness regime: $\mathcal{G}(\gamma)=\{\mu\}$.
\item[$T \; < \; T_c$:] There exist $2$d-manifolds where the system is in
  a unique phase, solid, liquid or gaseous, depending on the pressure $P$. These unicity manifolds
  have as boundaries $1$d-manifolds where two different phases coexists. These coexistence lines have as
  boundary a $0$d-manifold where all three phases coexist.
\end{description}
>From a physical point of view, a phase transition is the
transformation of $(P,T)$-variables that allows the passage from
one of the regions of uniqueness to another one, through a region
of non-uniqueness. Two kind of phase transitions can be
distinguished  here: When this $(P,T)$-transformation crosses a
coexistence line, one says that  a first order phase transition
occurs, whereas when one crosses the  {\em
  critical point} C from $T>T_c$, one says that it is  second order\footnote{This distinction between first or second order phase
  transitions is also mathematically characterized in terms of differentiable properties of the pressure introduced in Chapter 4,
  see \cite{Is,Simon,VEFS}.}. Let us focus on the latter, for which a quantity called {\em  correlation length} is introduced.
  Let $\mu$ be a Gibbs
measure, $\sigma\in\Omega$. For $T \neq T_c$, and $i \neq j \in S$, the covariance between
the random variables $\sigma_i$ and $\sigma_j$ is expected\footnote{and sometimes proved \cite{LebPe} using
among others {\em correlation inequalities}.}  to decay exponentially in such a way that one could define and
define a quantity $\xi(T)$
such that:
\begin{displaymath}
\mu(\sigma_i\sigma_j)-\mu(\sigma_i)\mu(\sigma_j) \; \cong \;
e^{-\frac{\mid i-j\mid}{\xi(T)}}
\end{displaymath}
where $\cong$ means logarithmic equivalence for large $\mid i - j \mid$. The quantity
$\xi(T)$ has the dimension of a distance and is interpreted as the {\em correlation length} of the system, beyond
which two spins are physically considered to be independent. It is
then considered as a natural scale of the system which enables us to measure the length with the unit $ 1 \; \xi$
instead of $1$ meter. We cannot do it at the
critical temperature because the decay of  correlation is not expected to be exponential and in some sense one has presumably
\begin{displaymath}
\lim_{T \to T_c} \xi(T)=+ \infty.
\end{displaymath}
This is interpreted as the absence of proper scale for the system
at the critical temperature. Physically, for the system of water,
we observe a ''milky'' water or critical opalescence, showing a
strong interaction between all the particles of the system that
creates a highly chaotic behavior. As a consequence of this
absence of proper scale, the behavior of the system at the
critical point should be the same at any scale, providing a tool
to  study these critical behaviors, which are ill-known and
difficult to observe: Natural transformations of the observation
scale seem to be an appropriate tool to understand it better,  the
critical point being considered in some sense  as a fixed  point
of these transformations. This has motivated the introduction of
the {\em renormalization group}, a semi-group of transformations
directly related to a change of scale of the system,  as a tool in
theoretical physics to study critical phenomena using change of
scales in particles systems, which  appeared to be rather powerful
 in these fields, see e.g. \cite{RGT2,DG,fis,RGT,GP,wil}.\\

  We shall describe it more precisely next section, but let us first
consider a scaling transformation $T:\Omega \longrightarrow
\Omega'$, where $\Omega$ and $\Omega'$ are the configuration spaces
at two different scales. Let $\mu$ be a measure describing an
equilibrium state of the system at the first scale, i.e. a Gibbs
measure on $\Omega$, and denote formally $H$ its Hamiltonian. The
transformation $T$ acts naturally on measures and we denote
$\mu'=T\mu \in \mathcal{M}^{+}_{1}(\Omega',\mathcal{F}')$. The
natural aim in our theory would be to obtain $\mu'$ as a Gibbs
measure and to define an Hamiltonian $H'$, image of the Hamiltonian
$H$ by a renormalization transformation on spaces of Hamiltonians
 for it,  in order to get the following
diagram  defined and commutative.
\begin{displaymath}
\epsfbox{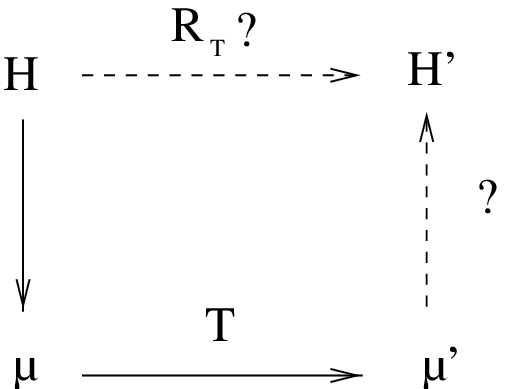}
\end{displaymath}
\begin{center}
Figure 1: A renormalization group transformation acting on measures.
\end{center}
In the late seventies and early eighties, Griffiths and
Pearce/Israel \cite{GP,Is2} discovered some pathologies of the
behavior of these image measures: It turned out that they often
are not Gibbsian and it came out as a surprise that one could
break the equilibrium properties of a state by only looking at it
at a different scale, and according to the ideology of the
renormalization-group theory, this should not be so. These
pathologies have been rigorously proven to exist and mostly
identified as the manifestation of non-Gibbsianness due to a
failure of the quasilocality property, mostly at low temperature,
in 1993 by van Enter, Fern\'andez and Sokal \cite{VEFS} in a
rather general and rigorous mathematical framework. Let us recall
that the Hamiltonian can be recovered from the Gibbs measure $\mu$
with the help of the Moebius inversion formula. When the
correlation do not decrease sufficiently fast, and in particular
at the critical point or in some phase transitions region, a
divergence might appear in this inversion, preventing a definition
of $H'$ from $\mu'$,  and leading to non-Gibbsianness of the
image-measure. Nevertheless, we shall see that this phenomenon
also hold in other parts of the phase diagram, sometimes far from
the critical point, and that it is already present after one
single change of scale for very simple scaling transformations. We
formalize all these heuristics more precisely now.

\section{RG pathologies and non-Gibbsianness}

In this section, we formalize mathematically the scaling transformations,
introduce the renormalization group transformations and describe how these
transformations could lead to non-Gibbsianness. As we shall see, this phenomenon
is often related to the occurrence of phase transitions in some {\em hidden or constrained system}, and to
get the main features of the phenomenon we describe precisely how this happens  for the simplest
 RG transformation, the so-called {\em decimation of the 2d Ising model at low temperature},
 the latter being low enough to get a phase transition that creates long-range dependencies
 leading to non-quasilocality in the mentioned hidden system. Thereafter, we shall give a
 (non-exhaustive) catalogue of other RG transformations of Gibbs measures that also lead to non-Gibsianness
 for similar reasons, but for which the proof is much more complicated, when it exists !

\subsection{Decimation of the $2d$ ferromagnetic Ising model \cite{VEFS}}

The basic example we describe here, which already captures the
main non-trivial features of the pathologies, concerns the
decimation with spacing $2$, which corresponds to the projection
of the 2d-Ising model on the sublattice of even sites, while we
shall give similar results for decimations at other dimensions and
for other spacings later on. More generally, the decimation
transformation  on $\mathbb Z^{2}$ with spacing $b$ is defined to
be the transformation
\begin{eqnarray*}
T_b \colon (\Omega,\mathcal{F}) & & \longrightarrow (\Omega',\mathcal{F}')=(\Omega,\mathcal{F})\\
\omega \; \; & & \longmapsto \omega'=(\omega'_i)_{i \in \mathbb{Z}^2}
\end{eqnarray*}
defined by $\forall i \in \mathbb Z^{2}$ by $\omega'_{i}=\omega_{bi}$, and denote simply $T=T_2$ the case
of a spacing $b=2$ described here. The following result is crucial to understand renormalization group
pathologies and the arising of non-Gibbsianness in equilibrium mathematical statistical mechanics.
This example, already described by Israel in \cite{Is2} to detect renormalization group pathologies,
 has been fully analyzed in \cite{VEFS}, in the seminal  description of RG pathologies in
 the realm
 of Gibbsianness vs. non-Gibbsianness framework.

\begin{theorem}\cite{VEFS}\label{Dec} Let $\beta  > \tilde{\beta}_c=\frac{1}{2} \cosh^{-1}\big(e^{2 \beta_c}\big)$
and denote by $\nu_\beta = T \mu_\beta$ the decimation of any Gibbs measure $\mu$ for the homogeneous ferromagnetic
n.n Ising model on $\Z^2$ with zero magnetic field. Then $\nu_\beta$ is not quasilocal, hence non-Gibbs.
\end{theorem}

To prove non-quasilocality of the renormalized measure, one exhibits
a so-called {\em bad configuration} where the conditional
expectation,  w.r.t. the outside of a finite set,  of a local
function, is essentially discontinuous, or equivalently is
discontinuous on a (non-negligible) neighborhood, as described in
Chapter 3. The role of a bad configuration is played  by a so-called
{\em alternating  configuration} $\omega'^{\rm{alt}}$ defined for
all $i=(i_1,i_2) \in \mathbb{Z}^2$ by
$\omega'^{\rm{alt}}_i=(-1)^{i_1+i_2}$. Computing the magnetization
under the image measure $\nu_\beta$, conditioned on the boundary
condition $\omega'^{\rm{alt}}$ outside the origin, will give
different limits when one approaches
 this configuration with all $+$ (resp. all $-$) arbitrarily far away, as soon as phase transition is possible
 for the Ising model on the so-called {\em decorated lattice}, a version of $\mathbb{Z}^2$ where even sites have
  been removed. The latter phase transition is shown to be possible as soon as the inverse temperature is larger
   than the above value $\tilde{\beta}_c$. The global neutrality of this bad configuration leaves the door open
   to such a phase transition to occur at low enough temperature, and is crucial in its badness. Let us formalize
    this a bit more, following the full proof given in \cite{VEFS}.\\

{\bf Proof of theorem \ref{Dec}:}

We denote by $\nu= T \mu$ this {\em decimated measure}:
$$\forall A' \in \mathcal{F}', \; \nu(A')=\mu(T^{-1}(A'))=\mu(A)$$ with the notation $A=T^{-1}(A')\; \in \mathcal{F}$.
In order to describe how a phase transition in some hidden system  gives rise to non-Gibbsianness, we also extend this
 decimation  on the ''even'' sites of $\mathbb Z^{2}$, i.e. $2 \Z^{2}$,  by:
\begin{displaymath}
T: i=2i' \longmapsto i'
\end{displaymath}
and on subsets: $\forall \Lambda \subset 2 \Z^{2}, T(\Lambda):= \Lambda'=\{i \in \mathbb Z^{2},\;2i \in \Lambda \} \subset \Z^{2}.$
We underline that $T$ maps the finite subsets of $2 \Z^{2}$ on the finite subsets of $\Z^{2}$, but the converse is not true:
 If one defines a {\em cofinite} subset of $S$ to be
the complement of any finite set, then the inverse
transformation $T^{-1}$ \emph{does not} map the \emph{co}finite subsets of $\Z^{2}$ on the \emph{co}finite subsets of $2 \Z^{2}$. For example, when $\Lambda'=\{{\mathbf 0}\}$ consists of the origin ${\mathbf 0}$ of $\Z^2$. Then $\Lambda'^{c}=\mathbb Z^2 \backslash \{{\mathbf 0}\}$
\begin{eqnarray*}
{\rm and} \; \; \; \; T^{-1}(\Lambda'^{c})&=&\{i=2i' \; {\rm s.t.}\; i' \in \Lambda'^{c}\}=\{x=2x',\;  x' \in \mathbb Z^{2},\; x' \neq 0 \}\\
&=&2\mathbb Z^2 \backslash \{{\mathbf 0}\}=\big[(\mathbb Z^2
\backslash 2\mathbb Z^2)\cup \{{\mathbf 0}\}\big]^{c}=:\lambda^c
\end{eqnarray*}
where $\lambda=(\mathbb Z^2 \backslash 2\mathbb Z^2)\cup
\{{\mathbf 0}\}$ is {\em not} a finite subset of $\mathbb Z^2$,
and this is actually the reason why non-Gibbsianness could occur:
Getting still an infinite-volume framework after this
conditioning,  this leaves the possibility of  phase transitions
for the measure conditioned on $\lambda$, giving rise to the
essential discontinuity: In order to prove that $\nu$ is not
Gibbsian, we prove that there exists $\Lambda' \in \mathcal{S}'$
and a function $f$ local on $\Omega'$ such that no version of
$\nu[f|\mathcal{F}_{\Lambda'^{c}}](\cdot)$ is quasilocal, that is
we want to find $\omega'$ in $\Omega'$ for which there exists  $f$
local on $\Omega'$ with
$\nu[f|\mathcal{F}_{\Lambda'^{c}}](\omega')$ \emph{essentially
  discontinuous}. Now, using the action of the measurable map $T$ on subsets and configurations,
  one easily gets, for any  $\Lambda' \in \mathcal{S}'$
\begin{displaymath}
\nu[A'|\mathcal{F}'_{\Lambda'c}](\omega')=\mu[A|\mathcal{F}_{\lambda^c}](T^{-1}\omega')
, \; \nu \textrm{-a.e.}(\omega'), \; \forall A' \in \mathcal{F}'.
\end{displaymath}
So we have to compute the conditional probabilities
$\mu[A|\mathcal{F}_{\lambda^c}]$ for $\lambda$ {\em non-finite},
and this is not given by the specification $\gamma$, which only
provides versions of conditional probabilities for the outside of
finite sets only\footnote{To describe such conditional
probabilities, on should use {\em global specifications} developed
in \cite{FP}.}. Thus  $\lambda=T^{-1}(\Lambda'^{c})$ \emph{is not}
a cofinite set,  as we show for $\Lambda'=\{\mathbf{0} \}$ and
illustrated on the figures 2 below. A short computation leads
indeed to
\begin{displaymath}
\lambda^{c} \;= \mathbb Z^{2} \backslash (2\Lambda') \;= \;T^{-1}(\Lambda') \cup (\mathbb Z^{2} \backslash 2\mathbb Z^{2})
\end{displaymath}


and thus  $\lambda^{c}$ consists of all the spins of $2\mathbb
Z^{2}$ except the origin: If we ''knew'' everything except the
origin on the decimated system $\Omega$, we should ''know'' the
spins on $2\mathbb Z^{2}$ except at the origin. Then $\lambda$ is
the origin \emph{plus} the sites which are not in $2\mathbb
Z^{2}$, as shown in figures 2 below\footnote{The letters denote
the value of spins on the underlying sites, when fixed,  and the ?
indicate that the spin over the underlying site is unknown.}:

\begin{center}
\begin{picture}(80,80)
\graphpaper[40](-1,-1)(80,80)
\put(35,35){\makebox(0,0){O}}
\put(44,44){\makebox(0,0){?}}
\put(-6,77){\makebox(0,0){a}}
\put(34,74){\makebox(0,0){b}}
\put(82,82){\makebox(0,0){c}}
\put(-6,37){\makebox(0,0){d}}
\put(82,42){\makebox(0,0){e}}
\put(-3,-6){\makebox(0,0){f}}
\put(36,-6){\makebox(0,0){g}}
\put(76,-6){\makebox(0,0){h}}
\end{picture}
\end{center}
\smallskip
\begin{center}
Figure 2: The configuration space after decimation, $\Omega'$.
\end{center}
\smallskip
\begin{center}
\begin{picture}(160,160)
\graphpaper[40](-2,-2)(160,160)
\put(74,74){\makebox(0,0){O}}
\put(82,82){\makebox(0,0){?}}
\put(-6,-6){\makebox(0,0){f}}
\put(78,-7){\makebox(0,0){g}}
\put(162,-6){\makebox(0,0){h}}
\put(34,34){\makebox(0,0){?}}
\put(74,34){\makebox(0,0){?}}
\put(114,34){\makebox(0,0){?}}
\put(-6,74){\makebox(0,0){d}}
\put(34,74){\makebox(0,0){?}}
\put(114,74){\makebox(0,0){?}}
\put(162,74){\makebox(0,0){e}}
\put(34,114){\makebox(0,0){?}}
\put(74,114){\makebox(0,0){?}}
\put(114,114){\makebox(0,0){?}}
\put(-6,162){\makebox(0,0){a}}
\put(-6,122){\makebox(0,0){?}}
\put(34,163){\makebox(0,0){?}}
\put(78,163){\makebox(0,0){b}}
\put(162,162){\makebox(0,0){c}}
\end{picture}
\end{center}
\medskip
\begin{center}
Figure 2b : The configuration space before decimation, $\Omega$.
\end{center}

To prove the failure of quasilocality, we thus have to
compute $\mu[\cdot|\mathcal{F}_{\lambda^{c}}](\omega)$ when $\omega \in T^{-1}(\omega')$,
with $\omega'$ in the neighborhood of a particular configuration in $\Omega'$.
Of course, we know that $\mu$ is a Gibbs measure for the $2d$-Ising model, so there exists
 $\Omega_\mu$ with $\mu(\Omega_\mu)=1$ s.t. for all
$\omega \in \Omega_{\mu}$, for all $\sigma \in \omega$ and $\Lambda \in \mathcal{S}$,
\begin{equation}
\mu \big[\sigma | \mathcal{F}_{\Lambda^c} \big](\omega):=\mu[\sigma |\sigma_{\Lambda^{c}}=\omega_{\Lambda^{c}}]=\frac{1}{\mathbf{Z}_{\Lambda}(\omega)} \exp\big(\sum_{\langle ij \rangle \subset \Lambda} \beta \sigma_{i}\sigma_{j}\ + \sum_{\langle ij \rangle,i \in \Lambda,j \in \Lambda^{c}} \beta \sigma_{i}\omega_{j}\big)
\end{equation}
but we want to study $\mu[\cdot |\mathcal{F}_{\lambda^{c}}]$ with
$\lambda$ \emph{non-finite}, which is not a finite-volume
probability, but on the contrary appears to be an {\em
infinite-volume Gibbs measure}, as generally proved  in
\cite{VEFS}:
\begin{lemma}\cite{VEFS} Let $\omega' \in \Omega'$ and let $\lambda$ a infinite subset of $\mathbb
Z^{2}$. Then the restriction of $\mu[\cdot
|\mathcal{F}_{\lambda^{c}}](\omega')$ to
$(\Omega_{\lambda},\mathcal{F}_{\lambda})$ is a Gibbs measure for a UAC
potential  $\Phi=\Phi(\lambda,\omega')$.
\end{lemma}

We shall not prove this lemma in the general case, but rather  directly establish the result for a
 particularly well-chosen configuration, for which a phase transition is possible for the resulting
  Gibbs measure on the smaller infinite-volume configuration space $(\Omega_\lambda, \mathcal{F}_\lambda)$.
  This particular configuration will now be in some neighborhood $\mathcal{N}_\Lambda$ of the (neutral) \emph{alternating}
  configuration
 $\omega'^{\rm{alt}}$ defined by
\begin{displaymath}
\forall i=(i_{1},i_{2}) \in \mathbb Z^{2},\;\omega'^{\rm{alt}}_{i}
,\;=(-1)^{i_{1}+i_{2}}.
\end{displaymath}

Denote $\mu^{\omega',\lambda}$  the restriction of $\mu[\cdot
|\mathcal{F}_{\lambda^{c}}](\omega')$ to
$(\Omega_{\lambda},\mathcal{F}_{\lambda})$, well-defined as a
probability measure by the existence of regular versions of
conditional probabilities as described in Chapter 3. As $\lambda$ is
fixed (we always take now $\Lambda'= \{{\bf 0}\}$), we forget it and
write $\mu^{\omega',\lambda}=\mu^{\omega'}$. To prove that it is a
Gibbs measure on $(\Omega_{\lambda},\mathcal{F}_{\lambda})$, we
consider  $\Delta \subset \lambda$ finite and a boundary condition
$\tau \in \Omega_{\lambda}$, which yields the following picture:
\smallskip
\begin{center}
\begin{picture}(240,240)
\put(30,30){\framebox(180,180)}
\put(80,225){\makebox(0,0){$\tau$}}
\put(215,30){\makebox(0,0){$\Delta$}}
\put(120,120){\makebox(0,0){O}}
\put(0,120){\makebox(0,0){$+$}}
\put(60,120){\makebox(0,0){$-$}}
\put(90,120){\makebox(0,0){$\sigma_{-1,0}$}}
\put(150,120){\makebox(0,0){$\sigma_{1,0}$}}
\put(180,120){\makebox(0,0){$-$}}
\put(240,120){\makebox(0,0){$+$}}
\put(0,150){\makebox(0,0){.}}
\put(60,150){\makebox(0,0){.}}
\put(90,150){\makebox(0,0){$\sigma_{-1,1}$}}
\put(120,150){\makebox(0,0){$\sigma_{0,1}$}}
\put(150,150){\makebox(0,0){$\sigma_{1,1}$}}
\put(180,150){\makebox(0,0){.}}
\put(240,150){\makebox(0,0){.}}
\put(0,90){\makebox(0,0){.}}
\put(60,90){\makebox(0,0){.}}
\put(90,90){\makebox(0,0){.}}
\put(120,90){\makebox(0,0){$\sigma_{0,-1}$}}
\put(150,90){\makebox(0,0){.}}
\put(180,90){\makebox(0,0){.}}
\put(240,90){\makebox(0,0){.}}
\put(0,60){\makebox(0,0){$-$}}
\put(60,60){\makebox(0,0){$+$}}
\put(90,60){\makebox(0,0){.}}
\put(120,60){\makebox(0,0){$-$}}
\put(150,60){\makebox(0,0){.}}
\put(180,60){\makebox(0,0){$+$}}
\put(240,60){\makebox(0,0){$-$}}
\put(0,180){\makebox(0,0){$-$}}
\put(30,180){\makebox(0,0){.}}
\put(60,180){\makebox(0,0){$+$}}
\put(90,180){\makebox(0,0){.}}
\put(120,180){\makebox(0,0){$-$}}
\put(150,180){\makebox(0,0){.}}
\put(180,180){\makebox(0,0){$+$}}
\put(210,180){\makebox(0,0){.}}
\put(240,180){\makebox(0,0){$-$}}
\put(0,30){\makebox(0,0){.}}
\put(240,30){\makebox(0,0){.}}
\put(0,210){\makebox(0,0){.}}
\put(30,210){\makebox(0,0){.}}
\put(60,210){\makebox(0,0){.}}
\put(90,210){\makebox(0,0){.}}
\put(120,210){\makebox(0,0){.}}
\put(150,210){\makebox(0,0){.}}
\put(180,210){\makebox(0,0){.}}
\put(210,210){\makebox(0,0){.}}
\put(240,210){\makebox(0,0){.}}
\put(0,240){\makebox(0,0){$+$}}
\put(30,240){\makebox(0,0){.}}
\put(60,240){\makebox(0,0){$-$}}
\put(90,240){\makebox(0,0){.}}
\put(120,240){\makebox(0,0){$+$}}
\put(150,240){\makebox(0,0){.}}
\put(180,240){\makebox(0,0){$-$}}
\put(210,240){\makebox(0,0){.}}
\put(240,240){\makebox(0,0){$+$}}
\put(0,0){\makebox(0,0){$+$}}
\put(30,0){\makebox(0,0){.}}
\put(60,0){\makebox(0,0){$-$}}
\put(90,0){\makebox(0,0){.}}
\put(120,0){\makebox(0,0){$+$}}
\put(150,0){\makebox(0,0){.}}
\put(180,0){\makebox(0,0){$-$}}
\put(210,0){\makebox(0,0){.}}
\put(240,0){\makebox(0,0){$+$}}
\put(0,30){\makebox(0,0){.}}
\put(30,30){\makebox(0,0){.}}
\put(60,30){\makebox(0,0){.}}
\put(90,30){\makebox(0,0){.}}
\put(120,30){\makebox(0,0){.}}
\put(150,30){\makebox(0,0){.}}
\put(180,30){\makebox(0,0){.}}
\put(210,30){\makebox(0,0){.}}
\put(30,60){\makebox(0,0){.}}
\put(30,90){\makebox(0,0){.}}
\put(30,120){\makebox(0,0){.}}
\put(30,150){\makebox(0,0){.}}
\put(30,180){\makebox(0,0){.}}
\put(30,210){\makebox(0,0){.}}
\put(210,30){\makebox(0,0){.}}
\put(210,60){\makebox(0,0){.}}
\put(210,90){\makebox(0,0){.}}
\put(210,120){\makebox(0,0){.}}
\put(210,150){\makebox(0,0){.}}
\put(210,180){\makebox(0,0){.}}
\put(210,210){\makebox(0,0){.}}
\put(75,90){\makebox(0,0){$\beta$}}
\put(73,90){\vector(-1,0){11}}
\put(77,90){\vector(1,0){11}}
\put(125,135){\makebox(0,0){$\beta$}}
\put(120,135){\vector(0,-1){11}}
\put(120,135){\vector(0,1){11}}
\end{picture}
\end{center}
\begin{center}
Figure 3: Configuration space $\Omega_{\lambda}$ with the alternating configuration in $\lambda^{c}$.
\end{center}
To check the D.L.R. equations for a suitable interaction, we have to compute, for $\omega'$ in the
 neighborhood $\mathcal{N}_\Lambda$ of $\omega'^{\rm{alt}}$, for   $\mu^{\omega'}$-a.e. $\tau$ and
 for all $\sigma_{\lambda} \in \Omega_{\lambda}$,
\begin{displaymath}
\mu^{\omega'}[\sigma |\mathcal{F}_{\lambda \backslash \Delta}](\tau) = \mu^{\omega'}[\sigma_{\lambda}|\sigma_{\lambda \backslash \Delta} =\tau_{\lambda \backslash \Delta}]=\sum_{\sigma_{\lambda^{c}} \in \Omega_{\lambda^{c}}} \mu[\sigma|\sigma_{\lambda \backslash \Delta}=\tau_{\lambda \backslash \Delta},\; \sigma_{\lambda^{c}}=\omega_{\lambda^{c}}]
\end{displaymath}
 with $\omega \in T^{-1}(\omega')$. We first assume that $\Lambda$ is big enough to contain the $\Delta$
 considered, in order to describe the resulting interaction. Later on, we shall take the infinite-volume
  for $\Delta$ and eventually encounter a volume where $\omega'$ is different to the alternating configuration, to eventually
  select different phases when possible, but for the moment the large $\Lambda$ allows us to work with the
  single $\omega'^{\rm alt}$ only but avoiding preventing possible conditioning with sets of measure zero.
   In the previous sum, only one term is not zero, when
$\sigma_{\lambda^{c}}=\omega_{\lambda^{c}}$, which is the alternating
configuration on $\lambda^{c}$. Hence,
\begin{displaymath}
\mu^{\omega'}[\sigma_{\lambda}|\sigma_{\lambda \backslash \Delta}=\tau_{\lambda \backslash \Delta}]=\mu[\sigma|\mathcal{F}_{\Delta^{c} \cup \lambda^{c}}](\tau_{\lambda} \omega_{\lambda^{c}}).
\end{displaymath}
But $\Delta^{c} \cup \lambda^{c}=(\Delta \cap \lambda)^{c}$, and
$\Delta \cap \lambda=\Delta$ is a finite subset of $\mathbb Z^{2}$,
so we now use the D.L.R. equations for $\mu$ to get, for $\mu$ -a.e.
$\tau_{\lambda} \omega_{\lambda^{c}} \in \Omega$:
$$
\mu[\sigma_{\lambda}|\sigma_{\Delta^{c}}= \tau_{\Delta_{c}}]= \frac{1}{\mathbf{Z}_{\Delta}^{\omega'}(\tau)} \exp\Big(\beta (\! \sum_{\langle ij \rangle \subset \Delta} \!  \! \sigma_{i}\sigma_{j}\ + \! \sum_{\langle ij \rangle,i \in \Delta,j \in \lambda^{c}} \! \! \sigma_{i}\omega_{j} + \! \sum_{\langle ij \rangle,i \in \Delta,j \in \lambda \cap \Delta^{c}} \! \! \sigma_{i}\tau_{j})\Big)
$$
where the normalization is the standard partition function. A
remarkable fact now, that will eventually lead us to consider the
Ising model on the decorated lattice, is that in the sum
$\sum_{\langle ij \rangle,i \in \Delta,j \in
  \lambda^{c}}J\sigma_{i}\omega_{j}$, the $j$'s are ''even'', i.e. $j=2k$
with $k \in \mathbb Z^{2}$ such that $\omega_{j}=\omega'_{k}$ is fixed in
the alternating configuration.  Then, we obtain
  the validity of the DLR equation for $\mu^{\omega'}$-almost $\tau \in
  \Omega_{\lambda}$, i.e  we have proved the :
\begin{lemma} Let $\omega'$ be the alternating configuration defined above and
assume that there exists $\omega \in T^{-1}(\omega')$ for which the D.L.R. equations
 for $\mu$ are valid with a n.n. potential appearing. Then $\mu^{\omega'}$, the restriction of
$\mu[\cdot |\mathcal{F}_{\lambda^{c}}](\omega)$ on
$(\Omega_{\lambda},\mathcal{F}_{\lambda})$ is a Gibbs measure for some
UAC potential.
\end{lemma}We do not need to give explicitly the potential but it will appear to be equivalent
to an Ising potential on the decorated lattice in the the
computation of the magnetization that we perform now. We shall then
observe (Figure 4) that the coupling due to ''even'' sites cancels
and we obtain a Gibbs measure for an Ising model on
$(\Omega_{\lambda},\mathcal{F}_{\lambda}$), with the same definition
of the nearest-neighbors as in $\mathbb Z^{2}$. To achieve this, we
just need to know that there is \emph{some} Gibbs measure for the
interaction of the previous equation. In case of phase transition,
we \emph{do not} know which it could be,
 and we shall prove that local variations in $\omega'$ could change drastically the selected phase.
  This will yield a non-Gibbsianness of the decimated measure.\\

{\sl Computation of the magnetization}\\

To prove a non quasilocality of $\nu$ at sufficiently low temperature, consider then the action of the conditional
probabilities on a local function chosen to
be characteristic of the phase transition mentioned
above. Namely, it should be an \emph{order parameter of the phase
  transition}\footnote{In statistical mechanics, an order parameter of
  a UAC potential which admit a family $\{\mu_j,j
  \in J\}$ of distinct Gibbs measures is a finite system
  $\{f_1,\ldots,f_n\}$ of local functions which discriminate these
  Gibbs measures by means of the associated expectations
  $\{\mu_j[f_1],\ldots,\mu_j[f_n]\}$ \cite{Ge}.} and
consider here the so-called \emph{magnetization}. Denote again its
origin by ${\mathbf 0}$ or $(0,0)$ and consider the local function $f \colon \Omega' \longrightarrow \mathbf{R};
\sigma' \longmapsto f(\sigma')=\sigma'_{\mathbf 0}$
and  to study $\nu[\sigma'_{{\mathbf 0}}|\mathcal{F}_{\Lambda'^{c}}](\omega')$ for  $\omega'$ in the neighborhood
of the alternating configuration, considering first that $\Lambda$ is big enough to feel $\omega'$ as the alternating
configuration itself . Then, $\nu$-a.s.  \begin{displaymath}
\nu[\sigma'_{\mathbf 0}|\mathcal{F}_{\Lambda'^{c}}](\omega')=\mu^{\omega'}[\sigma_{\mathbf 0}]
\end{displaymath}
as described in the previous section. We know that it is a Gibbs
  measure for some interaction, then by Lemma \ref{selextbc} there exists a
  sequence ($\nu_{R} \gamma_{\Lambda_{R}})_{R \in \mathbb N}$ whose weak
  limit is $\mu^{\omega'}$. For $R \in \mathbb{N}$, write, by a slight abuse of notation,
  $\Lambda_{R}$ be the intersection between $\lambda$ and the usual cube  of length $2R$ s.t.
   their exists a sequence $\nu_{R}$ with
\begin{displaymath}
\langle \sigma_{{\bf 0}} \rangle^{\omega'}:=\mu^{\omega'}[\sigma_{{\bf 0}}]=\lim_{R \rightarrow \infty}\langle \sigma_{{\bf 0}}\rangle^{\omega',\nu_{R}},
\end{displaymath}
where
\begin{displaymath}
\langle \sigma_{{\bf 0}} \rangle^{\omega',\nu_{R}}:=\int_{\Omega} \mu^{\omega'}[\sigma_{{\bf 0}}|\mathcal{F}_{\Lambda^{c}_{R}}](\tau_{R})d\nu_{R}[\tau_{R}]
\end{displaymath}
is the expectation of the spin at the origin when the boundary condition
which selects $\mu^{\omega'}$ has the law $\nu_{R}$.
Let us first fix one boundary condition $\tau_{R}$ and note $\langle
\cdot \rangle^{\omega',\tau_{R}}$ the expectation under the measure
$\mu^{\omega'}[\cdot |\mathcal{F}_{\Lambda^{c}_{R}}](\tau_{R})$. We
know that $\mu^{\omega'}$ is a Gibbs measure on
$(\Omega_{\lambda},\mathcal{F}_{\lambda})$, whose lattice  consists of all the non-even spins plus the origin.

In
order to study this measure on a more conventional lattice, let us
try to fix the spin at the origin. Define $L_{R}=\{i \in \Lambda_{R}$
s.t. $i_{1}\; {\rm and} \; i_{2}$ are both odd $\}$ and $H_{R}=\Lambda_{R}
\backslash L_{R}$. We have, using the notation $\kappa_{R}^{\lambda}(d\sigma_{\lambda})=\rho_{\Lambda_{R}}\otimes \delta^{\otimes \lambda \backslash \Lambda_{R}}_{\tau_{\lambda \backslash \Lambda_{R}}}(d\sigma_{\lambda})$,
\small
\begin{equation}\label{33}
\langle \sigma_{\mathbf{0}} \rangle^{\omega',\tau_{R}} =\frac{1}{\mathbf{Z}^{\omega',\tau_{R}}}\int_{\Omega_{\lambda}}\sigma_{\mathbf{0}}
e^{\beta(\sigma_{\mathbf{0}}-1)(\sum_{\langle i{\mathbf 0} \rangle} \sigma_{i})}e^{\sum_{\langle ij \rangle,i \in \Lambda_{R},j \in \lambda \backslash \Lambda_{R}} \beta \sigma_{i} \tau_{j}} \prod_{a \in L_{R}} (e^{\sum_{\langle ia \rangle \subset \Lambda_{R}} \beta \sigma_{a} \sigma_{i}})  \kappa_{R}^{\lambda}(d\sigma_{\lambda})
\end{equation}

\normalsize

where $\sum_{\langle i{\mathbf 0} \rangle}$ means that the sum is taken over all the spins attached to
 the origin and $\mathbf{Z}^{\omega',\tau_{R}}$ is a standard normalization. Using Fubini's theorem for positive functions,
 we integrate out  w.r.t.
the origin first to get (with $\lambda^{\star}=\lambda$ and
$\Lambda_{R}^{\star}=\Lambda_R \backslash \{{\mathbf 0}\}$)
\small
\begin{displaymath}
\langle \sigma_{{\mathbf 0}} \rangle^{\omega',\tau_{R}}=
\frac{1}{\mathbf{Z}^{\omega',\tau_{R}}}\Big(1-\int_{\Omega_{\lambda^{\star}}}e^{-2 \beta (\sum_{\langle
    i{\mathbf 0} \rangle} \sigma_{i})} e^{\sum_{\langle ij \rangle,i \in \Lambda_{R},j \in \lambda^{\star} \backslash \Lambda_{R}}\beta  \sigma_{i} \tau_{j}} \prod_{a \in L_{R}} (e^{\sum_{\langle ia \rangle \subset \Lambda_{R}} \beta \sigma_{a} \sigma_{i}}) \kappa_{R}^{\lambda^{\star}}(d\sigma_{\lambda^{\star}})\Big)
\end{displaymath}
\normalsize
with  the partition function
\begin{displaymath}
\mathbf{Z}^{\omega',\tau_{R}}=
1+\int_{\Omega_{\lambda^{\star}}}e^{-2 \beta (\sum_{\langle
i{\mathbf 0} \rangle} \sigma_{i})}  e^{\sum_{\langle ij \rangle,i \in \Lambda_{R},j \in \lambda \backslash \Lambda_{R}}\beta  \sigma_{i} \tau_{j}}\prod_{a \in L_{R}} \big(e^{\sum_{\langle ia \rangle \subset \Lambda_{R}} \beta \sigma_{a} \sigma_{i}}\big) \kappa_{R}^{\lambda^{\star}}(d\sigma_{\lambda^{\star}})
\end{displaymath}

 Hence, we only have to compute the expectation of
$e^{-2\beta (\sum_{\langle i{\mathbf 0} \rangle} \sigma_{i})}$
w.r.t. the Gibbs distribution with boundary condition $\tau_{R}$
for an Ising model on
($\Omega_{\lambda^{\star}},\mathcal{F}_{\lambda^{\star}})$ when
the spin is fixed to be ''$+$'' at the origin. We obtain this
model because of the very particular interaction we get with the
alternating configuration: The contributions of the ''even
sites'', which are fixed in the alternating configuration, cancel
each other. We have then the alternating configuration everywhere
on $2 \mathbb Z^{2}$ and an Ising distribution on the so-called
\emph{decorated lattice} $\lambda^{\star}$, without external
magnetic field as soon as $\Delta \subset \Lambda$.
 Denote $\mu^{+,\omega',\tau_{R}}$ this measure and $\langle \cdot \rangle^{+,\omega',\tau_{R}}$
 the expectation with respect to it, to get
\begin{equation}
\langle\sigma_{{\bf 0}}\rangle^{\omega',\tau_{R}}= \frac{1-\langle e^{-2 \beta (\sigma_{0,1}+\sigma_{1,0}+\sigma_{-1,0}+\sigma_{0,-1})}\rangle^{+,\omega',\tau_R}}{1+\langle e^{-2\beta(\sigma_{0,1}+\sigma_{1,0}+\sigma_{-1,0}+\sigma_{0,-1})}\rangle^{+,\omega',\tau_R}}.
\end{equation}
To get a more standard expression in terms of standard Ising
models,
 we now use the following trick, standard in statistical mechanics with $\pm 1$
 Ising spins, to reduce  $\langle\sigma_{0,1}\rangle^{+,\omega',\tau_{R}}$, the expectation of one spin attached to the origin:
\begin{lemma}
\begin{eqnarray*}
\langle \sigma_{0,1}\rangle^{+,\omega',\tau_{R}}&=&\Big \langle\tanh(J(\sigma_{1,1} + \sigma_{-1,1}))\Big\rangle^{+,\omega',\tau_{R}}=\Big\langle(\frac{1}{2} \tanh(2J))(\sigma_{1,1} + \sigma_{-1,1})\Big\rangle^{+,\omega',\tau_{R}}
\end{eqnarray*}
where $\sigma_{1,1}$ and $\sigma_{-1,1}$ are the spins attached to
$\sigma_{0,1}$.
\end{lemma}

This reduces our study to  the distribution of the spins in
$L_{R}$, that is in fact
 the {\em decorated lattice}, the lattice of spins whose coordinates are both odd.\\

\begin{center}
\begin{picture}(240,240)
\put(30,30){\framebox(180,180)}
\put(215,25){\makebox(0,0){$\Lambda_{R}$}}
\put(120,120){\makebox(0,0){$+$}}
\put(0,120){\makebox(0,0){$+$}}
\put(60,120){\makebox(0,0){$-$}}
\put(90,120){\makebox(0,0){.}}
\put(150,120){\makebox(0,0){.}}
\put(180,120){\makebox(0,0){$-$}}
\put(240,120){\makebox(0,0){$+$}}
\put(0,150){\makebox(0,0){.}}
\put(60,150){\makebox(0,0){.}}
\put(120,150){\makebox(0,0){$\sigma_{0,1}$}}
\put(150,150){\makebox(0,0){$\sigma_{1,1}$}}
\put(180,150){\makebox(0,0){.}}
\put(240,150){\makebox(0,0){.}}
\put(0,90){\makebox(0,0){.}}
\put(60,90){\makebox(0,0){.}}
\put(90,90){\makebox(0,0){.}}
\put(120,90){\makebox(0,0){.}}
\put(150,90){\makebox(0,0){.}}
\put(180,90){\makebox(0,0){.}}
\put(240,90){\makebox(0,0){.}}
\put(0,60){\makebox(0,0){$-$}}
\put(60,60){\makebox(0,0){$+$}}
\put(90,60){\makebox(0,0){.}}
\put(120,60){\makebox(0,0){$-$}}
\put(150,60){\makebox(0,0){.}}
\put(180,60){\makebox(0,0){$+$}}
\put(240,60){\makebox(0,0){$-$}}
\put(0,180){\makebox(0,0){$-$}}
\put(30,180){\makebox(0,0){.}}
\put(60,180){\makebox(0,0){$+$}}
\put(90,180){\makebox(0,0){.}}
\put(120,180){\makebox(0,0){$-$}}
\put(150,180){\makebox(0,0){.}}
\put(180,180){\makebox(0,0){$+$}}
\put(210,180){\makebox(0,0){.}}
\put(240,180){\makebox(0,0){$-$}}
\put(0,30){\makebox(0,0){.}}
\put(240,30){\makebox(0,0){.}}
\put(0,210){\makebox(0,0){.}}
\put(30,210){\makebox(0,0){.}}
\put(60,210){\makebox(0,0){.}}
\put(90,210){\makebox(0,0){.}}
\put(120,210){\makebox(0,0){.}}
\put(150,210){\makebox(0,0){.}}
\put(180,210){\makebox(0,0){.}}
\put(210,210){\makebox(0,0){.}}
\put(240,210){\makebox(0,0){.}}
\put(0,240){\makebox(0,0){$+$}}
\put(30,240){\makebox(0,0){.}}
\put(60,240){\makebox(0,0){$-$}}
\put(90,240){\makebox(0,0){.}}
\put(120,240){\makebox(0,0){$+$}}
\put(150,240){\makebox(0,0){.}}
\put(180,240){\makebox(0,0){$-$}}
\put(210,240){\makebox(0,0){.}}
\put(240,240){\makebox(0,0){$+$}}
\put(0,0){\makebox(0,0){$+$}}
\put(30,0){\makebox(0,0){.}}
\put(60,0){\makebox(0,0){$-$}}
\put(90,0){\makebox(0,0){.}}
\put(120,0){\makebox(0,0){$+$}}
\put(150,0){\makebox(0,0){.}}
\put(180,0){\makebox(0,0){$-$}}
\put(210,0){\makebox(0,0){.}}
\put(240,0){\makebox(0,0){$+$}}
\put(0,30){\makebox(0,0){.}}
\put(30,30){\makebox(0,0){.}}
\put(60,30){\makebox(0,0){.}}
\put(90,30){\makebox(0,0){.}}
\put(120,30){\makebox(0,0){.}}
\put(150,30){\makebox(0,0){.}}
\put(180,30){\makebox(0,0){.}}
\put(210,30){\makebox(0,0){.}}
\put(30,60){\makebox(0,0){.}}
\put(30,90){\makebox(0,0){.}}
\put(30,120){\makebox(0,0){.}}
\put(30,150){\makebox(0,0){.}}
\put(30,180){\makebox(0,0){.}}
\put(30,210){\makebox(0,0){.}}
\put(210,30){\makebox(0,0){.}}
\put(210,60){\makebox(0,0){.}}
\put(210,90){\makebox(0,0){.}}
\put(210,120){\makebox(0,0){.}}
\put(210,150){\makebox(0,0){.}}
\put(210,180){\makebox(0,0){.}}
\put(210,210){\makebox(0,0){.}}
\put(90,180){\makebox(0,0){$b$}}
\put(90,210){\makebox(0,0){$b'$}}
\put(90,150){\makebox(0,0){$b''$}}
\put(75,225){\makebox(0,0){$\tau$}}
\put(30,90){\makebox(0,0){$a$}}
\put(90,90){\makebox(0,0){$a'$}}
\put(150,90){\makebox(0,0){$b'$}}
\put(180,90){\makebox(0,0){$b$}}
\put(210,90){\makebox(0,0){$b''$}}
\end{picture}
\end{center}
\begin{center}
Figure 4 : Ising model on the decorated lattice $\lambda^{*}$
\end{center}

We then have to compute
$\langle\sigma_{1,1}\rangle^{+,\omega'\tau_{R}}$. As claimed
before, we can start integration with respect to the spins in
$H_{R}$, the sites of the decorated lattice which have exactly two
neighbors. We call $H^{0}_{R}=H_{R} \backslash \Gamma_{R}$ where
$\Gamma_{R}=\Lambda_R \backslash \Lambda_{R-1}$ is the boundary of
$\Lambda_{R}$. The sites in $H^{0}_{R}$ are those which have two
neighbors \emph{in} $\Lambda_{R}$. We also call $H^{1}_{R}=H_{R}
\cap \Gamma_{R}$ the set of the sites which have two neighbors in
the lattice $\lambda^{\star}$, one in $\Lambda_{R}$ and the other,
fixed by the boundary condition $\tau$, outside $\Lambda_{R}$.
Compute:
\begin{equation}\label{reduce}
\langle\sigma_{1,1}\rangle^{+,\omega'\tau_{R}}=
\frac{1}{\mathbf{Z}^{+,\omega',\tau_{R}}}
\int_{\Omega_{\lambda^{\star}}} \sigma_{1,1} \cdot
A^{0}_{R}(\sigma,d\sigma_{H_{R}^{0}}) \cdot
A^{1}_{R}(\sigma,d\sigma_{H_{R}^{1}}) \cdot
A_{R}(\sigma,d\sigma_{L_R})
\end{equation}
where
\begin{eqnarray*}
A^{0}_{R}(\sigma,d\sigma_{H_{R}^{0}})=\prod_{b \in H^{o}_{R}} e^{\beta  \sigma_{b} (\sigma_{b'} + \sigma_{b''})} \; \rho_{0}(d\sigma_{b}) \\
A^{1}_{R}(\sigma,d\sigma_{H_{R}^{1}})=\prod_{b \in H^{1}_{R}} e^{\beta  \sigma_{b} (\sigma_{b'} + \tau_{b''})} \;\rho_{0}(d\sigma_{b}) \\
A_R(\sigma,d\sigma_{L_R})= \prod_{a \in L_{R}}\; \rho_{0}[d\sigma_{a}] \otimes \delta^{\lambda \backslash \Lambda_{R}}_{\tau_{\lambda \backslash \Lambda_R}}(d\sigma_{\lambda \backslash \Lambda_R})
\end{eqnarray*}
where  for each $b \in H_{R}$, we have called $b'$ and $b''$ its neighbors in $L_{R}$ or filled by the boundary condition $\tau$ in
$\Lambda_{R+1}$, to get   for the integral (\ref{reduce})
\begin{displaymath}
\int_{\Omega_{L_{R}}} \sigma_{1,1}\Big(\prod_{b \in H_{R}} \int_{E} e^{\beta \sigma_{b} (\sigma_{b'} + \sigma_{b''})} \rho_{0}[d\sigma_{b}]\Big) \prod_{a \in L_{R}} \rho_{0}[d\sigma_{a}] \otimes \delta^{\otimes \lambda \backslash \Lambda_{R}}_{\tau_{\lambda \backslash \Lambda_R}}[d\tau_{\lambda \backslash \Lambda_R}].
\end{displaymath}
Now, using another standard trick on Ising spins, we calculate
\begin{displaymath}
\int_{E} e^{\beta \sigma_{b} (\sigma_{b'} + \sigma_{b''})}
\rho_{0}(d\sigma_{b})=\frac{e^{\beta (\sigma_{b'} +\sigma_{b''})}
+ e^{- \beta (\sigma_{b'} +\sigma_{b''})}}{2}
\end{displaymath}
in such a way that the contribution of the spins in $H_{R}$ does
not appear in the integral anymore, because the set $\{(b',b''),b
\in H_{R}\}$ is $L_{R}$. To get a more standard Ising
representation, we would like to obtain now a coupling interaction
between the spins in $L_{R}$. To do so, write
\begin{equation}\label{trick}
\frac{e^{\beta (\sigma_{b'} +\sigma_{b''})} + e^{-\beta (\sigma_{b'} +\sigma_{b''})}}{2}=Ke^{\beta' \sigma_{b'} \sigma_{b''}}
\end{equation}
where $K$  cancels by normalization. On the event
$\{\sigma_{b'}=+1,\sigma_{b''}=+1\}$, we should have
\begin{displaymath}
\cosh[2 \beta ]=Ke^{\beta'}
\end{displaymath}
and on the events  $\{\sigma_{b'}=-1,\sigma_{b''}=+1\}$ and $\{\sigma_{b'}=+1,\sigma_{b''}=-1\}$
\begin{displaymath}
1=Ke^{-\beta'}
\end{displaymath}
then, one could take $K=e^{\beta'}$ and $e^{2\beta'}=\cosh[2\beta]$ i.e.
 $\beta'=\frac{1}{2} \cosh^{-1}\big(e^{2 \beta}\big)$ in (\ref{trick})  so
\begin{displaymath}
\langle\sigma_{1,1}\rangle^{+,\omega'\tau_{R}} \; =
\frac{1}{\mathbf{Z}^{+,\omega',\tau_{R}}} \int_{\Omega_{L_{R}}}
\Big(\sigma_{1,1} e^{\beta' \sum_{\langle aa' \rangle \subset L_{R}}
  \sigma_{a} \sigma_{a'} + \beta' \sum_{\langle aa' \rangle,a \in L_{R},a' \in \lambda \backslash \Lambda_{R}} \sigma_{a} \tau_{a'}}\Big) \rho_{L_R}[(d\sigma_{L_R}).
\end{displaymath}

\begin{center}
\begin{picture}(240,240)
\put(30,30){\framebox(180,180)}
\put(215,20){\makebox(0,0){$\Lambda_{R}$}}
\put(120,120){\makebox(0,0){$+$}}
\put(0,120){\makebox(0,0){$+$}}
\put(60,120){\makebox(0,0){$-$}}
\put(180,120){\makebox(0,0){$-$}}
\put(240,120){\makebox(0,0){$+$}}
\put(90,150){\makebox(0,0){$*$}}
\put(150,150){\makebox(0,0){$*$}}
\put(90,90){\makebox(0,0){$*$}}
\put(150,90){\makebox(0,0){$*$}}
\put(0,60){\makebox(0,0){$-$}}
\put(60,60){\makebox(0,0){$+$}}
\put(120,60){\makebox(0,0){$-$}}
\put(180,60){\makebox(0,0){$+$}}
\put(240,60){\makebox(0,0){$-$}}
\put(0,180){\makebox(0,0){$-$}}
\put(60,180){\makebox(0,0){$+$}}
\put(120,180){\makebox(0,0){$-$}}
\put(180,180){\makebox(0,0){$+$}}
\put(240,180){\makebox(0,0){$-$}}
\put(30,210){\makebox(0,0){$*$}}
\put(90,210){\makebox(0,0){$*$}}
\put(150,210){\makebox(0,0){$*$}}
\put(210,210){\makebox(0,0){$*$}}
\put(0,240){\makebox(0,0){$+$}}
\put(60,240){\makebox(0,0){$-$}}
\put(120,240){\makebox(0,0){$+$}}
\put(180,240){\makebox(0,0){$-$}}
\put(240,240){\makebox(0,0){$+$}}
\put(0,0){\makebox(0,0){$+$}}
\put(60,0){\makebox(0,0){$-$}}
\put(120,0){\makebox(0,0){$+$}}
\put(180,0){\makebox(0,0){$-$}}
\put(240,0){\makebox(0,0){$+$}}
\put(30,30){\makebox(0,0){$*$}}
\put(90,30){\makebox(0,0){$*$}}
\put(150,30){\makebox(0,0){$*$}}
\put(210,30){\makebox(0,0){$*$}}
\put(30,150){\makebox(0,0){$*$}}
\put(210,30){\makebox(0,0){$*$}}
\put(210,90){\makebox(0,0){$*$}}
\put(210,150){\makebox(0,0){$*$}}
\put(210,210){\makebox(0,0){$*$}}
\put(90,215){\makebox(0,0){$b'$}}
\put(90,145){\makebox(0,0){$b''$}}
\put(75,225){\makebox(0,0){$\tau$}}
\put(30,90){\makebox(0,0){$a$}}
\put(95,90){\makebox(0,0){$a'$}}
\put(60,100){\makebox(0,0){$\beta'$}}
\put(60,90){\vector(-1,0){26}}
\put(60,90){\vector(1,0){26}}
\put(90,180){\vector(0,1){26}}
\put(90,180){\vector(0,-1){26}}
\put(75,180){\makebox(0,0){$\beta'$}}
\end{picture}
\end{center}
\begin{center}
Figure 5 : Ising model on $2\mathbb Z^{2}$ with coupling $\beta'$.
\end{center}

It is \emph{exactly} the magnetization of a ferromagnetic Ising model
at inverse temperature $\beta'$ on $2 \mathbb Z^{2}$, with the boundary condition
$\tau$ on $(\lambda \backslash \Lambda_{R}) \cap 2\mathbb Z^{2}$ and
without external field. When the temperature is low enough, we know by Theorem \ref{thmIsing}
that a phase transition holds and that the above magnetization is an order parameter, and this
will eventually lead to essential discontinuity as soon as $\beta'> \beta_c$,
which yields $\beta  > \tilde{\beta}_c=\frac{1}{2} \cosh^{-1}\big(e^{2 \beta}\big)$.

To rigorously get the essential discontinuity,  one should now do
the same computation when $\Lambda_n$ is bigger than $\Lambda$,
i.e for distinct neighborhoods of the alternating configuration,
where all pluses or all minuses far away will create a external
field that eventually selects the different phases. The procedure
is the same but one has to be careful in some computations, to
eventually give  rise to the essential discontinuity we seek for.
Thus, this failure of quasilocality comes directly from the
presence of a phase transition in some ''hidden system'', that of
the internal spins. This is carefully proved in detail in
\cite{VEFS} in the
\begin{lemma}[essential discontinuity]
Let $\beta>\frac{1}{2} \cosh^{-1}\big(e^{2 \beta_c}\big)$ and let $\omega'^{alt}$ be the alternating
configuration. $\forall \epsilon > 0, \; \forall \mathcal{N}$ neighborhood of $\omega'^{alt}$ , $\exists  R_{o}>0$ such that $\forall R>R_o$ , we can find $\mathcal{N}_{R,+}, \; \mathcal{N}_{R,-} \subset \mathcal{N}$ with $\nu[A_{R,+}]=\nu[A_{R,-}]>0$ and for $\nu$-a.e. $\omega'_{1} \in \mathcal{N}_{R,+}$, for $\nu$-a.e. $\omega'_2 \in \mathcal{N}_{R,-}$,
\begin{displaymath}
\nu[\sigma'_{{\bf 0}}|\mathcal{F}'_{\{{\bf 0}\}^{c}}](\omega'_1) \; - \; \nu[\sigma'_{{\bf 0}}|\mathcal{F}'_{\{{\bf 0}\}^{c}}](\omega'_2) \; > \epsilon.
\end{displaymath}
Thus, no version of the conditional probabilities of $\nu$ given
$\mathcal{F}'_{\{0\}^{c}}$ can be continuous.
\end{lemma}

This proves Theorem \ref{Dec}. This basic example expresses the
link between the pathology and the existence of a phase transition
in some ''hidden'' system. The same procedure has to be used for
more general RG transformations, but it is sometimes difficult or
even unknown to detect a bad configuration: Indeed, getting a
 bad configuration amounts to prove phase transitions, and this is sometimes, not to say often,  difficult or unknown,
  involving e.g. very sophisticated versions of the theory of Pirogov-Sinai.
  We give now a non-exhaustive small catalogue of results that have been proved
  during the last decades, many other examples are rigorously described in \cite{VEFS,VE}.

\subsection{General RG transformations, main examples and results}

In the general framework, we deal with two configuration spaces, a so-called {\em
  original} one $(\Omega,\mathcal{F},\rho)$ and a so-called {\em image}
  one $(\Omega',\mathcal{F}',\rho')$. Most of the time, the lattice $S'$
  of the image system is smaller, and of the same kind
  (e.g. $S=\mathbb{Z}^{d},\;S'=\mathbb{Z}^{d'},\; d \geq d'$); the above decimation transformation is e.g.
  sometimes described with $2\mathbb{Z}^d$ as image
  lattice, and the projection on an hyperplane or ''restriction to a layer'' is studied in this
  context of the renormalization group whereas it does not satisfy all
  the properties of our following formal definition. The extension of the definition to spaces of measures is standard.

\begin{definition}[R.G. kernels]
A  renormalization group transformation (R.G.T.) is a a probability
kernel $T$ from $(\Omega,\mathcal{F})$ to $(\Omega',\mathcal{F}')$ such that
\begin{enumerate}
\item $T$ carries $\mathcal{M}_{1,\rm{inv}}^+(\Omega)$ onto
$\mathcal{M}_{1,\rm{inv}}^+(\Omega)$. \item There exists sequences
of cubes $(\Lambda_n)_{n \in \mathbb{N}}$ and $(\Lambda'_n)_{n \in
\mathbb{N}}$, respectively finite subsets of $S$ and $S'$, such
that:
\begin{enumerate}
\item $\forall A' \in \mathcal{F}'_{\Lambda'_n}$, the function
$T(\cdot,A')$ is $\mathcal{F}_{\Lambda_n}$-measurable: The
behavior of the image spins in $\Lambda'_n$ depends only on the
original spins in $\Lambda_n$. \item $\limsup_{n \to \infty}
\frac{\mid \Lambda_n \mid}{\mid \Lambda'_n
  \mid} \; \leq K < \infty$.
\end{enumerate}
\end{enumerate}
\end{definition}

We give now a few examples where non-Gibbsianness has been proved
to arise. As already claimed, it is more illustrative than
exhaustive, a general method applying to many examples is
available in \cite{VEFS} and in related papers from our
bibliography. We distinguish two types of examples, those where
the transformation is deterministic, like the above decimation,
and the more general stochastic ones. We first extend the previous
results to more general decimations.
\subsubsection{Deterministic transformations}

A RGT is said to be {\em deterministic} when the probability kernel induced by $T$ is deterministic in the sense that
\begin{displaymath}
\forall A' \in \mathcal{F}', \forall \omega \in \Omega,
T(\omega,A')=\delta_{\omega'}(A')
\end{displaymath}
where the image $\omega'=t(\omega)$ is a function of the original
configuration $\omega$.
\begin{enumerate}

\item{{\bf  Decimation transformations in higher dimensions:}}

 It is thus a deterministic probability kernel
from $\Omega=E^{\mathbb{Z}^d}$ onto itself, with $t(\omega)=\omega'$ and $\omega'_i=\omega_{b i}$.
In the same spirit of the phenomenon observed for the 2d Ising model with the alternating configuration,
 but with much more difficult proofs in general, usually involving an heavy machinery and
 tricks to find a special configuration and to prove it is a point of essential discontinuity. Among others, one gets
\begin{theorem}\cite{VEFS}
Let $d \geq 2$ and $b \geq 2$. Then for all $\beta > \beta(d,b)$ sufficiently large,
for any Gibbs measure $\mu$ for the standard n.n. homogeneous Ising model on $\Z^d$ with coupling $J>0$ and magnetic
 field $h=0$, the decimated measure $\nu=T_b \mu$ is not quasilocal.
\end{theorem}

This result is also extended in some open region $(\beta,h)$ of the phase diagram, e.g.
 to small magnetic field at dimension $d \geq 3$, for an adapted special configuration, see Section 4.3.6. in \cite{VEFS}.
  It is also interesting in view of the historical aim of renormalization group, that iterates the transformations to reach
  the presumably fixed critical point, that the quasilocality property could be recovered after iterating this decimation
   transformation \cite{VEL}, and that other positive results on conservation on Gibbisanness in other parts
   of the phase diagram exist \cite{HKen}.
\item{{\bf  Deterministic majority-rule transformation for the Ising
  model:}}

The configuration spaces are still identical and are those of
the $d$-dimensio\-nal
Ising model $\Omega'=\Omega=\{-1,+1\}^{\mathbb{Z}^{d}},d\geq 1$. Let
  $b \geq 1$ be an integer and let $B_0 \in \mathcal{S}$ with $|
  B_0 |$ {\em odd}. Define, $\forall i \in \mathbb{Z}^{d}, \; B_i$
  to be $B_0$ translated by $b \cdot i$: $B_i=B_0+b \cdot i$. We call this subsets
  of $\mathbb{Z}^{d}$ {\em blocks}. The deterministic kernel
is the transformation $t(\omega)=\omega'$  defined by
\begin{displaymath}
\forall i \in \mathbb{Z}^{d}, \omega'_i=\frac{\sum_{j \in
    B_i}\omega_j}{\mid \sum_{j \in B_i}\omega_j \mid}.
\end{displaymath}

For these transformations, getting some bad configurations leading
to non-quasilocality is sometimes difficult, due to the constraint
it gives on the block, see \cite{VEFS}. One nevertheless proves
the
\begin{theorem}\cite{VEFS}
Let $J$ large enough, $\mu$ any Gibbs measure for the $2d$ Ising model with n.n.
 coupling $J$ and zero magnetic field. Let $T$ be the majority with blocks of size $|B_0|=7$.
 Then $\nu=\mu T$ is not quasilocal.
\end{theorem}
The result has been extended for smaller blocks using a
computer-assisted proof in \cite{Ken}. Such transformations belong
to a more general family of {\em block-spins transformations},
very useful in renormalization procedures or multi-scale analysis,
see e.g \cite{BCO,Sido}.

\item{{\bf Modified majority-rule on a Cayley tree with
overlapping blocks}:}

The Ising model on a Cayley tree has been introduced in the
previous chapter. We shall restrict ourself to the simplest
rooted-Cayley tree $\mathcal{T}^{2}_{0}$ \cite{BLG}  and let $\mu$
be any Gibbs measure for this model (we have seen in the previous
chapter that there always exists at least one Gibbs measure for
this model). We choose the root as the origin and we denote it
$r$. Define $\Omega=\{-1,+1\}^{\mathcal{T}^{2}_{0}}$  and
$\Omega'=\{-1,0,+1\}^{\mathcal{T}^{2}_{0}}$, and let $R$ be any
non negative integer to  define the closed ball of $S$ of radius
$R$ to be $V_R=\{i \in \mathcal{T}^{2}_{0} \mid d(r,i) \leq R \}$.
Denote also  its boundary by $W_R=\{i \in \mathcal{T}^{2}_{0} \mid
d(r,i) =R \}$ where $d$ is the canonical metric on
$\mathcal{T}^{2}_{0}$. We shall represent the vertices of
$\mathcal{T}^{2}_{0}$ by sequences of bits, defined by recurrence:
The representation of the origin $r$  is the void binary sequence,
and that of its neighbors are chosen to be $0$ and $1$. Now let
$R>0$ and let $i \in W_R$ with representation $i^{*}$. There are
only two sites $k$ and $l$ in $W_{R+1}$ at distance 1 from $i$. We
define then their representation to be $k^*=i^*0$ and $l^*=i^*1$.
We obtain a representation of all the vertices of
$\mathcal{T}^{2}_{0}$. We shall now write the same symbol $i$ for
the vertex or the binary representation $i^*$. Define
$C_r=\{r,0,1\}$and $\forall j \in \mathcal{T}^{2}_{0},\; j \neq
r$, the cell
\begin{displaymath}
C_j=\{j,j0,j1\}
\end{displaymath}
where $j0$ and $j1$ are the two neighbors of $j$ from the
''following'' level. For example, $C_0=\{0,00,01\}$. Define as
well $c_j=\textrm{Card}(C_j)$, with here $c_j=c=3$, and consider
now the deterministic transformation $t :  \omega \longmapsto
t(\omega)=\omega'$ where $\omega'$ is defined by
\begin{displaymath}
\omega'_j=\left\{
\begin{array}{lll}\; +1 \; & \textrm{iff}& \frac{1}{c}\sum_{i \in C_j}
\omega_i=+1\\
\; 0 \; & \textrm{iff}&\frac{1}{c} \mid \sum_{i \in C_j}
\omega_i \mid < 1\\
\; -1 \;& \textrm{iff}& \frac{1}{c}\sum_{i \in C_j} \omega_i=-1.
\end{array} \right.
\end{displaymath}

This could be seen as a static version of the voter model, where a
child votes like its parents when they agree. In this case, due to
the overlapping of the blocks, the failure of quasilocality occurs
at
 all temperatures, for the very particular null everywhere configuration. An interesting fact is that the
 set of bad configurations is topologically rather big but suspected to be of zero DLR-measure \cite{LN}.

\begin{theorem}\cite{LN}
Let $\mu$ be any Gibbs measure for the Ising model on
$\mathcal{T}^{2}_{0}$ and let $\nu$ be the image of $\mu$ by $T$. Then
$\nu$ is non quasilocal at any temperature and cannot be a Gibbs
measure.
\end{theorem}

It could be generalized
to non-rooted Cayley trees, and with other sizes of blocks, also as a stochastic transformation,
 modelling the fact that a child does not always vote like its parents,  as follows:
 Let $\epsilon \in [0,1]$ and $\xi$ be a Bernoulli random
variable with parameter $\epsilon$. Define the deterministic map
$t_{\epsilon} : \omega \longmapsto t(\omega)=\omega'$ where
$\omega'$ is
 defined for all $j \in \mathcal{T}_0^2$ by
\begin{displaymath}
\omega'_j=\left\{
\begin{array}{lll}\; +1 \; & \textrm{iff}& \frac{1}{c}\sum_{i \in C_j}
\omega_i=+1 \; \textrm{and} \; \xi = 1 \\
\; -1 \;& \textrm{iff}& \frac{1}{c}\sum_{i \in C_j}
\omega_i=-1 \; \textrm{and} \; \xi = 1 \\
\; 0 \; & \textrm{if} & \; \xi=0.
\end{array} \right.
\end{displaymath}

Its action is described by a probabilistic kernel $T_{\epsilon}$ defined by:
\begin{displaymath}
\forall A' \in \mathcal{F}', \forall \omega \in \Omega, T_{\epsilon}(\omega,A')=(1-\xi)\delta_{t_{\epsilon}(\omega)}(A')+ \xi\delta_{0}(A') .
\end{displaymath}
 It could be interesting to study the difference between the deterministic
 transformation and the stochastic ones because this could play a role on the degree
  of non-Gibbsianness of the image measure.

\item{{\bf Restriction of Ising model to a
layer}} \cite{Maes3,sch}:

This transformation is not properly speaking a R.G.T. in the sense
of our definition, because of the lack of strict locality.
Nevertheless, it is known to lead to non-Gibbsianness (see
\cite{sch}) and is a good example of a new kind of random fields,
the weakly Gibbsian measures, which will be introduced soon
\cite{Maes3}. The configuration spaces are
$\Omega=\{-1,+1\}^{\mathbb{Z}^{2}}$ and
  $\Omega'=\{-1,+1\}^{\mathbb{Z}}$. The transformation is
    deterministic and defined by $t(\omega)=\omega'$
where $\omega'$ is defined by $\forall i \in \mathbb{Z}^{d-1},
\omega'_i=\omega_{(i,0)}$ where $0$ denotes here the origin in
$\mathbb{Z}$. The interesting fact in this example, together with
the fact that it could seem very natural at a first sight to
consider the projected measure to be Gibbs, is that the original
proof relies on {\em wrong large deviation properties} of the
projected measure.
\end{enumerate}

\subsubsection{Stochastic transformations}
In contrast to the deterministic case, a stochastic
transformation could lead to different image configuration,  with a certain probability for each. We have already
seen an example of stochastic R.G.T. on the tree.

\begin{enumerate}
\item{{\bf Stochastic majority-rule for Ising model}:}

The definition is very similar to the deterministic one, except that
we deal with blocks $B_0$ with $| B_0 |$
even. The configuration spaces are
still $\Omega'=\Omega=\{-1,+1\}^{\mathbb{Z}^{d}},d\geq 1$. Let
  $b \geq 1$ be an integer and let $B_0 \in \mathcal{S}$ with $|
  B_0 |$ {\em even}. Define, $\forall i \in \mathbb{Z}^{d}, \; B_i$
  to be $B_0$ translated by $b \cdot i$: $B_i=B_0+b \cdot i$. Let $\xi$ be a
  Bernoulli random variable on $(\Omega,\mathcal{F})$ with parameter
  $p$. Most of the time, $p$ is considered to be $\frac{1}{2}$. The
  stochastic majority-rule is the transformation $T$ which transforms
  $\omega$ in $t(\omega)=\omega'$ with

\begin{displaymath}
\omega'_i=\left\{
\begin{array}{lllll}\; +1 \; & \textrm{if}&\sum_{j \in B_i}
\omega_j > 0 \\
\; -1 \;& \textrm{if}& \sum_{j \in B_i}
\omega_j < 0 \\
\; +1 \; & \textrm{if}& \sum_{j \in B_i}
\omega_j = 0 \; \textrm{and} \; \xi=+1\\
\; -1 \; & \textrm{if}& \sum_{j \in B_i}
\omega_j = 0 \; \textrm{and} \; \xi=0.\\
\end{array} \right.
\end{displaymath}

\item{{\bf Kadanoff transformations for the Ising model}:}

These transformations model a lot of interesting and historical
R.G.T. We shall not deal with them, but some are widely studied in
\cite{GP,Is2,VEFS}. Here again the blocks $B_i$ are defined
in the same way for $i \in \mathbb{Z}^d$, the configuration spaces are
$\Omega=\Omega'=\{-1,+1\}^{\mathbb{Z}^{d}}$ and $p$ is a strictly
  positive real.  The R.G.T. map is defined
  atom per atom by
\begin{displaymath}
T(\omega,\omega')=\prod_{i \in
  \mathbb{Z}^d}\frac{\exp(p\omega'_i\sum_{j \in B_i}\omega_j)}{2
  \cosh(p\sum_{j \in B_i} \sigma_j)}
\end{displaymath}

This transformation is also associated with stochastic evolutions
of Gibbs measures, as we shall see. They have been proved to lead
to non-Gibbsianness for $d \geq  2, b \geq 1$ and $p$ finite. It
also includes majority rules or decimations in the limit $p \to
\infty$ for suitable blocks.
\end{enumerate}
Many other examples are available in the literature, and as
claimed in \cite{VE}, the surprise is eventually not that they are
non-Gibbsian, but that it took so long to realize it, the set of
Gibbs measures being topologically very small \cite{Is3}. In the
same seminal paper \cite{VEFS},  positive general results are
given about the action of these transformations on Hamiltonians
and potentials, excluding various scenarii related to figure 1. We
only quote them, see the discussions in \cite{F} and \cite{VEFS}.

\subsubsection{Renormalization transformation on potentials}

As explained in the beginning of this chapter, the extension of
the renormalization transformations to potentials is not always
well-defined, whereas the extension of the R.G.T. to an action on
measures is standard and always possible. Nevertheless, two
positive results have been proved by van Enter {\em et al.}, and
we introduce them before describing the pathologies of the
renormalization group. We restrict ourselves to a space
$\mathcal{B}^{1}$ consisting of the translation-invariant,
continuous and uniformly absolutely convergent potentials. We
introduce first a relation instead of a function:

\begin{definition}[R.G.T. on interactions]
Let $T$ be a R.G.T. We define a renormalization group relation
$\mathcal{R}=\mathcal{R}_{T}$ on interactions by the relation
\begin{displaymath}
\mathcal{R}=\{(\Phi,\Phi') \in \mathcal{B}^{1} \times \mathcal{B}^{1}:
  \; \exists \mu \; \textrm{translation-invariant in} \;
  \mathcal{G}(\gamma^{\Phi}) \; \textrm{s.t.} \; \mu T \in
  \mathcal{G}(\gamma^{\Phi'})\}
\end{displaymath}
Where $\mu T$ is the image measure of $\mu$ by $T$.
\end{definition}

The next theorem tells us that $\mathcal{R}$ is single-valued and is
proved in \cite{VEFS}.

\begin{theorem}[first fundamental theorem of the renormalization
  group]
Let $\mu$ and $\nu$ be translation-invariant Gibbs measures with
respect to the same interaction $\Phi \in \mathcal{B}^{1}$ and let $T$
be an R.G-Transformation. The following results are true:
\begin{enumerate}
\item Either $\mu T$ and $\nu T$ are both non-quasilocal, or else there
exists a quasilocal specification $\gamma'$ with which both $\mu T$
and $\nu T$ are consistent.
\item Either $\mu T$ and $\nu T$ are both non-Gibbsian, or else there exists
a uniformly absolutely convergent potential $\Phi'$ for which both
$\mu T$ and $\nu T$ are Gibbs measures.
\end{enumerate}
\end{theorem}


\subsection{Stochastic evolution of Gibbs measures \cite{EFHR}}
Once these RG pathologies have been identified as the manifestation of non-Gibbsianness,
a natural source of examples to be investigated to detect similar phenomena concern stochastic Ising models,
introduced in the previous chapter. Of course, here it could not come out as a surprise due to the equilibrium
 considerations that led to the introduction of the Gibbs property: It should be natural to encounter such
 a phenomenon in the course of stochastic evolutions of Gibbs measures, and it is indeed the case during the
  heating of a low temperature Ising model, i.e. the stochastic evolution of a low temperature Gibbs measure
   for the ferromagnetic n.n. Ising model during a high temperature Glauber dynamics. It is not so simple to
   establish, and not always true; physical interpretations can be found in \cite{PO}. Nevertheless, not much
   is known in non-equilibrium statistical mechanics, so any information about Gibbsianness in transient regimes is welcome.
    The first systematical study of such phenomena has been made in \cite{EFHR}, although similar  investigations had been
    made earlier in \cite{Maes1,Lebsch}. Before describing a bit more the relationships with (stochastic) RG-transformations
    through the description of the an infinite-temperature Glauber dynamics, let us quote their general result.
\begin{theorem}\cite{EFHR}
Let $\Phi$ be a translation-invariant potential, $\mu$ a  corresponding translation-invariant Gibbs measure at inverse
 temperature $\beta$ and $S(t)$ the semi-group corresponding to the dynamics having $\mu$ has reversible measure. Denote
  by $\nu$ an initial translation-invariant distribution of the a priori configuration. Then
\begin{enumerate}
\item For all $\nu,\mu$, the time-evolved measure $\nu S(t)$ is
Gibbs for small times $t \leq t_0(\beta)$. \item If $\mu,\nu$
corresponds to high or infinite temperature Gibbs measures, then
the time evolved measure is Gibbs for all times $t$. \item If
$\nu$ is a low temperature Gibbs measure for some t.i. potential
whereas $\mu$ is a high temperature Gibbs measure,
 then the time evolved measure is Gibbs for large $t$. When $\nu$ is not a zero temperature Gibbs measure and $\mu$
 corresponds to a high temperature with a small magnetic filed, Gibbsianness is recovered for larger times.
\end{enumerate}
\end{theorem}

For the sake of simplicity, we describe the results for infinite
temperature Glauber dynamics of low temperature phases of the Ising model
at dimension $d \geq 2$. Starting from the +-phase $\mu_\beta^+$ of the Ising model
at low enough temperature $\beta^{-1} > 0$, we apply a stochastic spin-flip dynamics at rate 1,
 independently over the sites. The time evolved measure is then

\begin{equation}\label{timeevolv}
\mu_{\beta, t}(\eta):= \sum_{\sigma \in \Omega}
\mu_\beta^+(\sigma) \prod_{i \in \mathbb{Z}^d} \frac{e^{\eta_i
\sigma_i h_t}}{2 \cosh h_t}, \; \; \rm{with} \; \;h_t=\frac{1}{2}
\log \frac{1+e^{-2t}}{1-e^{-2t}}.
\end{equation}
The product  kernel in (\ref{timeevolv}) is a special case of a
Glauber dynamics for infinite temperature \cite{EFHR,LNR}, its particular form in terms of a dynamical magnetic
field $h_t$ being obtained by a tricky use of the small size of $E=\{-1,+1\}$. This last particular form allows
to interpret these dynamics as a Kadanoff-like transformation. It is
known that the time-evolved measure $\mu_{\beta,t}$ tends to a
spin-flip invariant product measure on $\{-1,+1\}^{\mathbb{Z}^d}$,
with $t \uparrow \infty$, which is trivially quasilocal and Gibbs.
Nevertheless, the Gibbs property is lost during this evolution and
recovered only  at equilibrium:

\begin{theorem}\cite{EFHR}\label{sto}
Assume that the initial temperature $\beta^{-1}$ is smaller than the critical temperature of the n.n. Ising model
for $d \geq 2$. Then there exists $t_0(\beta) \leq t_1(\beta)$
such that:
\begin{enumerate}
\item $\mu_{\beta,t}$ is a Gibbs measure for all $0 \leq t
<t_0(\beta)$. \item $\mu_{\beta,t}$ is {\bf not} a  Gibbs measure
for all $0 < t_1(\beta)\leq t < + \infty$.
\end{enumerate}
\end{theorem}
Non-Gibbsianness is here related to the possibility of a phase
transition in some constrained model\footnote{The constrained model
is a three dimensional Random Field Ising Model, due to the
randomness of the dynamical field. The possible occurrence of phase
transitions for this model has been proved in \cite{BK}.}. There
remains a large interval of time where the validity of the Gibbs
property of the time-evolved measure remains unknown for this
lattice model \cite{FdH}. This has motivated the study of similar
phenomena for mean-field models in \cite{KLN}, where the sharpness
of the Gibbs/non-Gibbs transition has been proved. This study has
required the introduction of the new notion of Gibbsianness for
mean-field models, see e.g. \cite{HK,K3}, and the relationships
between lattice and mean-field results encourages us to investigate
it further on. Moreover, Gibbianness for short times has been
established for more general local stochastic evolutions in
\cite{LNR} and the large deviation properties has even been proved
to be conserved during the evolution for the special case of the
Glauber dynamics in \cite{LNR2}. Other sources of non-Gibbsianness
during stochastic evolutions, but concerning the stationary and not
the transient regime, have been
 investigated in \cite{Lebsch,Maes4,FT} for e.g. discrete dynamics or probabilistic cellular
 automata. Similar considerations have led to a short review about
 the relationships between non-gibbsianness and disordered systems
 in \cite{VEK}.

\subsection{Joint measure of short range disordered systems \cite{K1}}

This example has allowed substantial progress in the Dobrushin
program of restoration of Gibbsianness, and reinforces our
philosophy of focusing on continuity properties of conditional
probabilities rather than on convergence properties of a
potential, as we shall see next section. Non-Gibbsianness has here
also been very useful to explain pathologies in
 the so-called {\em Morita approach to disordered systems}, see \cite{kuhn,K4} in the proceedings volume \cite{ELNR}.
 The {\em Random Field Ising Model} (RFIM) is an Ising model where the magnetic field $h$ is replaced
 by (say i.i.d. $\pm 1$) random variable $\eta_i$ of common law $\mathbb{P}$ at each site of the lattice.
 For a given $\eta=(\eta_i)_{i \in S}$, whose law is also denoted by $\mathbb{P}$, the corresponding (''quenched'')
 Gibbs measures depend on this disorder and are denoted by $\mu[\eta]$. The Morita approach \cite{kuhn} considers the joint
  measure ''configuration-disorder'', formally defined by $K(d \eta, d \sigma)=\mu[\eta](d \sigma) \mathbb{P}(d \eta)$,
  to be Gibbs for a potential of the joint variables, but it has been proved in \cite{K1} that this measure can be
   non-Gibbs for $d \geq 2$ and for a small disorder. The mechanism, although more complicated, is similar to the
   previous examples, and the arising of non-quasilocality is made possible when a ferromagnetic ordering is
   itself possible in the quenched system, and thus the conditions on $d$ and on the disorder are those required
    for such a phase transition to hold in \cite{BK}. A diluted and simpler version of this phenomenon concerns the
     {\em GriSing} random filed, whose link with Griffiths's singularities is also very relevant for its similarities
     with RG pathologies \cite{VEMSS}.

\subsection{Other sources of Non-Gibbsianness}

Soon after the detection of the renormalization group pathologies as
the manifestation of the occurrence of non-Gibbsianness, the latter
phenomenon has been detected in many other areas of probability
 theory and statistical mechanics, like Hidden Markov models, Random-cluster model,
 convex combinations of product measures, etc. see \cite{VEFS} or references in \cite{ELNR}.\\

Before using these examples to emphasize how important are continuity properties of conditional probabilities
 in the Gibbs formalism, we describe recent extensions of the Gibbs property within the so-called
 {\em Dobrushin program of restoration of Gibbsianness}.

\section{Generalized Gibbs measures}
In 1995, in view of the RG pathologies described in \cite{VEFS}
and in the physics literature, Dobrushin launched a program of
restoration of Gibbsianness consisting in two parts \cite{DOB1}:
\begin{enumerate}
\item To give an alternative (weaker) definition of Gibbsianness
that would be stable under scaling transformations. \item To
restore the thermodynamic properties of these measures in order to
get a proper definition of equilibrium states.
\end{enumerate}
\subsection{Dobrushin Program of restoration of Gibbsianness, Part I}

The first part of this program mainly yields two different
restoration notions, focusing either on a relaxation of the
convergence properties of the potential, leading to {\em weak
Gibbsianness}, or
 on a relaxation on the topological properties of conditional probabilities, leading to {\em almost Gibbsianness} or
  {\em almost quasilocality}. The first one appeared to be weaker than the latter and to be reminiscent to the notion
   chosen to describe systems with hard-core exclusion, or unbounded spins \cite{Lebo}. It express consistency w.r.t.
    an almost surely convergent potential:

\begin{definition}[Weakly Gibbs] A probability measure $\mu \in \mathcal{M}_1(\Omega)$ is said to be {\em weakly Gibbs}
if there exists a potential $\Phi$ and a tail-measurable set $\Omega_\Phi$ on which  $\Phi$ is convergent with
 $\mu(\Omega_\Phi)=1$ such that $\mu \in \mathcal{G}(\gamma^{\Phi})$.
\end{definition}
Tail-measurability is required to insure that the partition function
is well-defined. Weak Gibbsianness has been proved  for most of the
renormalized measures of the previous section
\cite{BKL,LM,MRM,Maes3} and relies on the existence of the already
mentioned relative energies \cite{Maes1}, that usually enables to
prove the almost sure convergence of a telescoping potential of
Kozlov type. A non-Gibbsian measure arising in stochastic evolutions
has also been proved to be weakly Gibbsian \cite{Maes4} and joint
measures of disordered systems too \cite{K2}, and it is actually not
common  in our context to find a transformation of a Gibbs measure
that is not weakly Gibbs, although examples such as convex
combinations of product measures exist \cite{MRM,LMV}. Moreover, in
such a case, the almost sure convergence of the potential does not
tell much about the crucial continuity properties of conditional
probabilities. This has motivated the second main restoration
notion:

\begin{definition}[Almost Gibbs] A probability
measures $\mu$ is {\em almost Gibbs} if its finite-volume
conditional probabilities are continuous functions of the boundary
conditions, except on a set of $\mu$-measure zero, i.e. if there
exists a specification $\gamma$ such that $\mu \in
\mathcal{G}(\gamma)$ and $\mu(\Omega_\gamma)=1$.
\end{definition}

An almost-sure version of  Kozlov-Sullivan's  use of the
inclusion-exclusion principle in the Gibbs representation theorem
proves that {\tt almost Gibbs implies weak Gibbs} \cite{Ko,Su,MRM},
but the converse is not true (see e.g. \cite{lef2}). The decimated
measure has been proved to be almost Gibbsian in \cite{FLNR}, and
the method applies to other renormalized measures \cite{FLNR2}, but
not for the projection to a layer for which the problem is open.
More interestingly, the contrary has been proved \cite{K2} for the
joint measure of the RFIM, which even has
 a set of bad configurations with full measure. The peculiarity of this example appeared to be a good advert
  for the importance of quasilocality in the characterization of equilibrium states for lattice spin systems,
  and to discriminate the weak Gibbs restoration from  almost Gibbs one, due to the consequences it has on the
   thermodynamic properties of the corresponding measures, in the second part of the Dobrushin program.\\

Other restoration notions (robust Gibbsianness \cite{VEL}, fractal
quasilocality \cite{LN}) exist and an intermediate between the
almost Gibbs and weak Gibbs has been introduced recently and seem
to be a very relevant starting definition of generalized Gibbs
measure, called {\em Intuitively weak Gibbsianness} in \cite{VEV}:

\begin{definition}[Intuitively weak Gibbs] A weakly Gibbsian measure $\mu \in \mathcal{M}_1^+(\Omega)$
is said to be {\em intuitively weak Gibbs} if there exists a set $\Omega_1 \subset \Omega_\Phi$ with
$\mu(\Omega_1)=1$ and s.t. for all $\sigma \in \Omega$,
\be \label{IWG}
\gamma_\Lambda^\Phi(\sigma_\Lambda | \omega_{\Lambda_n \setminus \Lambda} \eta_{\Lambda^c}) \;\mathop{\longrightarrow}\limits_{\Lambda \uparrow \mathcal{S}}\;\gamma_\Lambda^\Phi(\sigma | \omega)
\ee
for all $\omega,\eta \in \Omega_1$.
\end{definition}
The notion is intermediate between weak and almost Gibbsianness, the difference between all
these notions being that  the convergence (\ref{IWG}) holds \cite{VEV}:
\begin{itemize}
\item For all $\omega$ and all $\eta$ when $\mu$ is Gibbs (quasilocal).
\item For $\mu$-a.e $\omega$ and all $\eta$ when $\mu$ is almost Gibbs.
\item For $\mu$-a.e $\omega$ and $\mu$-a.e $\eta$ when $\mu$ is intuitively weak Gibbs.
\end{itemize}

The thermodynamics properties of Gibbs measures have been partially restored
for almost Gibbsian measures, whereas weak Gibbsianness seems to be indeed too weak
 to be a satisfactory notion from this point of view.

\subsection{Dobrushin program of restoration, part II}
This second part  aims at the restoration of the thermodynamic
properties of Gibbs measures, mostly in terms of a variational
principle that allows to identify them as the description of the
states which minimize the free energy of the system, thus
equilibrium states in virtue of the second law of thermodynamics.
First, one would also like to recover well-defined thermodynamic
functions at infinite-volume. In the weakly Gibbsian context, this
can be done directly but their existence has to be restricted to
typical boundary conditions, see \cite{lef1,Maes3,Maes1}, whereas
in the case of renormalized measures, the useful notion of {\em
asymptotically decoupled measures} introduced by Pfister \cite{P}
insures their existence and the validity of a large deviation
principle:

\begin{definition}
A measure $\mu \in \mathcal{M}_{1,\rm{inv}}^+(\Omega)$ is {\em asymptotically decoupled}
if there exists functions $g:\mathbb{N} \longrightarrow \mathbb{N}$ and
 $c: \mathbb{N} \longrightarrow [0,\infty)$ s.t.
$$
\lim_{n \to \infty} \frac{g(n)}{n}=0 \; \; \; {\rm and} \; \; \; \lim_{n \to \infty} \frac{c(n)}{|\Lambda_n|}=0
$$
and for all $n \in \mathbb{N},\; A \in \mathcal{F}_{\Lambda_n}, \; B \in \mathcal{F}_{\Lambda_{n+g(n)}^c}$
$$
e^{-c(n)} \; \mu(A) \mu(B) \; \leq \; \mu(A \cap B) \; \leq e^{c(n)} \; \mu(A) \mu(B).
$$
\end{definition}
This class strictly contains the set of all Gibbs measures and with our definition, this property is
 conserved under the local RG transformations introduced in these lectures. Thus thermodynamic
  functions exist for our RGT, but we emphasize the fact that the existence of relative entropy
   is still an open problem for the projection of the Ising model to a layer \cite{FLNR,FLNR2,Maes3}.\\

Thus, our renormalized measures satisfy the thermodynamical variational principle whereas the
specification independent one has been fully restored for the decimated measure in \cite{FLNR},
 using ideas taken from the concept of {\em global specification} \cite{FP}, and as a corollary this
 has established the almost Gibbsianness of this measure. Keeping the same notation and
  writing $\nu^+$ and $\nu^-$ the decimation of the $+$ and $-$ phases of the 2d Ising model, one has:
\begin{theorem}\label{FLNR}\cite{FLNR}
Consider the decimation of the Ising model at   $\beta > \beta_c$.
Then
\begin{enumerate}
\item For every $\mu \in \mathcal{M}_{1,\rm{inv}}^+(\Omega), \; h(\mu|\nu^+)$  exists and $h(\nu^-|\nu^+)=0$.
\item $\nu^- \in \mathcal{G}_{\rm inv}(\gamma^+)$,  where $\gamma^+$ specifies $\nu^+$, and they are almost Gibbs.
\item If $h(\mu | \nu^+)=0$ and $\mu(\Omega_{\gamma^+})=1$ then $\mu \in \mathcal{M}_{1,\rm{inv}}^+(\Omega)$.
\end{enumerate}
\end{theorem}
The proof relies on the general criterion given in Chapter 4, and a partial converse
 statement, to restore the second part of the variational principle, has been established in \cite{KLNR},
  using again the general criterion of Chapter 4.\\

It is moreover  established there that the joint measure of the RFIM, one of the  physically relevant
 known examples of weakly Gibbsian measure that is not almost Gibbs, also  provides an example of two
 (not almost Gibbs) measures described by a different system of conditional probabilities that are equilibrium
 states w.r.t each other. In this situation, we have two candidates to represent equilibrium states, corresponding to different
  interactions, that saturate the variational principle of each other, which can be easily seen to be physically
   irrelevant. The fact that this happens for a weakly and non almost Gibbsian
measure clearly indicates that one has to insist on continuity
properties of the conditional probabilities in order to restore
the Gibbs property in the framework of the Dobrushin program.
Together with the Gibbs representation theorem, this also emphasizes the
relevance of the description of  Gibbs measures in the quasilocal
framework, i.e.  in terms of topological properties of conditional
probabilities rather than in terms of potentials. These
observations have motivated a new ''DLR-like'' approach to
mean-field models, initiated recently in \cite{HK,K3,KLN}.

{\bf Acknowledgements:}\\
 The author thanks the members of the
departments of mathematics of the universities UFMG and UFGRS who
organized and attended to these lectures, A.C.D. van Enter
(Groningen) for its legendary careful readings and suggestions and
S. Friedli (Belo Horizonte) for many comments, readings and
discussions.


\addcontentsline{toc}{section}{\bf
References}
\begin{thebibliography}{14}

\bibitem{Aiz}  M. Aizenman. Translation-invariance and instability of phase coexistence in the two-dimensional
 Ising system. {\em Comm. Math. Phys.} {\bf 73}, no 1:83--94, 1980.

\bibitem{Bei} H. van Beijeren. Interface sharpness in the Ising model. {\em Comm. Math. Phys.} {\bf 40}1--6, 1975.


\bibitem{BR}  P. Berti, P. Rigo. 0-1 laws for regular conditional distributions. {\em Ann. Prob.} {\bf 35}:649--662, 2007.

\bibitem{BR2} P. Berti, P. Rigo. A conditional 0-1 law for the symmetric $\sigma$-field. To appear to {\em J. Theo. Proba.}, 2007.

\bibitem{BCO} L. Bertini, E. Cirillo, E. Olivieri. Renormalization group in the uniqueness region:
weak Gibbsianity and convergence.  {\em Comm. Math. Phys.}  {\bf 261}, no. 2:323--378, 2006.

\bibitem{Bil2} P. Billingsley. {\em Ergodic theory and information}. John Wiley, New York-London-Sydney, 1965.

\bibitem{Bill} P. Billingsley. {\em Probability and measures}. 3rd ed. Wiley series in probability and mathematical statistics. New York : Wiley, 1968 and 1995.





\bibitem{BLG}
P.M. Bleher, N.N. Ganihodgaev. \newblock On pure phases of the Ising
model on the Bethe lattice. \newblock \emph{Theo. Prob. Appl.} {\bf
  35}:1-26, 1990.


\bibitem{Bod} T. Bodineau. Translation-invariant Gibbs states for the Ising model.  {\em Prob. Th. Relat. Field.}  {\bf 135},  no. 2, 153--168, 2006.

\bibitem{Bo}
L. Boltzmann. \newblock Le\c cons sur la th\'eorie des gaz,
Gauthier-Villars, Paris, tome I 1902, tome II 1905. R\'e\'edition
Jean Gabay, Paris, 1987.


\bibitem{Bo2}
L. Boltzmann. \newblock The second law of Thermodynamics. In {\em Ludwig Boltzmann, Theoretical Physics and
Philosophical Problems}, B. Mc Guinness eds, Reidel, 1974.

\bibitem{BZ} A. Bovier, M. Zahradn\'ik. The low temperature phase of Kac-Ising models. {\em  J. Stat.Phys.}
{\bf 87}, nos 1/2:311--332, 1997.

\bibitem{BK} J. Bricmont, A. Kupiainen. Phase transition in the 3d Random Field Ising Model. {\em Comm. Math. Phys.} {\bf 142}:539--572, 1988.

\bibitem{BKL} J. Bricmont, A. Kupiainen, R. Lefevere. Renormalization group pathologies and the definition of Gibbs states.
{\em Comm. Math. Phys.} {\bf 194}, no. 2:359--388, 1998.

\bibitem{BFS} R.  Burton, C.-E. Pfister, J. Steif.
The Variational Principle for Gibbs States Fails on Trees. {\em Mark. Proc. Relat. Fields} {\bf  1}, no. 3:387--406, 1995.

\bibitem{Sido} F. Camia, C.M. Newman,  V. Sidoravicius. A particular bit of universality: scaling limits of some dependent percolation models.  {\em Comm. Math. Phys.} {\bf 246}, no. 2:311--332, 2004.

\bibitem{RGT2} J. Cardy. {\em Scaling and renormalization in statistical physics}. Cambridge Lecture Notes in Physics {\bf 5}. Cambridge University Press, Cambridge, 1996.

\bibitem{CP} M. Cassandro, E. Presutti. Phase transitions in Ising systems with long but finite range interactions. {\em Mark. Proc. Relat. Fields} {\bf 2}, no 2:241--262, 1996.

\bibitem{Com} F. Comets. Large Deviation
Estimates for a Conditional Probability Distribution. Applications
to Random Interaction Gibbs Measures. {\em Prob. Th. Relat.
Fields.} {\bf 80}:407-432, 1986.



\bibitem{C} I. Csisz\'ar. \newblock Information-type measures of
difference of probability distributions and indirect observations.
\newblock \emph{Studia Sci. Math. Hungar.} {\bf 2}:299--318, 1967.


\bibitem{CH} K.L. Chung. Markov chains with stationary transition probabilities. Second edition, Springer-Verlag, New-York, 1987.

\bibitem{AMasi} A. De Masi. Systems with long range interactions. {\em Progress in Probability} {\bf 54}:25--81, Birkhauser, 2003.

\bibitem{Maes4} J. Depoorter, C. Maes. Stavskaya's measure is weakly Gibbsian.  {\em Mark. Proc. Relat. Fields}  {\bf 12}, no 4: 791--804, 2006.



\bibitem{DG} C. Domb, M.S. Green (Eds). {\em Phase transitions and Critical Phenomena}, Vol. 6, Academic Press, NY, 1976.

\bibitem{Diac} P. Diaconis. Recent progress on de Finetti`s notions of exchangeability. {\em Bayesian statistics} {\bf 13}:111--125, Oxford University Press, 1988.


\bibitem{DOB0} R.L. Dobrushin. Existence of a phase transition in the two-dimensional and three-dimensional Ising models.  {\em Dokl. Akad. Nauk SSSR}  {\bf 160}:1046--1048 (Russian); translated as {\em Soviet Physics Dokl.}{\bf 10}:111--113, 1965.


\bibitem{DOB} R.L. Dobrushin. Description of a random field by means of conditional probabilities and the conditions governing its regularity. {\em Theor. Prob. Appl.} {\bf 13}:197--224, 1968.

\bibitem{DOB01} R.L. Dobrushin. Gibbs states describing coexistence of phases for a three-dimensional Ising model. {\em Theo. Prob. Appl.} {\bf 17}, no 4:582--600, 1972.

\bibitem{DOBPe} R.L. Dobrushin, E.A. Pecherski. A criterion for the uniqueness of Gibbsian fields in the non-compact case. In {\em Probability theory and Mathematical Statistics}, Lecture Notes in Mathematics {\bf 1021}:97--110, Springer-Verlag, Berlin, 1983.


\bibitem{DOB1} R.L. Dobrushin, S.B. Shlosman. \newblock Gibbsian description of ''non Gibbsian'' field, \newblock {\em Russian Math Surveys} {\bf 52}, 285-297, 1997. Also ''Non-Gibbsian'' states and their description. {\em Comm. Math. Phys.}{\bf 200},no 1: 125-179, 1999.

\bibitem{Dy} E.B. Dynkin. Sufficient statistics and extreme points. {\em Ann. Proba.} {\bf 6}, No. 5:705--730, 1978.

\bibitem{Dys} F.J. Dyson. An Ising ferromagnet with discontinuous long-range order.  {\em Comm. Math. Phys.} {\em 21}:269--283, 1971.


\bibitem{VE} A.C.D. van Enter. On the possible failure of the Gibbs property for measures on lattice systems. Disordered systems and statistical physics: rigorous results (Budapest, 1995).  {\em Mark. Proc. Relat. Field.}  {\bf 2},  no. 1:209--224, 1996.


\bibitem{VEF} A.C.D. van Enter, R. Fern{\'a}ndez. A remark on different norms and analycity for many particle interactions. {\em J. Stat. Phys.} {\bf 56}: 965--972, 1989.

\bibitem{EFHR} A.C.D. van Enter, R. Fern{\'a}ndez, F. den Hollander, F. Redig. Possible loss and recovery of Gibbsianness during the stochastic evolution of Gibbs measures.  {\em Comm. Math. Phys.} {\bf 226}, no. 1, 101--130, 2002.

\bibitem{VEFS}
A.C.D. van Enter, R.~Fern{\'a}ndez, A.D. Sokal. \newblock
Regularity properties and pathologies of position-space
renormalization-group transformations: Scope and limitations of
{G}ibbsian theory. \newblock {\em J. Stat. Phys.} {\bf
72}:879-1167, 1993.

\bibitem{VEK}  A.C.D. van Enter, C.  K\"ulske. Two connections between random systems and non-Gibbsian measures.
{\em J. Stat. Phys.} {\bf 126}, no. 4-5:1007--1024, 2007.

\bibitem{ELNR} A.C.D. van Enter, A. Le Ny, F. Redig (eds). Proceedings of the workshop
"Gibbs vs. non-Gibbs in statistical mechanics and related fields"
(Eurandom 2003). {\em Mark. Proc. Relat. Fields} {\bf 10}, no 3,
2004.

\bibitem{VEL} A.C.D. van Enter, J. L\"orinczi. Robustness of the non-Gibbsian property: some examples.  {\em J. Phys. A}  {\bf 29}, no. 10:2465--2473, 1996.
\bibitem{VEMSS} A.C.D. van Enter, C. Maes, R.H.  Schonmann, S.B. Shlosman. The Griffiths singularity random field.  {\em in} On Dobrushin's way. From probability theory to statistical physics, {\em Amer. Math. Soc. Transl.} Ser. 2, {\bf 198}: 51--58, Amer. Math. Soc., Providence, RI, 2000.

\bibitem{VEV} A.C.D. van Enter, E.A. Verbitskiy.  On the variational principle for generalized Gibbs measures.  {\em Mark. Proc. Relat. Fields} {\bf 10}, no 3:411--434, 2004.
\bibitem{fell} W. Feller. \newblock An introduction to probability theory and its applications, volume I. \newblock Wiley publications in statistics, 1957.

\bibitem{F} R. Fern\'andez. Gibbsianness and non-Gibbsianness in Lattice random fields.
In {\em Mathematical Statistical Physics. Proceedings of the
$83$rd Les Houches Summer School (july 2005)}, Elsevier, A.
Bovier, A.C.D. van Enter, F. den Hollander, F. dunlop eds., 2006.

\bibitem{FP} R. Fern\'andez, C.-E. Pfister. Global specifications
and non-quasilocality of projections of Gibbs measures. {\em Ann.
Proba.} {\bf 25}, no 3:1284-315, 1997.

\bibitem{FLNR} R. Fern\'andez, A. Le Ny, F. Redig.
Variational principle and almost quasilocality for renormalized
measures. {\em J. Stat. Phys.} {\bf 111}, nos 1/2:465--478, 2003.

\bibitem{FLNR2}  R. Fern\'andez, A. Le Ny, F. Redig.
Restoration of Gibbsianness for projected and FKG renormalized
measures. {\em Bull. Braz. Math. Soc.} {\bf 34}:437-55, 2003.

\bibitem{FT} R. Fern\'andez, A. Toom.  Non-Gibbsianness of the invariant measures of non-reversible cellular automata with totally asymmetric noise. {\em Geometric methods in dynamics.} II.  Ast本isque  No. 287, 2003.

\bibitem{fis} M.E. Fisher. \newblock Scaling, universality and renormalization group
theory, \newblock in {\em Critical phenomena (Stellenbosch 1982)},
Lecture Notes in Physics no 186:1-139, F.J.Hahne ed., Springer-Verlag,
Berlin, 1983.

\bibitem{FS} J. Fr\"ohlich, T. Spencer. The phase transition in the one-dimensional Ising model with $1/r^2$ interaction energy.  {\em Comm. Math. Phys.}  {\bf 84}, no. 1:87--101, 1982
.
\bibitem{RGT} K. Gawedski. Rigorous renormalization group at work. {\em Physica} {\bf 140A}: 78--84, 1986.
\bibitem{Gi} J.W. Gibbs. {\em Elementary principles in statistical mechanics}.
Yale University Press, New Haven, 1902.

\bibitem{Ge} H.O. Georgii. {\em Gibbs Measures and Phase Transitions}.
 Walter de
Gruyter (de Gruyter Studies in Mathematics, Vol.\ 9), Berlin--New
York, 1988.

\bibitem{GeHi} H.O. Georgii, Y. Higuchi. Percolation and number of
phases in the two-dimensional Ising model. Probabilistic
techniques in equilibrium and non-equilibrium statistical physics.
{\em J. Math. Phys.} {\bf 41}, no 3:1153--1169, 2000.

\bibitem{Grif} R.B. Griffiths. Peierls proof of spontaneous magnetization in a two-dimensional Ising ferromagnet.  {\em Phys. Rev.} {\bf 2}, 136:A437--A439, 1964.

\bibitem{GP} R.B. Griffiths, P.A. Pearce. \newblock  Mathematical properties of
position-space renormalization-group transformations,
\emph{J. Stat. Phys.} {\bf 20}: 499-545, 1979.


\bibitem{hig}
Y. Higuchi. \newblock Remarks on the limiting Gibbs states on a
(d+1)-tree. \newblock \emph{Publ. RIMS, Kyoto Univ.} {\bf
13}:335-348, 1977.

\bibitem{Hig} Y. Higuchi.  On the absence of non-translation-invariant Gibbs states for the two-dimensional Ising model.
 {\em Random fields, Vol I, II (Esztergom, 1979)}, Colloq. Math. Soc. Janos Bolyai {\bf 27}:517--534, 1981.


\bibitem{HK} O. H\"aggstr\"om, C. K\"ulske. Gibbs property of the fuzzy Potts model on trees and in mean-field. {\em Mark. Proc. Relat. Fields} {\bf 10} no  3:477--506, 2004.

\bibitem{HKen} K. Haller, T. Kennedy. Absence of renormalization group pathologies near the critical temperature. Two examples.  {\em J. Stat. Phys.}  {\bf 85}, no. 5-6:607--637, 1996.

\bibitem{FdH} F. den Hollander. Gibbs under stochastic dynamics ? {\em Mark. Proc. Relat. Fields} {\bf 10}, no  3:507--516, 2004.

\bibitem{jof} D. Ioffe. \newblock Extremality of the disordered state for the Ising
model on general trees.  {\em Progress in probability} {\bf 40}:3-14, B. Chauvin, S. Cohen, A. Rouault eds, 1996.


\bibitem{Is2} R.B. Israel. \newblock Banach algebras and Kadanoff transformations, in
\emph{Random Fields (Esztergom, 1979)} J. Fritz, J.L. Lebowitz, and
D.Sz\'asz eds, Vol II, pp.593-608, 1981.


\bibitem{Is} R.B. Israel. {\em Convexity in the theory of
lattice gases}. Princeton university Press, 1986.

\bibitem{Is3} R.B. Israel. Some generic results in Mathematical
Physics. {\em Mark. Proc. Rel. Field.} {\bf 10}, no 3:517--523,
2004.

\bibitem{Ja} E.T. Jaynes. {\em E. T. Jaynes: papers on probability, statistics and statistical physics}.
 Reprint of the 1983 original. Edited and  introduced by R. D. Rosenkrantz. Pallas Paperbacks, 50. Kluwer Academic Publishers Group, Dordrecht, 1989.

\bibitem{Kac} M. Kac, G. Uhlenbeck, P. Hemer. On the van der waals theory of vapor-liquid equilibrium. {\em J. Math. Phys} {\bf 4}:216--228, 1963.



\bibitem{Ken} T. Kennedy. Majority rule at low temperatures on the square and triangular lattices.  {\em J. Stat. Phys.} {\bf 86}, no. 5-6, 1089--1107, 1997.

\bibitem{Kh} A.I. Khinchin. Mathematical foundations of information theory. New York : Dover, 1957.

\bibitem{KW} H.A. Kramers, G.H. Wannier. Statistics of the
two-dimensional ferromagnet I, II. {\em Phys. Rev.} 2, no {\bf
60}:252--262 and 263--276, 1941.

\bibitem{KO} S. Kobe. Ernst Ising, 1900-1998. {\em Braz. J. Phys.} {\bf 30}, no.4: 649-654, 2000.

\bibitem{Ko} O. Kozlov. Gibbs description of a system of random variables. {\em Problems Inform. Transmission.} {\bf 10}:258--265, 1974.


\bibitem{K1} C. K\"ulske. (Non-) Gibbsianness and phase transition in random lattice spin models.
{\em Mark. Proc. Rel. Field.} {\bf 5}:357--383, 1999.

\bibitem{K2} C. K\"ulske. Weakly Gibbsian representation for joint measures of quenched lattice spin models. {\em Prob. Th. Relat. Fields} {\bf 119}:1--30, 2001.

\bibitem{K3} C. K\"ulske. Analogues of non-Gibbsianness in joint-measures of disordered mean-field models. {\em J. Stat. Phys.} {\bf 112}, no 5/6:1079--1108, 2003.

\bibitem{K4} C. K\"ulske. How non-Gibbsianness helps a metastable Morita minimizer to provide a stable free energy.  {\em Mark. Proc. Relat. Field.}  {\bf 10}, no. 3: 547--564, 2004.

\bibitem{KLN} C. K\"ulske, A. Le Ny. Spin-flip dynamics
of the  Curie-Weiss model: Loss of Gibbsianness with possibly
broken symmetry.  {\em Comm. Math. Phys.} {\bf 271}, vol 2:431--454, 2007.

\bibitem{KLNR} C. K\"ulske, A. Le Ny, F. Redig.
Relative entropy and variational properties of generalized
Gibbsian measures. {\em Ann. Proba.} {\bf 32}, no. 2:1691--1726,
2004.

\bibitem{kuhn} R. K\"uhn. Gibbs vs Non-Gibbs in the Equilibrium
Ensemble approach to disordered systems. {\em Mark. Proc. Relat.
Field.}  {\bf 10}, no. 3: 523--546, 2004.
\bibitem{Lan}
O.E. Lanford. Entropy and Equilibrium States in Classical Statistical Mechanics. {\em In} Statistical Mechanics and Mathematical Problems, Battelle Seattle 1971 Rencontres, Lectures Notes in Physics no 20, Springer-Verlag, Berlin etc., 1973.


\bibitem{LR} O.E. Lanford, D. Ruelle. Observables at infinity and states with short range correlations in statistical mechanics.  {\em Comm. Math. Phys.} {\bf  13}:194--215, 1969.



\bibitem{LN} A. Le Ny. Fractal failure of quasilocality for a majority rule transformation on a tree.  {\em Lett. Math. Phys.}  {\bf 54}, no. 1:11--24, 2000.
\bibitem{LNR} A. Le Ny, F. Redig. Short times conservation of Gibbsianness under local stochastic evolutions. {\em J. Stat. Phys.} {\bf 109}, nos 5/6:1073--1090, 2002.

\bibitem{LNR2} A. Le Ny, F. Redig. Large deviation principle at fixed time in Glauber evolutions.  {\em Mark. Proc. Relat. Fields} {\bf  10}, no 1:65--74, 2004.


\bibitem{Lebo} J. Lebowitz. Statistical mechanics of systems of unbounded spins.  {\em Comm. Math. Phys.} {\bf  50}:195-218, 1976.

\bibitem{LebPe} J. Lebowitz, O.E Penrose. On the exponential decay of correlations. {\em Comm. Math. Phys.} {\bf  39}:165--184, 1974.

\bibitem{LebMaes} J. Lebowitz, C. Maes. Entropy: a dialogue.  {\em Entropy}: 269--276, Princeton Ser. Appl. Math., Princeton Univ. Press, Princeton, NJ, 2003.

\bibitem{Lebsch} J. Lebowitz, R.H. Schonmann. Pseudo-free energies and Large deviations for Non Gibbsian FKG measures. {\em Prob. Th. Relat. Field.} {\bf 77}:49--64, 1988.

\bibitem{lef1} R. Lefevere. Variational principle for some renormalized measures.  {\em J. Stat. Phys.} {\bf 96}, nos 1-2, 109--133, 1999.

\bibitem{lef2}
R. Lefevere. Weakly Gibbsian measures and quasilocality: a
long range pair-interaction example. {\em J. Stat. Phys.} {\bf 95}, no 3/4:785-789, 1999.

\bibitem{LP} J.T. Lewis, C.-E. Pfister. Thermodynamic probability theory: some aspects of large deviations. {\em Uspekhi Mat. Nauk}  {\bf 50}  no. 2(302):47--88;  translation in  {\em Russian Math. Surveys}  {\bf 50}, no. 2:279--317, 1995.

\bibitem{LPS}
J.T. Lewis, C-E. Pfister, W.G. Sullivan. Entropy, concentration of probability and conditional limit theorems. {\em Mark. Proc. Relat. Fields} {\bf 1}, no 3:319-386, 1995.

\bibitem{ligg} T. M. Liggett. {\em Interacting particle systems}. Grundlehren der Mathematischen Wissenschaften [Fundamental Principles of Mathematical Sciences, 276. Springer-Verlag, New York, 1985.

\bibitem{LM} J. L\"orinczi, C.  Maes. Weakly Gibbsian measures for lattice spin systems.  {\em J. Stat. Phys.}  {\bf 89},  no. 3-4:561--579, 1997.

\bibitem{LMV} J. L\"orinczi, C. Maes, K. Vande Velde. Transformations of Gibbs measures.  {\em Prob. Theo. Relat. Field.}  {\bf 112},  no. 1, 121--147, 1998.



\bibitem{MRM} C. Maes, F. Redig, A. Van Moffaert. Almost Gibbsian versus Weakly Gibbsian.
{\em Stoch. Proc. Appl.} {\bf 79}, no 1:1-15, 1999.

\bibitem{Maes3}
C.~Maes, F.~Redig, A.~Van Moffaert.
The restriction of the Ising model to a Layer.
{\em J. Stat. Phys.} {\bf 96}, nos 1/2:69-107, 1999.


\bibitem{Maes1}
C.~Maes, K.~Vande Velde. Relative energies for non-Gibbsian states.  {\em Comm. Math. Phys.} {\bf 189}, no 2:277-286, 1997.


\bibitem{Dell} P.-A. Meyer. {\em Probabilities and potential}. North Holland mathematics studies, Amsterdam: North-Holland, 1978.
\bibitem{Ons} L. Onsager. Crystal statistics. I. A two-dimensional model with an order-disorder transition.  {\em Phys. Rev.} {\bf 65} no 2:117--149, 1944.




\bibitem{Peie} R. B. Peierls. On Ising`s model of ferromagnetism. {\em Proc. Cambridge Philos. Soc.} {\bf 32}:477--481, 1936.

\bibitem{P} C.-E. Pfister. Thermodynamical aspects of classical lattice systems. In {\em In and out of equilibrium.
Probability with a physical flavor}. Progress in Probability,
Vladas Sidoravicius ed., Birkha\"user, pp 393--472, 2002.

\bibitem{PO}  A. Petri, M. De Oliveira. \textit{Temperature of non-equilibrium lattice systems}.
Intern. J. Mod. Phys. C \textbf{17}, no 12:1703--1715, 2006.



\bibitem{PS}
S.A. Pirogov, Y. Sinai. Phase Diagram of Classical Lattice Systems.
{\em Theor. Math. Phys} {\bf 25}:1185-1192, 1975 and {\bf 26}:39-49,
1976.

\bibitem{Pres2}
C.Preston. {\em Gibbs states on Countable sets}.
Cambridge tracts in Math., no 68, Cambridge University
Press, London, New York, 1974.

\bibitem{Pres}
C.Preston. \newblock {\em Random Fields}. \newblock Lecture Notes in Mathematics 534, \newblock Springer-Verlag, 1976.



\bibitem{Pr} B. Prum. {\em Processus sur un r\'{e}seau et mesures de
Gibbs}. Techniques stochastiques, Masson, Paris, 1986.

\bibitem{Rue} D. Ruelle. Equilibrium states of infinite systems in statistical mechanics.  {\em Mathematical aspects of statistical mechanics (Proc. Sympos. Appl. Math., New York, 1971)}. SIAM-AMS Proceedings, Vol. V:47--53, Amer. Math. Soc., Providence, R. I., 1972.


\bibitem{Ru} L. Russo. The infinite cluster method in the two-dimensional Ising model.  {\em Comm. Math. Phys.}  {\bf 67}, no. 3, 251--266, 1979.

\bibitem{sch}
R.H. Schonmann. Projection of Gibbs measures may be
non-Gibbsian. {\em Comm. Math. Phys.} {\bf 124}:1-7, 1989.

\bibitem{Simon} B. Simon. {\em The statistical mechanics of lattice gases}. Vol. I. Princeton Series in Physics. Princeton University Press, Princeton, NJ, 1993.





\bibitem{spi2}
F. Spitzer. Phase transition in one-dimensional nearest-neighbors
systems, {\em J. Funct. Anal.} {\bf 20}:240-255, 1975.

\bibitem{Stoy} J. Stoyanov. {\em Counterexamples in probability}. Wiley series in probability and mathematical statistics, Wiley, 1987.

\bibitem{Su2} W.G. Sullivan. \newblock Finite range random fields and energy fields. {\em J. Math. Anal. Appl.} {\bf 44}:710-724, 1973.


\bibitem{Su} W.G. Sullivan. \newblock Potentials for almost Markovian random fields,
  \newblock {\em Comm. Math. Phys.} {\bf 33}:61-74, 1976.

\bibitem{Will} D. Williams. {\em Probability with martingales}. Cambridge mathematical textbooks, Cambridge University Press, 1991.

\bibitem{wil}
K.G. Wilson. \newblock The renormalization group: Critical phenomena
and the Kondo problem, \newblock {\em Rev. Mod. Phys.} {\bf
  47}:773-840, 1975.

\bibitem{Xu} S. Xu. An ergodic process of zero divergence distance from the class of all stationary processes. {\em J. Theo. Prob.} {\bf 11}, no 1:181--195, 1988.
\bibitem{Ya} C.N. Yang. The spontaneous magnetization of a two-dimensional Ising model.  {\em Phys. Rev.} {\bf 85} no 2: 808--816, 1952.











\end{thebibliography}
\end{document}